\documentclass{amsart}
\usepackage{amsfonts}
\usepackage{amssymb}
\usepackage[pctex32]{graphics}
\usepackage[all,poly,graph]{xy}
\font\gotb=eufm10 at 14pt
\font\got=eufm10
\font\gots=eufm10 at 7pt
\def\g2{ \hbox{\got g}_2}
\def\f4{\hbox{\got f}_4}
\def\e6{\hbox{\got e}_6}
\def\fs4{\hbox{\gots f}_4}
\def\F4{\hbox{\got F}_4}
\def\Fs4{\hbox{\gots F}_4}

\def\id{\hbox{id}}
\def\h{\hbox{\got h}}

\def\L{{\mathcal
L}}
\def\imag{i}
\def\t{ {\rm tr}}
\def\T{ {\rm Tr}}

\def\N{ {\rm N}}
\def\Q{ {\rm Q}}
\def\W{\mathcal W}
\def\sitem#1{\newline #1}

\def\Tau{\mathcal T}
\def\H{ {  H}}
\def\K{ {  K}}
\def\X{\hbox{\got X}}
\def\Xs{\hbox{\gots X}}

\def\C{\mathbb C}
\def\Z{\mathbb Z}
\def\R{\mathbb R}
\def\aut{\mathop{\rm
aut}}
\def\der{\mathop{\rm Der}}
\def\Sym{ {\rm Sym}}
\def\Ad{ {\rm Ad}}

\def\In{ {\rm In}}
\def\s{ {\rm sl}}

\def\GL{ {\rm GL}}
\def\SO{ {\rm SO}}
\def\gl{ {\rm gl}}
\def\span#1{\langle #1 \rangle}

\def\rank{\mathop{\rm rank}}
\def\fix{\mathop{\rm fix}}
\def\diag{\mathop{\rm diag}}
\def\a{\alpha}
\def\b{\beta}
\def\al{\ifcase\xypolynode\or F \or A\or B\or C\or D\or G\fi}
\def\ala{\ifcase\xypolynode\or a \or b\or c\or d\or g\or f\fi}

\def\sup#1{^{(#1)}}
\def\om{\omega}
\def\si{\sigma}
\def\To{\hbox{\got T}}
\def\So{\hbox{\got S}}
\def\Co{\hbox{\got C}}
\def\No{\hbox{\got N}}
\def\Nos{\hbox{\gots N}}
\def\tor#1{\To^{\langle #1\rangle}}
\def\sor#1{\So^{\langle #1\rangle}}
\newtheorem{te}{Theorem}
\newtheorem{pr}{Proposition}
\newtheorem{lm}{Lemma}
\newtheorem{co}{Corollary}

\title[Gradings on Albert algebra and on $\fs4$]{Gradings on the Albert algebra and on $\hbox{\gotb f}_4$}
\author{Cristina
Draper Fontanals}
\address{Cristina Draper Fontanals: Departamento de
Matem\'atica Aplicada\\ Campus de El Ejido, S/N, 29071 M\'alaga,
Spain.}
\author{C\'andido Mart\'{\i}n Gonz\'alez}
\thanks{Supported by the Spanish MCYT
projects MTM2004-06580-C02-02 and MTM2004-08115-C04-04, and by the
Junta de Andaluc\'{\i}a PAI projects FQM-336 and FQM-1215}
\address{C\'andido Mart\'{\i}n Gonz\'alez:
Departamento de \'Algebra, Geometr\'{\i}a y Topolog\'{\i}a\\
Campus de Teatinos, S/N. Facultad de Ciencias, Ap.\,59, 29080
M\'alaga, Spain.}

\begin{document}
\maketitle

 \begin{abstract}
We study group gradings on the Albert algebra and on the simple
exceptional Lie algebra $\frak{f}_4$ over algebraically closed
fields of characteristic zero. The immediate precedent of this
work is \cite{g2} where we described (up to equivalence) all the
gradings on the exceptional simple Lie algebra $\frak{g}_2 $.   In
the cases of the Albert algebra and $\frak{f}_4$, we look for the
nontoral gradings finding that there are only eight nontoral
nonequivalent gradings on the Albert algebra (three of them being
fine) and nine on $\frak{f}_4$ (also three of them fine).
\end{abstract}

\section{Introduction}
The interest on gradings on simple Lie and Jordan algebras has
been remarkable in the last years. The gradings of finite
dimensional simple Lie algebras, ruling out $\frak a_l$, $\frak
d_4$ and the exceptional cases, are described in \cite{Ivan2}. The
gradings on simple Jordan algebras of type $H_n(F)$ and $H_n(Q)$
are given in the same reference, for an algebraically closed field
$F$ of characteristic zero and a quaternion algebra $Q$. In that
work, the authors use their previous results in \cite{Ivan15}
about gradings of associative algebras $M_n(F)$. In \cite{Ivan1}
all gradings on the simple Jordan algebras of Clifford type have
been described. The fine gradings on $\frak a_l$ have been
determined in \cite{LGII} solving the related problem of finding
maximal abelian groups of diagonalizable automorphisms of the
algebras (not only in $\GL(n,\C)$ but also in ${{\rm O}}(n,\C)$
for $n\ne 8$ and ${{\rm SP}}(2n,\C)$). General notions about Lie
gradings are considered in \cite{Zass}, and the real case is
treated in \cite{LGIII}. Notice that all the mentioned works make
use of techniques related to the associative case. The first
studies of gradings on exceptional Lie algebras are \cite{g2} and
\cite{otrog2}, which describe the group gradings on $\g2$. To
continue the study of exceptional Lie and Jordan algebras, we
aboard in this paper the task of describing nontoral group
gradings on $\f4$ and on
 the Albert algebra. The notion of group grading is
closely related to that of commuting set of semisimple
automorphisms (or equivalently abelian subgroup of semisimple
automorphisms) of the algebra. But  $F_4$, the automorphism group
of the Albert algebra $J$, is isomorphic to the automorphism group
of the Lie algebra $\f4=\der(J)$. So automorphism information can
be transferred from one to the other context. That is the reason
to study both algebras jointly.

We have ruled out the study of toral   gradings on $\f4$ (and then
on $J$) because of the overwhelming proliferation of nonequivalent
cases, and the fact that, their determination (though tedious)
follows from mechanical coarsenings of the Cartan grading. On the
contrary, the nontoral gradings on a simple Lie algebra are not
compatible with the root system. So they could lead us to new ways
of looking at the algebra. This is specially true in the case of
fine gradings, a fact that explains the activity around this
subject (see, for instance, \cite{Phys2} or \cite{Phys3}). In
other papers of Lie gradings (\cite{o4}), the aim is to study the
Lie gradings without the restriction of a grading group,  but we
have preferred not to adopt this approach because of \cite{Alb2}.
It is still an open question the existence of a grading group on
any finite dimensional graded simple Lie algebra over an
algebraically closed field of characteristic zero.

The techniques employed to search the gradings have been of very
different nature. It could be said that this is a
multidisciplinary field nowadays. There is a first tool, very
  intuitive, which is the
usage of models of the algebra. But, although it provides most of
the existing gradings (actually all of them), this fact cannot be
proved without a more powerful tool. This is why we exploit the
benefit of working inside the normalizer of a maximal torus of the
automorphism group. This turns out to be quite technical and less
intuitive (some computer aided arguments have been essential).
However this approach allows to confirm with full precision all
the hypothesis about how many nontoral gradings appear, the
relation among them and other aspects.

A summary of the contents of the work follows in the next
paragraphs.

Section 2 presents some preliminaries, and compiles basic facts on
gradings and related topics. In Subsection 2.1 we introduce some
terminology of algebraic groups needed for our approach. In 2.2
  we devote some attention to recall  the most known  model of
  the Albert algebra, $J=\H_3(C)$ for a Cayley algebra
  $C$, described in \cite{Schafer}.
Each model of a given algebra has its own advantages (and
drawbacks) to present particular gradings, thus
 in 2.4 Tits construction of the Albert algebra is also
recalled. In 2.3   we consider  the notion of toral grading. In
the Lie algebras case, a grading is toral when its homogeneous
components are sum of root spaces. In general a grading is toral
if it is produced by a set of automorphisms contained in a torus
of the automorphism group. At this point we fix a maximal torus of
$F_4$, and characterize the torality of a grading in different
terms.

In Section 3 we start inducing gradings on $J$ from gradings on
the related Cayley algebra $C$ such that $J=\H_3(C)$. This arises
from the well-known possibility of extending automorphisms of $C$
to elements in $F_4$ ($G_2\subset F_4$). It turns out that the
unique (up to equivalence) nontoral grading of $\g2$, induces a
nontoral grading on $\f4$ (a $\Z_2^3$-grading). Thus we are led to
the birth of our first nontoral grading on $J$. This grading will
  induce a numerous family of gradings whose relatives come from a mixing
process described in Subsection 3.3. Previously, we consider all
the gradings of the subalgebra $\H_3(F)$, where $F$ stands for the
ground field. But automorphisms of $\H_3(F)$ can also be nicely
extended to automorphisms of $J$. And curiously all of these
commute with automorphisms of $J$ coming from $C$. Thus by
crossing the nontoral $\Z_2^3$-grading on $J$ coming from $C$ with
all the gradings detected in $\H_3(F)$ we obtain a family of six
nontoral gradings on $J$ described in 3.3. But, as it is pointed
out at the end of such subsection, one of these six nontoral
gradings admits a proper nontoral coarsening. Thus we are led to a
set of seven nonequivalent nontoral gradings on $J$.

Unfortunately we are not done with this set of gradings. To detect
the remaining nontoral gradings on $J$ we need to invoke some
other model of $J$ different from the usual $\H_3(C)$. But at this
point,  Tits construction comes in our help to provide a
$\Z_3^3$-grading with comes from the natural embedding of
automorphisms of the algebra $  M_3(F)$ into $F_4$. The origin of
this $\Z_3^3$-grading is a  known   nontoral $\Z_3^2$-grading of $
M_3(F)$ which can be lifted and finally refined to the nontoral
$\Z_3^3$-grading on $J$. In this way we get a set of eight
pairwise nonequivalent nontoral gradings. Furthermore, any
nontoral grading is equivalent to some of these. This is one of
our main results presented at the end of Section 3 (though the
proof will have to be postponed to the final sections of the
work).

In Section 4 we focus on $\f4$. How can we present a concrete
grading on a Lie algebra? One of the  ways in which a grading on
an algebra can be  given is to provide the set of automorphisms
inducing the grading. Specially because the automorphisms can be
given in terms of toral elements and the Weyl group. This happens
because of the relevant fact that the quasitorus inducing the
grading is always contained in the normalizer of a maximal torus.
Consequently, our first task is to fix a particular representation
of its Weyl group $\W$ and provide a set of representatives of
conjugacy classes in $\W$. Next we present in Subsection 4.1 a
maximal torus of $\F4:=\aut(\f4)$ nicely related to the maximal
torus previously presented in $F_4:=\aut(J)$. Also the action of
$\W$ on the maximal torus is described.

Section 5 is intended to provide the main results on quasitori
which will allow the classification of fine and nontoral gradings.
We introduce a family of quasitori $A(j,t)\subset\F4$ and
immediately proceed with the study of those which are nontoral.
Theorem~\ref{main1} describes the maximal  quasitori of $\F4$ up
to conjugacy which in particular yields the classification of fine
gradings on $\f4$ up to equivalence (Corollary~\ref{fiinas}).
Finally in Subsection 5.2 we get a more detailed classification up
to conjugacy of nontoral quasitori of $\F4$. This result contained
in Theorem~\ref{main2} might be of independent interest. As a
byproduct, we get an exhaustive set of representatives of 9
equivalence classes of nontoral gradings on $\f4$.

In Section 6 we revisit the gradings on the Albert algebra to
provide a proof that, up to equivalence, the nontoral gradings on
this algebra are the eight ones given in Theorem~\ref{notenJ}.

Having given the quasitori which induce all the nontoral gradings
on $\f4$ we look more closely at the fine gradings in Section 7.
There, we describe those gradings in a twofold way. First we give
the homogeneous components of each grading in a previously fixed
basis of the algebra. This description depends on computational
methods (simultaneous diagonalization algorithms).  In spite of
their computational nature, these descriptions may be interesting
for applications in which explicit calculations are needed. The
second description procedure for the gradings is the exhibition of
models, which make the grading appear in a natural setting with no
appeal to a particular basis. For instance, the fine
$\Z_3^3$-grading on $\f4$ comes from a $\Z_3^3$-grading in $\e6$
which can be described in a very convenient way by using the known
model based in the $\Z_3$-grading $\e6= \s(X_1)\oplus
\s(X_2)\oplus \s(X_3)\oplus\ X_1\otimes X_2\otimes X_3\ \oplus\
X_1^*\otimes X_2^*\otimes X_3^* $ with zero homogeneous component
three copys of the algebra of type $\frak a_2$
\cite[p.\,85]{Adams}. The fine $\Z_2^5$-grading may be seen
directly in $\f4=\der(J)$ but there is also an easy way to see it
in the Tits unified construction for the exceptional Lie algebras
(see, for instance, \cite[p.\,122]{Schafer}), so as the
$\Z_2^3\times\Z$ fine grading. Finally, in 7.4 we justify why this
philosophy can be valid for describing most of the gradings on a
simple Lie algebra. To illustrate  how this works, we describe the
$\Z_2^2\times\Z_4$ and $\Z_2^2\times\Z_8$-nontoral gradings on
$\f4$, and also the $\Z_2^3\times\Z_8$-fine grading equivalent to
the $\Z_2^3\times\Z$-one, by using the model of $\f4$ described in
\cite{mio}, which is  based in an initial $\Z_4$-grading.

\section{Preliminary definitions and results}

\subsection{Group gradings}
Our aim is the study of group gradings on certain nonassociative
algebras over fields. If $V$ is such an algebra and $G$ is an
abelian group, we shall say that the decomposition $V=\oplus_{g\in
G}V_g$ is a \emph{$G$-grading} whenever for all $g,h\in G$,
$V_gV_h\subset V_{gh}$.

In this work we shall use some notions borrowed from the theory of
algebraic groups.  Since we only need linear algebraic groups, all
concepts must be understood in that context. Notice that the group
of automorphisms of the algebra $V$ is an algebraic linear group.
The ground field $F$ will be supposed to be algebraically closed
and of characteristic zero  throughout  this work. There is a deep
relationship between gradings on $V$ and quasitori of the group of
automorphisms $\aut (V)$, according to \cite[\S 3, p.\,104]{enci}.
Following this reference, a commutative algebraic group whose
identity component is an algebraic torus is called an algebraic
{\em quasitorus}. An algebraic linear group is a quasitorus if and
only if there is a basis relative to which the elements of the
quasitorus are simultaneously diagonalizable. If $S$ is a finitely
generated abelian group, then its group of characters $\X(S)=
\hom(S,F^\times)$ is a quasitorus and reciprocally, the group of
characters of a quasitorus turns out to be a finitely generated
abelian group.

If $V=\oplus_{g\in G}V_g$ is a $G$-grading, the map
$\psi\colon\X(G)\to\aut(V)$ mapping each $\alpha\in\X(G)$ to the
automorphism $\psi_\alpha\colon V\to V$ given by $V_g\ni x\mapsto
\psi_\alpha(x):=\alpha(g)x$ is a group homomorphism. In particular
$\psi(\X(G))$ is a quasitorus. And conversely, if $Q$ is a
quasitorus and $\psi\colon Q\to\aut(V)$ is a homomorphism,
$\psi(Q)$ is formed by semisimple automorphisms
(\cite[p.\,99]{Humphreys}) and we have a $\X(Q)$-grading
$V=\oplus_{g\in{\Xs}(Q)}V_g$ given by
$$
V_g=\{x\in V\mid \psi(q)(x)=g(q)x\ \forall q\in Q\}.
$$

When we speak in this paper about a $G$-grading $V=\oplus_{g\in
G}V_g$, we mean that $G$ is generated by the set $\{g\in G\mid
V_g\ne0\}$,
\def\sop{\mathop{\hbox{Supp}}}
called the \emph{support} of the grading and denoted by $\sop(G)$.
In terms of algebraic groups, this is equivalent to the fact that
the homomorphism $\psi\colon Q\to \aut(V)$ is injective. Let us
see it. If $G$ is generated by the support and $\psi(q_0)=\id_V$
for some $q_0\in Q$, $x=\psi(q_0)x=g(q_0)x$ for any $x\in V_g$, so
$g(q_0)=1$ for any $g$ in the support. But since $G$ is generated
by the support, we have $g(q_0)=1$ for any $g\in G$ and this
implies $q_0=1$. \vbox{\hfill
 \xy <1cm, 0cm>:
 \POS(0,0) *+{Q}
 \ar @{->}+(2.5,0)*+{Q'}
 \ar @{->}+(2.5,-2)*+{\aut(V)}
 \POS(2.4,-1.7)\ar @{<..}+(0,1.5)
 \POS(1.2,0.2)*+{_{\pi}}
 \POS(2.8,-0.8)*+{{\psi'}}
 \POS(0.6,-1)*+{\psi}
\endxy}
\vskip -2.6cm
 \vbox{\noindent\hsize=8.7cm

Reciprocally if $\psi$ is a monomorphism, the subgroup $S$ of $G$
generated by the support is $S=\X(Q')$ for some quasitori $Q'$.
The inclusion $i\colon S\to G$ induces by duality an epimorphism
$\pi\colon Q\to Q'$. The grading induced by $S$ comes from a
homomorphism $\psi'\colon Q'\to\aut(V)$ making commutative the
diagram on the right. Since $\psi$ is a monomorphism, so is it
$\pi$. Hence $\pi$ is an isomorphism and by duality the same can
be said about $i$. Thus $S=G$.}

We say that two gradings $V=\oplus_{g\in G}X_g=\oplus_{g'\in
G'}Y_{g'}$ are \emph{isomorphic} if there is $f\in\aut(V)$ and
$\alpha\colon G\to G'$ a group isomorphism such that
$f(X_g)=Y_{\a(g)}$ for any $g\in G$. So, if $\psi\colon
\X(G)\to\aut(V)$ and $\psi'\colon \X(G')\to\aut(V)$ are the
corresponding homomorphisms, and we take $\alpha^*\colon \X(G')\to
\X(G),\alpha^*(\b)=\b\a$, and $\Ad( f)\colon\aut(V)\to\aut(V)$
given by $\Ad( f)(\rho)=f\rho f^{-1}$, the previous condition is
equivalent to the commutativity $\Ad( f)\ \psi'=\psi\alpha^*$.

We say that two gradings $V=\oplus_{g\in G}X_g=\oplus_{g'\in
G'}Y_{g'}$ are \emph{equivalent} if the sets of homogeneous
subspaces are the same up to isomorphism, that is, there are an
automorphism $f\in\aut(V)$ and a bijection between the supports
$\alpha\colon \sop(G)\to \sop(G')$ such that $f(X_g)=Y_{\a(g)}$
for any $g\in \sop(G)$. Our objective is to classify gradings up
to equivalence. A convenient invariant for equivalence is that of
\emph{ type}. Suppose we have a grading on a finite dimensional
algebra, then for each positive integer $i$ we will denote
(following \cite{Hesse}) by $h_i$ the number of homogeneous
components of dimension $i$. In this case we shall say that the
grading is {of type} $(h_1,h_2,\dots,h_l )$, for $l$ the greatest
index such that $h_l\ne0$.   Of course the number $\sum_i i h_i$
agrees with the dimension of the algebra.

Another key notion is that of a {\em coarsening} of a given
grading.  Thus consider an $F$-algebra $V$, a $G$-grading
$V=\oplus_{g\in G}X_g$ and an $H$-grading $V=\oplus_{h\in H}Y_h$.
We shall say that the $H$-grading is a coarsening of the
$G$-grading if and only if  each nonzero homogeneous component
$Y_h$ with $h\in H$  is a direct sum of some homogeneous
components $X_g$. In this case we shall also say that the
$G$-grading is a {\em refinement} of the $H$-grading. Notice that
there is not a relationship between $G$ and $H$, but if
$\psi\colon \X(G)\to\aut(V)$ and $\psi'\colon \X(H)\to\aut(V)$ are
the homomorphisms producing the above gradings, any automorphism
in $\psi'(\X(H))$ acts as a scalar multiple of the identity in
$X_g$ and hence, all of them commute with $\psi(\X(G))$.

The concept of \emph{universal grading group}   is fundamental to
obtain the coarsenings of a given grading. Though this notion is
given in the context of simple Lie algebras in \cite[Section
2.2]{g2}, it can be translated to more general settings.
 Consider as before a nonassociative $F$-algebra $V$, but with a grading
 $V=\oplus_{i\in I}V_i$ which is not supposed to be a group
 grading, that is, each $V_i$ is nonzero and for  any $i,j\in I$
 there is $k\in I$ such that $V_iV_j\subset V_k$.
  Define the group $G_I:=\Z(I)/M$ where $\Z(I)$ is the free
$\Z$-module generated by $I$ and $M$ the $\Z$-submodule generated
by the elements $i+j-k$ such that $0\ne V_iV_j\subset V_k$. We
shall denote by $\bar x$ the equivalence class of $x\in\Z(I)$ in
the quotient $G_I$. We get a $G_I$-grading on $V$ given by
$V=\oplus_{\bar i\in G_I} V'_{\bar i}$ where $V'_{\bar
i}=\sum_{\bar j=\bar i}V_j$. In general this new grading is not
equivalent to the original one. For instance we could have
$G_I=0$, in which case the new grading is the trivial one. But if
the original grading of $V$ is equivalent to an abelian group
grading, then it is easily proved that the canonical map $I\to
G_I$ given by $i\mapsto\bar i$ is injective and the new grading is
equivalent to the original one. Moreover any coarsening of this
grading comes from a group epimorphism as in \cite[Proposition
2(2), p.\,90]{g2}. More precisely, if $V=\oplus_{g\in G}X_g$ is a
$G$-grading and $V=\oplus_{h\in H}Y_h$ is a coarsening, then
 if $G$ is the universal grading group,
there is an epimorphism $p\colon G\to H$ such that
$Y_h=\oplus_{p(g)=h}X_g$. It is this universal property which
suggests the term universal when applied to $G_I$. Getting back to
our context, (abelian) group gradings, we can characterize the
universal grading group in terms of linear algebraic groups. If
$\psi\colon Q\to\aut(V)$ is a grading with $G=\X(Q)$, then $G$ is
the universal grading group if and only if for any coarsening
$\psi'\colon Q'\to\aut(V)$, there is a monomorphism $i\colon Q'\to
Q$ such that $\psi'=\psi i$. This is easily proved applying
duality to the above universal property of $G$. From another
viewpoint, $G$ is the universal group of a grading if and only if
$\X(G)$ is a maximal element in the set of quasitori of $\aut(V)$
producing exactly the same grading.

A group grading is {\em fine} if its unique refinement is the
given grading. If $G$ is the universal group of a grading
$\psi\colon \X(G)\to\aut(V)$, the grading is fine if and only if
$\psi(\X(G))$ is a maximal abelian subgroup of semisimple
elements, which is usually called a \emph{ MAD }(\lq \lq maximal
abelian diagonalizable") in papers about fine gradings, like
\cite{Phys1}. Besides, each MAD $Q\subset\aut(V)$ produces a
$\X(Q)$-fine grading on $V$ such that  $\X(Q)$ is the universal
group of this grading. The converse is   true in the sense just
mentioned but notice that there are fine $G$-gradings such that
$\psi(\X(G))$ is not a MAD, for instance the
$\Z_2^3\times\Z_8$-grading on $\f4$ described in 7.4 or the
$\Z$-grading on $\g2$ described in \cite[Theorem~2,(4)]{g2}. In
particular from the above, the number of fine gradings on $V$ up
to equivalence is the same than the number of MAD's of $\aut(V)$
up to isomorphism (and less than the number of fine gradings up to
isomorphism, in general).

Other notations which will be  used along this paper are the
following.  For a linear algebraic group $G$, and a subset
$S\subset G$, the centralizer of $S$ in $G$ will be denoted by
$\Co_G(S)$. Analogously, by $\No_G(S)$ we shall mean the
normalizer of $S$ in $G$.

\subsection{About the Albert algebra and $\f4$}

Consider the Cayley $F$-algebra $C$, which under our hypothesis on
the ground field   must be isomorphic to Zorn matrices algebra.
Take the standard involution $x\mapsto\bar x$, the norm $n\colon
C\to F$ given by $n(x):=x\bar x$ and the trace $\t\colon C\to F$
defined as $\t(x):=x+\bar x$. Recall that the polar form $f\colon
C\times C\to F$ of $n$ is the symmetric bilinear form
$f(x,y):=\frac{1}{2}\big(n(x+y)-n(x)-n(y)\big)$.  The Albert
algebra $J=\H_3(C) =\{x=(x_{ij})\in M_3(C)\mid x_{ij}=\overline{
x_{ji}}\}$
 is the exceptional  reduced Jordan algebra,
that is, the set of matrices of the form
\begin{equation}\label{uno} \begin{pmatrix}\alpha & o_1 & o_2\cr
\overline{o_1} & \beta & o_3\cr \overline{o_2} &\overline{o_3} &
\gamma
\end{pmatrix} \end{equation} where $\alpha,\beta,\gamma\in F$ and $o_i\in C$
($i=1,2,3$). Since our base field is of characteristic zero, we
can shelter on the linear theory of Jordan algebras and the
product in $J$ may be defined as $x\cdot y:=\frac{1}{2}(xy+yx)$
where juxtaposition stands for the usual matrix product. This
simple Jordan algebra is exceptional in the sense that it is not a
subalgebra of the symmetrization of any associative algebra. It
will be convenient to introduce some notations for further
reference. Thus the element in $J$ obtained in (\ref{uno}) for
$\alpha=1$, $\beta= \gamma=o_i=0$ will be denoted by $E_{1}$
  while the one obtained for $\beta=1$,
$\alpha=\gamma=o_i=0$ will be denoted by $E_{2}$. In a similar way
we can define $E_{3}$ so that $1:=\sum_1^3E_{i}$ is the unit of
$J$. Next, for any $a\in C$ we define by $a\sup 1$ the element in
$J$ obtained  making $\alpha=\beta= \gamma=o_1=o_2=0$, $o_3=a$ in
(\ref{uno}). Define by $a\sup 2$ the one obtained for
$\alpha=\beta=\gamma=o_1=o_3=0$ and $o_2=\bar a$. Finally denote
by $a\sup 3$ the element arising for
$\alpha=\beta=\gamma=o_2=o_3=0$, and $o_1=a$. The multiplication
table of the commutative algebra  $J$ may be summarized in the
following relations:
$$\begin{array}{lll} E_i^2=E_i,\quad  & E_i\,a\sup i=0,\qquad  & a\sup ib\sup
i=f(a,b)(E_j+E_k),\\
E_iE_j=0, & E_i\,a\sup j=\frac{1}{2}a\sup j, & a\sup ib\sup
j=\frac{1}{2}(\bar b\bar a)\sup k, \end{array}$$
 where $(i,j,k)$
is any cyclic permutation of $(1,2,3)$ and $a,b\in C$. Following
Schafer (\cite[(4.41), p.\,109]{Schafer}), any $x\in J$
 satisfies a cubic equation
 $
 x^3-\T(x)x^2+\Q(x)x-\N(x)1=0
$
 where $\T(x),\Q(x),\N(x)\in F$. Moreover, the inversibility of
 $x$ (in Jordan context) is equivalent to the fact that
 $N(x)\ne0$.

For further reference we fix first $(e_1,e_2,u_1,u_2,
u_3,v_1,v_2,v_3)$ the {\em standard basis} of the Cayley algebra
$C$, defined for instance in \cite[Section 3]{g2}, given by the
relations
$$\begin{array}{rll}
e_1u_j&=u_j&=u_je_2,\\
e_2v_j&=v_j&=v_je_1,\end{array}\
\begin{array}{rll}
u_iu_j&=v_k&=-u_ju_i, \\
-v_iv_j&=u_k&=v_jv_i,\end{array}\ \begin{array}{rl} u_iv_i&=e_1,\\
v_iu_i&=e_2,\end{array}
$$
where $e_1$ and $e_2$ are orthogonal idempotents, again $(i,j,k)$
is any cyclic permutation of $(1,2,3)$, and the remaining
relations are null. Thus we can fix
 our {\em standard basis} of the Albert algebra:
$$
B=(E_1,E_2,E_3,e_1\sup{3},e_2\sup{3},u_1\sup{3},u_2\sup{3},
u_3\sup{3},v_1\sup{3},v_2\sup{3},v_3\sup{3},
e_2\sup{2},e_1\sup{2},-u_1\sup{2},-u_2\sup{2}, -u_3\sup{2},$$
$$-v_1\sup{2},-v_2\sup{2},-v_3\sup{2},
 e_1\sup{1},e_2\sup{1},u_1\sup{1},u_2\sup{1},
u_3\sup{1},v_1\sup{1},v_2\sup{1},v_3\sup{1}).
$$

Let us define now the group $F_4:=\aut(J)$ and its Lie algebra
$\f4=\der(J)$. Recall (see for instance \cite[p.\,285]{Jac}) that
the automorphism group $F_4 $ and the automorphism group
$\F4:=\aut(\f4)$ are isomorphic via the map $\Ad\colon F_4\to\F4 $
such that $\Ad(f)d:=fdf^{-1}$ for any $f\in F_4$ and $d\in\f4$.
This isomorphism of algebraic groups provides a tool for
translating gradings from the Albert algebra to $\f4$ and
conversely. However, unfortunately this translating tool does not
preserve equivalence. This phenomenon is similar to the one
explained in \cite[Section 4]{g2} in the similar context produced
by the analogue group isomorphism $\Ad\colon
G_2=\aut(C)\to\aut(\g2)=\aut(\der(C))$.

\subsection{Maximal torus of $F_4$}

Let us fix certain maximal torus of $F_4$. We use the standard
basis of $J$ above defined. Define now the maximal torus $\To_0$
of $F_4$ whose elements are the automorphisms of $J$ which are
diagonal relative to $B$. This is isomorphic to $(F^\times)^4$ and
it is easy to see that the matrix of any such automorphism
relative to $B$ is
\begin{equation}\label{maxtor}
\hbox{diag}\Big(1,1,1,\alpha ,\frac{1}{\alpha },\beta ,\gamma
,\frac{{\delta }^2}{\alpha \,\beta \,\gamma },\frac{1}{\beta
},\frac{1}{\gamma },
  \frac{\alpha \,\beta \,\gamma }{{\delta }^2},\delta ,\frac{1}{\delta },\frac{\alpha \,\beta }{\delta },\frac{\alpha \,\gamma }{\delta },
  \frac{\delta }{\beta \,\gamma },\frac{\delta }{\alpha \,\beta },
  \frac{\delta }{\alpha \,\gamma },\frac{\beta \,\gamma }
  {\delta },\frac{\delta }{\alpha },$$
  $$\frac{\alpha }{\delta },\frac{\beta }{\delta },\frac{\gamma }{\delta },\frac{\delta }{\alpha \,\beta \,\gamma },\frac{\delta }{\beta },\frac{\delta }{\gamma },
  \frac{\alpha \,\beta \,\gamma }{\delta }\Big)
\end{equation}
for some $\alpha$, $\beta$, $\gamma$, $\delta\in F^\times$.
  Define now $t_{\alpha,\beta,\gamma,\delta}$ as the automorphism
  in $\To_0$ whose matrix relative to $B$ is just the above one.

 Consider a grading of an algebra $A$  given by a group
  homomorphism $\rho\colon \X(G)\to\aut(A)$.
  This grading is
  said to be \emph{toral} if $\rho(\X(G))$ is contained in some torus
  of the algebraic group $\aut(A)$.

Notice that  by means of  the above isomorphism $\Ad\colon
F_4\to\F4 $, a grading on $J$ is toral if and only if it is
applied to a toral grading on $\f4$. This does not imply that the
number of nontoral gradings (up to equivalence) on $J$ and on
$\f4$ is necessarily the same, because the mechanism of
translating gradings does not preserve equivalence.

  An useful characterization of
  the torality of a grading on $J$ is the following:

\begin{pr}\label{toynoto}  A grading on the Albert algebra $J$ is toral if and
only if the elements of the standard basis  of $J$ (or any
conjugated basis) are homogeneous.
\end{pr}
Proof.  Consider a toral grading induced by automorphisms
$\{t_i\}_{i\in I}$ all of them contained in the previous maximal
torus $\To_0$ of $J$ (which supposes no restriction because any
other maximal tori is conjugated to $\To_0$). Then obviously the
simultaneous diagonalization of $J$ relative to the family
$\{t_i\}_{i\in I}$ provides the original grading and the elements
in the standard basis are homogeneous. Conversely, if this holds,
then any element in this basis is an eigenvector of any of the
grading automorphisms. Thus these automorphisms are in the maximal
torus specified before and the grading is toral.\qed

In particular, the proof of this proposition implies that if a
grading on $J$ is toral, then there are three orthogonal
idempotents contained in the zero component (because up to
conjugacy, $\{E_1,E_2,E_3\}$ is contained in such component).
Another way of checking the torality is to look at the rank of the
zero part of the induced grading on $\f4$, it will be toral in
case this rank is 4 (see \cite[Subsection 2.4]{g2}).

By using  this set of idempotents, we can give more information
about the general form of any semisimple automorphism of $J$,
because they are the building blocks of the grading sets. First
notice that there is a group monomorphism $D_4\to F_4$. Indeed, if
$U\in\mathop{\rm O}(C,n)=\{g\in\gl(C)\mid n(x,y)=n(g(x),g(y))\
\forall x,y\in C\}$, there are $U',U''\in\mathop{\rm O}(C,n)$ such
that $U(xy)=U'(x)U''(y)$ for any $x,y\in C$. This is called
\textsl{global triality principle} in \cite[Th.\,3,
p.\,90]{octon}.
 Let $\Psi_U\colon
J\to J$ given by
\begin{equation}\label{psiu}
 \Psi_U(E_i)=E_i,\,\Psi_U(x^{(i)})=(f_i(x))^{(i)}\text{ for
 }f_1=U,\,f_2=U'',\,f_3=U'.
\end{equation}
It is  easy to check that  $\Psi_U $ is an automorphism of $J$
fixing each idempotent. Conversely,    if $\psi $ is an
automorphism of $J$ fixing each idempotent, there exists
$U\in\mathop{\rm O}(C,n)   $ such that $\psi=\Psi_U$.   Moreover,
any semisimple automorphism  in $F_4$ is toral, since $F_4$ is a
connected group, so that  up to conjugacy it is contained in
$\To_0$, fixes each $E_i$ and is of the form $\Psi_U$ for some
$U\in\mathop{\rm O}(C,n)   $.

  Besides, there are two
outstanding automorphisms which deserve some consideration. The
  first one $\theta\colon J\to J$ applies
$E_i\mapsto E_{i+1}$ and $x^{(i)}\mapsto x^{(i+1)}$ cyclically.
The second one $\vartheta\colon J\to J$ fixes $E_3$ permuting
$E_1$ and $E_2$, and acts as the identity on elements $x\sup{3}$
while $x\sup{1}\mapsto x\sup{2}\mapsto x\sup{1}$. These
automorphisms fix the set $\{E_1,E_2,E_3\}$ and are of order three
and two respectively. Our previous argument shows that any
semisimple element in $F_4$ fixing the set $\{E_1,E_2,E_3\} $ is a
composition of one automorphism in $ \{\id,
\theta,\theta^2,\vartheta,\vartheta\theta,\vartheta\theta^2\}$
with one in $\{\Psi_U\mid U\in\mathop{\rm O}(C,n)\}   $.

\subsection{Tits construction of the Albert algebras}\label{Tits}
There is another way in which the Albert algebra   can be
constructed. Let us start with the $F$-algebra $A=  M_3(F)$ and
denote by $\T_A,\Q_A,\N_A\colon A\to F$ the coefficients of the
generic minimal polynomial such that
$x^3-\T_A(x)x^2+\Q_A(x)x-\N_A(x)1=0$ for all $x\in A$. Recall that
if $x=(x_{ij})\in A$ then $\T_A(x)= \sum_1^3 x_{ii}$,
$$\Q_ A(x)= -x_{12}x_{21} + x_{11}x_{22}-x_{13}
x_{31}-x_{23}x_{32}+ x_{11}x_{33}+x_{22} x_{33},$$  and
$\N_A(x)=\det(x)$. Define also the quadratic map $\sharp\colon
A\to A$ given by $x^\sharp:=x^2-\T_A(x)x+Q_A(x)1$. For any $x,y\in
A$ denote   $x\times y:=(x+y)^\sharp-x^\sharp-y^\sharp$, and
$x^*:=\frac{1}{2}x \times 1=\frac{1}{2}\T_A(x) 1-\frac{1}{2}x$.
Finally consider the Jordan algebra $A^+$ whose underlying vector
space agrees with that of $A$ but whose product is $x\cdot
y=\frac{1}{2}(xy+yx)$. Next, define in $A^3:=A\times A\times A$
the product
$$\begin{array}{c}(a_1,b_1,c_1)(a_2,b_2,c_2):=\\
\big(a_1\cdot a_2+(b_1c_2)^*+(b_2c_1)^*,
a_1^*b_2+a_2^*b_1+\frac12(c_1\times c_2),
c_2a_1^*+c_1a_2^*+\frac12(b_1\times b_2)\big).\end{array}$$
\smallskip Then $A^3$ with this product is isomorphic to
$J=\H_3(C)$. This is the so called Tits construction of the Albert
algebra.   This allows us to identify $J$ with the algebra $A^3$
in the rest of this section. For further reference we shall recall
that the norm $N$ (module the identification of $J$ with the Tits
construction) is given by
 \begin{equation}\label{rel}
\N(a,b,c)=\N_A(a)+\N_A(b)+\N_A(c)-\T_A(abc)
\end{equation}
 for any
$a,b,c\in A$ (see \cite[p.\,525]{Boi}).

One of the relevant facts on Tits construction from our viewpoint
is that it allows to embed $\aut(A)$ in $F_4$ via  the map
$\aut(A)\to F_4$ given by $f\mapsto f^\bullet$ where
$f^\bullet\colon J\to J$ is the automorphism such that
$f^\bullet(x,y,z):=(f(x),f(y),f(z))$. As a further consequence we
will be able to get gradings on $J$ coming from gradings in the
associative algebra $A$ via this monomorphism of algebraic groups.

\section{Inducing gradings on the Albert algebra}\label{type}

 Some remarkable subalgebras of $J=\H_3(C)$ induce gradings on
$J$ by means of a particular embedding of its automorphism group
in $F_4$. For instance, the automorphism group of $\H_3(F)$ can be
considered as a subgroup of $F_4$ providing thus a source of
gradings on $J$: those whose grading automorphisms come from
automorphisms of $\H_3(F)$. More generally, if $V$ is an algebra
such that $\aut(V)$ is a subgroup of $F_4$, then any grading in
$V$ will also  induce a grading on $J$. We shall see that this
happens for instance for the octonion $F$-algebra. It turns out
that this idea for inducing gradings on $J$ provides a great
number of the gradings existing on $J$.

\subsection{Gradings from octonions}\label{gfo}
Let $C$ be the Cayley algebra. We consider as in the previous
section the standard basis of $C$. We shall also need the maximal
torus of $G_2:=\aut(C)$ given by the automorphisms $t_{\a,\b}$
whose matrix relative to the standard basis is
$$\diag(1,1,\alpha,\beta,(\a\b)^{-1},\a^{-1},\b^{-1},\a\b).$$
It has been first proved in \cite{Alb1}, and then in
\cite[Subsection 3.3]{g2}, that up to equivalence the unique
nontoral grading on $C$ is the $\Z_2^3$-grading whose order-two
grading automorphisms are $\{t_{1,-1},t_{-1,1},f_0\}$ where $f_0$
is the automorphism whose matrix relative to the standard basis is
$$\small
\begin{pmatrix}
 0 & 1 & 0 & 0 & 0 & 0 & 0 & 0\cr
 1 & 0 & 0 & 0 & 0 & 0 & 0 & 0\cr
 0 & 0 & 0 & 0 & 0 & 1 & 0 & 0\cr
 0 & 0 & 0 & 0 & 0 & 0 & 1 & 0\cr
 0 & 0 & 0 & 0 & 0 & 0 & 0 & -1\cr
 0 & 0 & 1 & 0 & 0 & 0 & 0 & 0\cr
 0 & 0 & 0 & 1 & 0 & 0 & 0 & 0\cr
 0 & 0 & 0 & 0 & -1& 0 & 0 & 0
\end{pmatrix}.
$$
 The construction $J=\H_3(C)$ of the Albert algebra  is
particularly interesting for extending derivations and
automorphisms from $C$ to $J$. More precisely if $f\in\aut(C)$ is
an automorphism of $C$ then we can construct the automorphism
$\hat f$ of $J$ fixing the idempotents $E_{i}$ and such that $\hat
f(o\sup i):=f(o)\sup i$ for $i=1,2,3$ and any $o\in C$. This
provides a monomorphism of algebraic groups $i\colon G_2\to F_4$
such that $f\mapsto\hat f$. By differentiating at $1$ we get a
monomorphism of Lie algebras $\text{d}i(1)\colon
\der(C)=\g2\to\der(J)=\f4$ mapping each derivation $d\in\der(C)$
to the derivation $\hat d\in\der(J)$ annihilating the idempotents
and making $\hat d(o^i)=d(o)^i$, for $i=1,2,3$. This, of course,
has an immediate application to gradings: any grading on $C$
induces a grading on $J$. Indeed, a $G$-grading on $C$ comes from
an algebraic group homomorphism $\rho\colon \X(G)\to\aut(C)$,
therefore $i\rho\colon \X(G)\to\aut(J)$ provides a $G$-grading on
the Albert algebra. This device gives a first source of gradings
on $J$, namely, all those coming from  gradings on $C$. Since $i$
maps tori of $\aut(C)$ to tori of $\aut(J)$, toral gradings on $C$
induce toral gradings on $J$. The unique nontoral grading on $C$
up to equivalence provides a $\Z_2^3$-grading on $J$ induced by
the automorphisms
$\{\widehat{t_{1,-1}},\widehat{t_{-1,1}},\widehat{f_0}\}$ whose
homogeneous spaces are
 \begin{eqnarray}\label{grad1}
 J_{1,1,1}&=&\span{E_1, E_2,
E_3, 1\sup{3}, 1\sup{2},
  1\sup{1}},\cr
 J_{1,1,-1}&=&\span{
 (-{e_1} + {e_2})\sup{3}, (-{e_1} +
 {e_2})\sup{2}, (-{e_1} + {e_2})\sup{1}},\cr
 J_{1,-1,1}&=&\span{({u_2} + {v_2})\sup{3}, ({u_2} + {v_2})\sup{2},
 ({u_2} + {v_2})\sup{1}},\cr
 J_{-1,1,1}&=&\span{({u_1} + {v_1})\sup{3}, ({u_1} + {v_1})\sup{2},
 ({u_1} + {v_1})\sup{1}},\cr
 J_{1,-1,-1}&=&\span{(-{u_2} + {v_2})\sup{3},
 (-{u_2} + {v_2})\sup{2}, (-{u_2} + {v_2})\sup{1}},\cr
 J_{-1,1,-1}&=&\span{(-{u_1} + {v_1})\sup{3},
 (-{u_1} + {v_1})\sup{2}, (-{u_1} + {v_1})\sup{1}},\cr
 J_{-1,-1,1}&=&\span{(-{u_3} + {v_3})\sup{3}, (-{u_3} + {v_3})\sup{2},
 (-{u_3} + {v_3})\sup{1}},\cr
 J_{-1,-1,-1}&=&\span{({u_3} + {v_3})\sup{3},
 ({u_3} + {v_3})\sup{2}, ({u_3} + {v_3})\sup{1}}.
 \end{eqnarray}
 For this grading we have $h_i=0$ except $h_3=7$, $h_6=1$. Thus
 the algebra is of type $(0,0,7,0,0,1)$.
It is easy to prove that the subalgebra of $\f4$ whose elements
are those $d\in\f4$ such that
$[d,\widehat{t_{1,-1}}]=[d,\widehat{t_{-1,1}}]=[d,\widehat{f_0}]=0$
is three-dimensional. Thus, the grading on $\f4$ induced by
$\{\Ad(\widehat{t_{1,-1}}),\Ad(\widehat{t_{-1,1}}),\Ad(\widehat{f_0})\}$
is nontoral since its zero homogeneous component has dimension
$3$, hence its rank is less than $4=\rank(\f4)$. We summarize the
results in this subsection in the following:
\begin{pr}\label{dece} The unique nontoral grading on $J$ coming from a
grading on $C$ is the above grading \textrm{(\ref{grad1})} up to
equivalence.
\end{pr}

\subsection{Gradings on $\H_3(F)$}
It is known from \cite[p.\,184-185]{JacJor} that any automorphism
of $\H_3(F)$ is of the form $\In(p)\colon x\mapsto pxp^{-1}$ where
$p$ can be taken in the group $\SO(3):=\SO(3,F)$ consisting of
those $3\times 3$ matrices $x$ with entries in $F$ such that
$xx^t=1$  and $\det(x)=1$ (the notation $x^t$ meaning the
transpose of the matrix $x$). Thus we have an algebraic group
isomorphism $\In\colon\SO(3)\to\aut(\H_3(F))=:G$. Now consider the
algebra $  M_3(C)$ with the usual matrix product. This contains,
as a subalgebra, $ M_3(F)$ and any element in this subalgebra
associates with any two other elements in $M_3(C)$. Since
$\H_3(F)$ is a (Jordan) subalgebra of $J=\H_3(C)$, for any
$p\in\SO(3)$, the map $\In(p)\colon J\to J$ such that $x\mapsto
pxp^{-1}$ (products in $ M_3(C)$) is an automorphism of $J$. Thus
we have an algebraic group monomorphism $\In\colon\SO(3)\to  F_4$
and thus a monomorphism $\SO(3)\cong G\to F_4=\aut(J)$ so that we
can extend automorphisms from $\H_3(F)$ to $J$ (in fact this
monomorphism maps $\In(p)$
 seen as an element in $G$ to $\In(p)$ as element in $F_4$).
  So we have constructed
 a tool for translating gradings from $H_3(F)$ to $J$.
This suggests the convenience of describing group gradings on
$\H_3(F)$. The   work \cite{Ivan2} contains an exhaustive and deep
description of gradings on some simple Jordan and Lie algebras, in
particular the gradings on $\H_n(F)$ could be obtained along the
line of this work. But in our particular setting it is worth to
find this description by using elementary algebraic group theory.
This will provide a comfortable landscape for testing geometric
tools and will provide a more self contained exposition.

It is well known that a maximal torus $P$   of $\SO(3)$ is given
by the matrices of the form
 $$p_{\a,\b}:=\begin{pmatrix}1 & 0 &0\cr
 0 & \alpha & \beta\cr
 0 & -\beta & \alpha
 \end{pmatrix}$$
 with $\a,\b\in F$   such that $\alpha^2+\beta^2=1$. Then denote by
 $\tau_{\a,\b}:=\In(p_{\a,\b})$ the corresponding element in
 $G=\aut(\H_3(F))$.
 The set of all $\tau_{\a,\b}$ is a maximal torus $\Tau$ of $G$ and the
 set of eigenvalues of $\tau_{\a,\b}$ is $S_{\a,\b}=\{1,z,z^{-1},z^2,z^{-2}\}$
 for
 $z:=\a+i\b$. Supposing $\vert S_{\a,\b}\vert=5$, we find for $\tau_{\a,\b}$ the
 following eigenspaces
\begin{equation}\label{gr1}\begin{array}{ll}
 \H_3(F)_1=\span{E_1,E_2+E_3}, &
 \H_3(F)_z=\span{-i\sup{3}+1\sup{2}},\\
 \H_3(F)_{z^{-1}}=\span{i\sup{3}+1\sup{2}}, &
 \H_3(F)_{z^2}=\span{-i E_{2}+i E_3+1\sup{1}},\\
 \H_3(F)_{z^{-2}}=\span{i E_{2}-i E_3+1\sup{1}}, &
 \end{array}\end{equation}
where the subindex   indicates the eigenvalue of $\tau_{\a,\b}$.
This gives a $\Z$-grading of $\H_3(F)$, with $n$-th component
 $\H_3(F)_{z^n}$.
  This is toral and fine (as it is produced by the whole torus $\Tau$).
 Any other toral grading of $\H_3(F)$ is a coarsening of this.

 For $\vert S_{\a,\b}\vert<5$ we have the following excluding possibilities:
\begin{itemize} \item $1=z$ which gives the trivial grading.  \item $1=z^2$
which excluding the previous case implies $z=-1$.  This is the $\Z_2$-grading
induced by the involutive automorphism $\tau_{-1,0}$ and it is given by
\begin{eqnarray}\label{gr2}
\H_3(F)_1=\span{E_1,E_2,E_3,1\sup{1}},\quad
\H_3(F)_{-1}=\span{1\sup{2},1\sup{3}}.  \end{eqnarray}
\item $z=z^{-2}$ implying $z^3=1$. Ruling out previous cases, this
induces a $\Z_3$-grading coming for instance from
$\tau_{-1/2,\sqrt{3}/2}$. The grading is
\begin{eqnarray}\label{gr3} \H_3(F)_1&=&\span{E_1,E_2+E_3},\cr
\H_3(F)_{z}&=&\span{1\sup{2}-i\sup{3},iE_2-iE_3+1\sup{1}},\cr
\H_3(F)_{z^2}&=&\span{1\sup{2}+i\sup{3},-iE_2+iE_3+1\sup{1}}.
\end{eqnarray}
\item $z^2=z^{-2}$ implying $z^4=1$. This gives a $\Z_4$
-grading coming from $\tau_{0,1}$. The grading is
\begin{equation}\begin{array}{ll}\label{gr4}
\H_3(F)_1=\span{E_1,E_2+E_3}, &
\H_3(F)_z=\span{-i\sup{3}+1\sup{2}},\\
\H_3(F)_{z^{-1}}=\span{i\sup{3}+1\sup{2}}, &
\H_3(F)_{z^2}=\span{E_{2}- E_3,1\sup{1}}.
\end{array}\end{equation}
\end{itemize}

The gradings in (\ref{gr1})-(\ref{gr4})   are therefore the unique
cyclic (hence necessarily toral, see Appendix) gradings. To find
the rest of the gradings on $\H_3(F)$ we compute the centralizers
of the grading automorphisms producing the previous gradings. The
computations of the various centralizers are easy taking advantage
of the isomorphism $\SO(3)\to G$. For any $\tau_{\a,\b}$ with
$\vert S_{\a,\b}\vert=5$ we have $\Co_G(\tau_{\a,\b})=\Tau$ (for
this, we only need to prove that the centralizer of $p_{\a,\b}$ in
$\SO(3)$ is the maximal torus $P$), so that  the grading
(\ref{gr1}) is fine, as mentioned. The centralizer of
$\tau_{-1,0}$ has two connected components: the identity component
is the maximal torus $\Tau$ and
$\Co_G(\tau_{-1,0})/\Tau\cong\Z_2$. Working in $\SO(3)$, the
identity component of the centralizer of $p_{-1,0}$ is the torus
$P$ while its other component is $sP$ where $s=-E_1+1\sup{1}$.
Taking any $\tau\in\Tau$ the grading induced by
$\{\tau_{-1,0},\tau\}$ is some of (\ref{gr1})-(\ref{gr4}). But if
we consider the grading $\{\tau_{-1,0},\In(s)\}$, we get the
$\Z_2\times\Z_2$-grading given by
$\H_{1,1}=\span{E_1,E_2+E_3,1\sup{1}}$,
$\H_{1,-1}=\span{1^{(2)}-{1^{(3)}}}$, $\H_{-1,1}=\span{E_2-E_3}$,
$\H_{-1,-1}=\span{1^{(2)}+{1^{(3)}}}$, which is isomorphic to:
\begin{equation}\label{gr5}\begin{array}{ll}
 \H_{1,1}=\span{E_1,E_2,E_3},\quad
 &\H_{1,-1}=\span{1\sup{1}},\\
 \H_{-1,1}=\span{1\sup{2}},\quad
 &\H_{-1,-1}=\span{1\sup{3}}.
\end{array}\end{equation}
 On the other hand
 $\Co_G(\span{\tau_{-1,0},\In(s)})=\Co_G(\tau_{-1,0})$,
 which implies that (\ref{gr5}) is fine (and nontoral taking into
 account \cite[Theorem~1]{g2}).  The grading (\ref{gr3}) is
 produced by $t_{-1/2,\sqrt{3}/2}$ whose centralizer is $\Tau$.
 The grading (\ref{gr4}) is produced by $\tau_{0,1}$ whose
 centralizer is again $\Tau$.

 Summarizing all the above results we claim:
 \begin{te}
  Any  nontrivial grading
 on $\H_3(F)$ is equivalent to one of the gradings
 \textrm{(\ref{gr1})-(\ref{gr5})} above.
\end{te}

  The gradings induced on the Albert algebra by the
monomorphism $\aut(H_3(F))\to F_4$ are all of them toral, but
these automorphisms of $J$ coming from $H_3(F)$ commute with the
automorphisms coming from $C$, fact which will provide larger
abelian sets of semisimple automorphisms and hence a source of
nontoral gradings on $J$.

\subsection{A family of nontoral gradings on the Albert algebra}
We now construct a machinery for building gradings on the Albert
algebra $J=\H_3(C)$ by mixing gradings on $C$ with that on
$\H_3(F)$. We must start with the simple observation that
$J=\H_3(C)\cong\H_3(F)\otimes F\oplus \K_3(F)\otimes C_0$ (tensor
product of $F$-spaces) where $\K_3(F)$ is the subspace of $3\times
3$ skewsymmetric matrices with entries in $F$ and $C_0:=\{x\in
C\mid \t(x)=0\}$ the subspace of trace zero elements in $C$. The
above isomorphism is given by $E_i\mapsto E_i\otimes 1$,
$1\sup{i}\mapsto 1\sup{i}\otimes 1$ and for $x\in C_0$,
$x\sup{i}\mapsto (e_{jk}-e_{kj})\otimes x$, being $(i,j,k)$ any
cyclic permutation of $(1,2,3)$ and $e_{ij}$ the elementary
$(i,j)$-matrix in $M_3(F)$. Taking in $J'=\H_3(F)\otimes F\oplus
\K_3(F)\otimes C_0$, subspace of $ M_3(F)\otimes C $, the product
$(c\otimes x)\cdot (d\otimes y)=\frac12((c\otimes x)(d\otimes
y)+(d\otimes y)(c\otimes x))$ for $(c\otimes x)(d\otimes
y)=cd\otimes xy$, $J'$ is a Jordan subalgebra of $ M_3(F)\otimes C
$ such that the previous vector space isomorphism between $J$ and
$J'$ is an algebra isomorphism.
 Module this identification of $J$ with
$J'$, the embedding of $G_2$ in $F_4$ described in Subsection
\ref{gfo} can be seen in the following way. Given $f\in G_2$, take
$\hat f$ the restriction of $\id\otimes f\in\gl( M_3(F)\otimes C)$
to $J'$, which is an automorphism of $J'$ (notice that $C_0$ is
$f$-invariant for any $f\in G_2$).
 On the other hand given any automorphism $g$ of
$\H_3(F)$, $g$ is the restriction to $\H_3(F)$ of an automorphism
$g$ of $ M_3(F)$ commuting with the transposition involution.
Hence $g(\K_3(F))\subset\K_3(F)$ and we can define $\tilde g$ as
the restriction of $g\otimes \id\in\gl( M_3(F)\otimes C)$ to $J'$,
which is an automorphism of $J'$. A trivial though remarkable fact
is the commutativity $\hat f\tilde g=\tilde g\hat f$ for any $f\in
G_2$, $g\in \aut(\H_3(F))$. Thus we have:
\begin{te}
Let $\{f_1,\ldots,f_k\}\subset G_2$ and $\{g_1,\ldots,g_n\}\subset
\aut(\H_3(F))$ be commutative sets of diagonalizable
automorphisms. Then
$\{\widehat{f_1},\ldots,\widehat{f_k},\widetilde{g_1},\ldots,
\widetilde{g_n}\}\subset F_4$ is a commutative set of
diagonalizable automorphisms of $J$. In particular if $C$ is
graded by a group $G_1$ and $\H_3(F)$ is graded by a second group
$G_2$, then the Albert algebra $J$ is $G_1\times G_2$-graded.
\end{te}

It is always the case that the grading induced by
$\{\widehat{f_1},\ldots,\widehat{f_k},\widetilde{g_1},\ldots,
\widetilde{g_n}\}$ is a refinement of the one given by
$\{\widehat{f_1},\ldots,\widehat{f_k}\}$. Besides if one of the
gradings $\{\widehat{f_1},\ldots,\widehat{f_k}\}$ or
$\{\widetilde{g_1},\ldots, \widetilde{g_n}\}$ is nontoral, the
refinement is also nontoral.
 These results  allow us to combine gradings on $C$ with
gradings on $\H_3(F)$. Thus, if we pick the (unique up to
equivalence) nontoral grading on $C$ given by
$\{t_{1,-1},t_{-1,1},f_0\}$ and any of the gradings
(\ref{gr1})-(\ref{gr5}) plus the trivial grading,  which are given
respectively by: $\{\tau_{\a,\b}\}$ (with $\vert
S_{\alpha\beta}\vert=5$), $\{\tau_{-1,0}\}$,
$\{\tau_{-1/2,\sqrt{3}/2}\}$, $\{\tau_{0,1}\}$,
$\{\tau_{-1,0},\In(s)\}$ and $\{1\}$, we get six nontoral gradings
on $J$ which are given in the next

\begin{co}
The following six nonequivalent gradings of the Albert algebra are
nontoral:
\begin{itemize}
\item[a)] $\{\widehat{t_{1,-1}},\widehat{t_{-1,1}},\widehat{f_0},
\widetilde{\tau_{\a,\b}}\}$ with $\vert S_{\a,\b}\vert=5$. This is
a $\Z_2^3\times\Z$-grading.
\item[b)] $\{\widehat{t_{1,-1}},\widehat{t_{-1,1}},\widehat{f_0},
\widetilde{\tau_{-1,0}}\}$. This is a $\Z_2^4$-grading.
\item[c)] $\{\widehat{t_{1,-1}},\widehat{t_{-1,1}},\widehat{f_0},
\widetilde{\tau_{-\frac{1}{2},\frac{\sqrt{3}}{2}}}\}$. This is a
$\Z_2^3\times\Z_3$-grading.
\item[d)] $\{\widehat{t_{1,-1}},\widehat{t_{-1,1}},\widehat{f_0},
\widetilde{\tau_{0,1}}\}$. This is a $\Z_2^3\times\Z_4$-grading.
\item[e)] $\{\widehat{t_{1,-1}},\widehat{t_{-1,1}},\widehat{f_0},
\widetilde{\tau_{-1,0}},\widetilde{\In(s)}\}$. This is a
$\Z_2^5$-grading.
\item[f)] $\{\widehat{t_{1,-1}},\widehat{t_{-1,1}},\widehat{f_0}\}$.
This is a $\Z_2^3$-grading.
\end{itemize}
\end{co}
We now describe explicitly the six previous gradings. \medskip

a) If $\epsilon=\pm 1$ then the grading is {\small
 \begin{eqnarray}\label{nt1}
 J_{0000}=\span{E_1, E_2 + E_3},& J_{000\epsilon}=
 \span{{-i\ \epsilon\sup{3} + 1\sup{2}}},\cr
 J_{001\epsilon}=\span{\imag\epsilon({e_1} -  {e_2})\sup{3} + ({e_1} -
 {e_2})\sup{2}}, & J_{0010}=\span{(-{e_1} +
 {e_2})\sup{1}},\cr
 J_{010\epsilon}=\span{-i\epsilon({u_2}+{v_2})\sup{3} - ({u_2} +
 {v_2})\sup{2}}, & J_{0100}=\span{({u_2} + {v_2})\sup{1}},\cr
 J_{100\epsilon}=\span{-i\epsilon({u_1} +{v_1})\sup{3} - ({u_1} +
 {v_1})\sup{2}}, & J_{1000}=\span{({u_1} +
 {v_1})\sup{1}},\cr
 J_{011\epsilon}=\span{i\epsilon({u_2} - {v_2})\sup{3} + ({u_2} -
 {v_2})\sup{2}}, & J_{0110}=\span{({u_2} - {v_2})\sup{1}},\cr
 J_{110\epsilon}=\span{i\epsilon({u_3} - {v_3})\sup{3} + ({u_3} -
 {v_3})\sup{2}}, & J_{1100}=\span{({u_3} - {v_3})\sup{1}},\cr
 J_{101\epsilon}=\span{i\epsilon({u_1} - {v_1})\sup{3} + ({u_1} -
 {v_1})\sup{2}}, &
 J_{1010}=\span{({u_1} - {v_1})\sup{1}},\cr
 J_{111\epsilon}=\span{-i\epsilon({u_3}+{v_3})\sup{3} - ({u_3} +
 {v_3})\sup{2}}, &J_{1110}=\span{({u_3} + {v_3})\sup{1}},\cr
 J_{0002}=\span{-i (E_2 - E_3) +
 1\sup{1}},& J_{000-2}=\span{-i (E_2 - E_3) -
 1\sup{1}}.
 \end{eqnarray}}
 This grading is a $\Z_2^3\times\Z$-grading of type
 $(25,1)$.\medskip

 b) This is the $\Z_2^4$-grading
 {\small\begin{eqnarray}\label{nt2}
 J_{0000}=\span{E_1, E_2, E_3, 1\sup{1}}, &
 \ J_{0001}=\span{1\sup{3}, 1\sup{2}},\cr
 J_{0010}=\span{({e_1} -
 {e_2})\sup{1}}, & J_{0100}=\span{({u_2} + {v_2})\sup{1}},\cr
 J_{1000}=\span{({u_1} + {v_1})\sup{1}},& J_{1100}=\span{({u_3} -
 {v_3})\sup{1}},\cr
  J_{1010}=\span{({u_1} - {v_1})\sup{1}},
 & J_{1001}=\span{({u_1} + {v_1})\sup{3}, ({u_1} +
 {v_1})\sup{2}}, \cr
 J_{0110}=\span{(-{u_2} + {v_2})\sup{1}},
 & J_{0101}=\span{({u_2} + {v_2})\sup{3}, ({u_2} +
 {v_2})\sup{2}},\cr
 J_{0011}=\span{({e_1} - {e_2})\sup{3}, ({e_1} -
 {e_2})\sup{2}}, & J_{1110}=\span{({u_3} + {v_3})\sup{1}},\cr
J_{1101}=\span{({u_3} - {v_3})\sup{3}, ({u_3} - {v_3})\sup{2}}, &
J_{1011}=\span{({u_1} - {v_1})\sup{3}, ({u_1} - {v_1})\sup{2}},\cr
J_{0111}=\span{ ({u_2} - {v_2})\sup{3}, (-{u_2} + {v_2})\sup{2}},
& J_{1111}=\span{({u_3} + {v_3})\sup{3}, ({u_3} + {v_3})\sup{2}},
 \end{eqnarray}}
 which is of type $(7,8,0,1)$.\medskip

 c) This is the $\Z_2^3\times\Z_3$-grading given by
 {\small
 \begin{eqnarray}\label{nt3}
 J_{0000}=\span{E_1, E_2 + E_3}, & J_{1000}=\span{({u_1} +
 {v_1})\sup{1}},\cr
J_{0100}=\span{({u_2} +
 {v_2})\sup{1}}, & J_{0010}=\span{({e_1} - {e_2})\sup{1}},\cr
 J_{0110}=\span{({u_2} - {v_2})\sup{1}}, & J_{1010}=
 \span{({u_1} - {v_1})\sup{1}},\cr
 J_{1100}=\span{({u_3} - {v_3})\sup{1}},&
  J_{1110}=\span{({u_3} + {v_3})\sup{1}},\cr
  J_{0001}=\span{-i 1\sup{3} +
    1\sup{2}, i E_2 - i E_3 +
    1\sup{1}}, &
  J_{1001}=\span{i( {u_3} + {v_3})\sup{3} - ({u_3} +
  {v_3})\sup{2}},\cr
  J_{0101}=\span{i({u_2} +{v_2})\sup{3} + ({u_2} + {v_2})\sup{2}},
  & J_{0011}=\span{i({e_1} - {e_2})\sup{3} + ({e_1} -
  {e_2})\sup{2}},\cr
  J_{0111}=\span{i({u_2} - {v_2})\sup{3} + ({u_2} -
  {v_2})\sup{2}}, & J_{1011}=
  \span{i( {u_1} - {v_1})\sup{3} + ({u_1} - {v_1})\sup{2}},\cr
  J_{1101}=\span{i({u_3} - {v_3})\sup{3} + ({u_3} -
  {v_3})\sup{2}}, &J_{1111}=\span{i({u_3} +{v_3})\sup{3} +
  ({u_3} + {v_3})\sup{2}},\cr
  J_{0002}=\span{i 1\sup{3} +
    1\sup{2}, i E_2 - i E_3 -
    1\sup{1}}, &
    J_{1002}=\span{i({u_1} + {v_1})\sup{3} - ({u_1} +
    {v_1})\sup{2}},\cr
    J_{0102}=\span{i({u_2} + {v_2})\sup{3} - ({u_2} +
    {v_2})\sup{2}}, &
    J_{0012}=\span{-i({e_1} -{e_2})\sup{3} + ({e_1} -
    {e_2})\sup{2}},\cr
    J_{0112}=\span{i({u_2} - {v_2})\sup{3} - ({u_2} -
    {v_2})\sup{2}}, &
    J_{1012}=\span{-i({u_1} - {v_1})\sup{3} + ({u_1} -
    {v_1})\sup{2}},\cr
    J_{1102}=\span{i({u_3}-{v_3})\sup{3} - ({u_3} -
    {v_3})\sup{2}}, &
    J_{1112}=\span{i({u_3} + {v_3})\sup{3} - ({u_3} +
    {v_3})\sup{2}},
\end{eqnarray}
 } which is of type $(21,3)$.\medskip

 d) This is the $\Z_2^3\times\Z_4$-grading
{\small
 \begin{eqnarray}\label{nt4}
 J_{0000}=\span{E_1, E_2 + E_3}, & J_{1000}=\span{({u_1} +
 {v_1})\sup{1}},\cr
 J_{0100}=\span{({u_2} + {v_2})\sup{1}}, &
 J_{0010}=\span{({e_1} - {e_2})\sup{1}},\cr
 J_{0110}=\span{({u_2} - {v_2})\sup{1}}, &
 J_{1010}=\span{({u_1} - {v_1})\sup{1}},\cr
  J_{1100}=\span{({u_3} - {v_3})\sup{1}},&
  J_{1110}=\span{({u_3} + {v_3})\sup{1}},\cr
  J_{0001}=\span{-i 1\sup{3} + 1\sup{2}}, &
  J_{1001}=\span{i({u_1} +{v_1})\sup{3} + ({u_1} +
  {v_1})\sup{2}},\cr
  J_{0101}=\span{i ({u_2}+{v_2})\sup{3} + ({u_2} + {v_2})\sup{2}},
  & J_{0011}=\span{i ({e_1} - {e_2})\sup{3} + ({e_1} -
  {e_2})\sup{2}},\cr
  J_{0111}=\span{i ({u_2} - {v_2})\sup{3} + ({u_2} -
  {v_2})\sup{2}}, &
  J_{1011}=\span{{i ({u_1} -{v_1})\sup{3} +({u_1}-
  {v_1})\sup{2}}},\cr
  J_{1101}=\span{i ({u_3} - {v_3})\sup{3} + ({u_3} -
  {v_3})\sup{2}}, &
  J_{1111}=\span{i ({u_3}+{v_3})\sup{3} + ({u_3} +
  {v_3})\sup{2}},\cr
  J_{0002}=\span{E_3 - E_2, 1\sup{1}}, &
  J_{0003}=\span{i 1\sup{3} + 1\sup{2}},\cr
  J_{1003}=\span{i ({u_1} + {v_1})\sup{3} - ({u_1} +
  {v_1})\sup{2}}, &
  J_{0103}=\span{i ({u_2} + {v_2})\sup{3} - ({u_2} +
  {v_2})\sup{2}},\cr
  J_{0013}=\span{-i ({e_1} - {e_2})\sup{3} + ({e_1} -
  {e_2})\sup{2}}, &
  J_{0113}=\span{-i ({u_2} -{v_2})\sup{3} + ({u_2} -
  {v_2})\sup{2}},\cr
  J_{1013}=\span{-i ({u_1} -{v_1})\sup{3} + ({u_1} -
  {v_1})\sup{2}}, &
  J_{1103}=\span{-i ({u_3} -{v_3})\sup{3} + ({u_3} -
  {v_3})\sup{2}},\cr
  & J_{1113}=\span{i ({u_3} + {v_3})\sup{3} - ({u_3} +
  {v_3})\sup{2}},
 \end{eqnarray}}
which is of type $(23,2)$.\medskip

e) This is the $\Z_2^5$-grading given by {\small
 \begin{eqnarray}\label{nt5}
 J_{00000}=\span{E_1, E_2 + E_3, 1\sup{1}}, &
 J_{00010}=\span{-1\sup{3} + 1\sup{2}},\cr
 J_{00001}=\span{E_3 - E_2}, &
 J_{10010}=\span{({u_1}+{v_1})\sup{3} + ({u_1} +
 {v_1})\sup{2}},\cr
 J_{10001}=\span{({u_1} + {v_1})\sup{1}}, &
 J_{01010}=\span{({u_2}+{v_2})\sup{3} + ({u_2} +
 {v_2})\sup{2}},\cr
 J_{01001}=\span{({u_2} + {v_2})\sup{1}}, &
 J_{00110}=\span{({e_1} - {e_2})\sup{3} + ({e_1} -
 {e_2})\sup{2}},\cr
 J_{00101}=\span{(-{e_1} + {e_2})\sup{1}}, &
 J_{00011}=\span{1\sup{3} + 1\sup{2}},\cr
 J_{00111}=\span {({e_2} - \ {e_1})\sup{3} + ({e_1} -
 {e_2})\sup{2}}, & J_{01011}=\span {({u_2} + \ {v_2})\sup{3} -
 ({u_2} + {v_2})\sup{2}},\cr
 J_{01101}=\span{({u_2} -{v_2})\sup{1}}, &
 J_{01110}=\span{({u_2} - \ {v_2})\sup{3} + ({u_2} -
 {v_2})\sup{2}},\cr
 J_{01111}=\span {(-{u_2} + \ {v_2})\sup{3} + ({u_2} -
 {v_2})\sup{2}}, & J_{10011}=\span{({u_1} + \ {v_1})\sup{3} -
  ({u_1} + {v_1})\sup{2}},\cr
  J_{10101}=\span {({u_1} - \ {v_1})\sup{1}}, &
  J_{10110}=\span {({u_1} - \ {v_1})\sup{3} + ({u_1} -
  {v_1})\sup{2}},\cr
  J_{10111}=\span {(-{u_1} + \ {v_1})\sup{3} + ({u_1} -
  {v_1})\sup{2}}, &
  J_{11001}=\span {({u_3} - \ {v_3})\sup{1}},\cr
  J_{11010}=\span {({u_3} - \ {v_3})\sup{3} + ({u_3} -
  {v_3})\sup{2}}, &
  J_{11011}=\span{(-{u_3} + \ {v_3})\sup{3} + ({u_3} -
  {v_3})\sup{2}},\cr
  J_{11101}=\span {({u_3} + \ {v_3})\sup{1}}, &
  J_{11110}=\span {(-{u_3} - \ {v_3})\sup{3} - ({u_3} +
  {v_3})\sup{2}},\cr &
  J_{11111}=\span {({u_3} + \ {v_3})\sup{3} - ({u_3} +
  {v_3})\sup{2}},
 \end{eqnarray}}
 which is of type $(24,0,1)$.\medskip

 f) This is the $\Z_2^3$-grading (\ref{grad1}), which is nontoral of type
 $(0,0,7,0,0,1)$, as mentioned in Proposition~\ref{dece}.
\smallskip

 Before finishing this subsection we would like to exhibit another
 (nontoral) grading obtained as a coarsening of the one in case d) above, by
 removing $\hat f_0$ in the set of grading automorphisms. Thus consider
 the grading
 on the Albert algebra induced by the automorphisms
 $\{\widehat{t_{1,-1}},\widehat{t_{-1,1}},
\widetilde{\tau_{0,1}}\}$. This is a $\Z_2^2\times\Z_4$-grading
whose homogeneous component $J_{ijk}$ is just $J_{ij0k}\oplus
J_{ij1k}$ in (d), that is,
\begin{eqnarray}\label{coar}
 J_{000}=\span{E_1, E_2 + E_3, ({e_2} - {e_1})\sup{1}}, &&
 J_{001}=\span{-ie_1\sup{3} + e_2\sup{2},
 -i {e_2}\sup{3} + e_1\sup{2}},\cr
 J_{002}=\span{E_3 - E_2, 1\sup{1}}, &&
 J_{003}=\span{i e_1\sup{3} + e_2\sup{2},
  i e_2\sup{3} + e_1\sup{2}},\cr
  J_{010}=\span{u_2\sup{1}, v_2\sup{1}}, &&
  J_{011}=\span{-i u_2\sup{3} - u_2\sup{2},
  -i v_2\sup{3} - v_2\sup{2}},\cr
 J_{013}=\span{i u_2\sup{3} - u_2\sup{2},
  i v_2\sup{3} - v_2\sup{2}}, &&
  J_{100}=\span{u_1\sup{1}, v_1\sup{1}},\cr
  J_{101}=\span{-i u_1\sup{3} - u_1\sup{2},
  -i v_1\sup{3} - v_1\sup{2}}, &&
  J_{103}=\span{i u_1\sup{3} - u_1\sup{2},
  i v_1\sup{3} - v_1\sup{2}},\cr
  J_{110}=\span{u_3\sup{1}, v_3\sup{1}}, &&
  J_{111}=\span{-i u_3\sup{3} - u_3\sup{2},
  -i v_3\sup{3} - v_3\sup{2}},\cr
  J_{113}=\span{i u_3\sup{3} - u_3\sup{2},
  i v_3\sup{3} - v_3\sup{2}}. &&
\end{eqnarray}
As the   subalgebra of fixed elements of $\f4$ by $\{
 \Ad(\widehat{t_{1,-1}}), \Ad(\widehat{t_{-1,1}}), \Ad(\widetilde{\tau_{0,1}})
 \}$ has   rank $3$,   the induced
grading on $\f4$    is nontoral and our grading  on the Albert
algebra is also nontoral, of type $(0,12,1)$.

\subsection{Gradings from $ M_3(F)$}
Up to the moment we have detected seven equivalence classes of
nontoral gradings on $J$, all of them coming from the refinements
of the nontoral grading on $C$ by gradings on $\H_3(F)$. In order
to find new nontoral gradings on the Albert algebra we need to
look at $J$ from another  point of view, that is, we can use a
different model of $J$ which provides a new perspective. For
instance, Tits construction of $J$ recalled in Subsection
\ref{Tits}.

 Thus let us consider again the associative algebra
$A:= M_3(F)$ and the monomorphism $\iota\colon\aut(A)\to F_4$
such that $f\mapsto f^\bullet$ as described in \ref{Tits}. If
$\{f_i\}$ is a finite commutative family of semisimple
automorphisms of $A$, the same is true for the family
$\{f_i^\bullet\}$. Hence, for a $G$-grading on $A$ given by a
group homomorphism $\rho\colon \X(G)\to\aut(A)$ we immediately can
define the grading on $J$ given by $\iota\rho\colon \X(G)\to F_4$.
Consider now the $\Z_3^2$-grading on $A$ produced by the commuting
automorphisms $f:=\In(p)$ and $g:=\In(q)$ where
 {\def\U{A}
$p=\diag(1,\om,\om^2)$ being $\om$ a primitive cubic root of the
unit and
$$\tiny
q=\begin{pmatrix}0 & 1 & 0 \cr
  0 & 0 & 1 \cr
  1 & 0 & 0
  \end{pmatrix}.$$
 These automorphisms of $A$ are semisimple of order $3$. The group
 they generate, $\span{f,g}$, is usualy called Pauli group.
The simultaneous diagonalization of $A$ relative to $\{f,g\}$
yields $\U=\oplus_{i,j=0}^2\U_{i,j}$ for
{\small\begin{eqnarray}\nonumber \U_{00}=\span{1_A}, &
\U_{01}=\span{\om^2e_{11}-\om e_{22}+e_{33}}, & \U_{02}=\span{-\om
e_{11}+\om^2e_{22}+e_{33}},\cr
 \U_{10}=\span{e_{13}+e_{21}+e_{32}}, & \U_{11}=\span{\om^2e_{13}-\om
 e_{21}+e_{32}}, & \U_{12}=\span{-\om e_{13}+\om^2 e_{21}+e_{32}},\cr
 \U_{20}=\span{e_{12}+e_{23}+e_{31}}, & \U_{21}=\span{\om^2e_{12}-\om
 e_{23}+e_{31}}, & \U_{22}=\span{-\om e_{12}+\om^2e_{23}+e_{31}},
\end{eqnarray}}

\noindent where again  $e_{ij}$ denotes the elementary
$(i,j)$-matrix in $M_3(F)$. Thus we have a $\Z_3^2$-nontoral
grading on $A$ (since any maximal torus of $A$ fixes a frame of
idempotents and so any toral grading has a zero component of
dimension at least 3).
{\def\A{J}
 Next we can consider the grading induced on $J$ by $\{f^\bullet,
g^\bullet\}$. If we make a simultaneous diagonalization of $J$
relative to these automorphisms we get the $\Z_3^2$-grading
$\A=\oplus_{i,j=0}^2 \U_{i,j}^3$, which has $9$ summands of
dimension $3$ each one. This $\Z_3^2$-grading on $J$ is obviously
toral according to Lemma~\ref{Kasper}.   Let us consider a third
semisimple automorphism $\phi$ of order $3$ in the centralizer of
$\{f^\bullet, g^\bullet\}$. This will allow us to refine the
previous $\Z_3^2$-grading on $J$ to a $\Z_3^3$-grading. So
consider
  $\phi\in\aut(J)$
given by $\phi(a_0,a_1,a_2)=(a_0,\om a_1,\om^2 a_2)$ where $\om$
is as before a primitive cubic root of the unit. It is clear that
$\{f^\bullet,g^\bullet,\phi\}$ is a commutative set of semisimple
automorphisms of $J$. Making again a simultaneous diagonalization
of $J$ relative to $\{f^\bullet,g^\bullet,\phi\}$ we get
$J=\oplus_{i,j,k=0}^2J_{i,j,k}$ with
\begin{eqnarray}\nonumber
 \A_{i,j,0}=\U_{ij}\times 0\times 0,\cr
 \A_{i,j,1}=0\times\U_{ij}\times 0,\cr
 \A_{i,j,2}=0\times 0\times\U_{ij},
\end{eqnarray}
so that we have  $27$ one-dimensional homogeneous components. In
particular this $\Z_3^3$-grading on $J$ is fine  and nontoral
(otherwise $J_{0,0,0}$ would contain three orthogonal idempotents,
by Proposition~\ref{toynoto}). Consequently, the subgroup
$\span{f^\bullet,g^\bullet,\phi}$ of $F_4$ is maximal among the
abelian subgroups of $F_4$ whose elements are semisimple (MAD).
Observe also that the generators of the subspaces $A_{ij}$ of $A$
are invertible elements in $A$, hence taking into account
(\ref{rel}), the generators of the homogeneous components
$J_{i,j,k}$ are also invertible in $J$. Thus we have found a basis
of invertible homogeneous elements in the Albert algebra.

\smallskip

Since we are describing gradings on the Albert algebra in the
usual standard basis and this last $\Z_3^3$-grading has been given
in a different one, we are now giving the mentioned grading
relative to some standard basis. We take, for instance, the
grading:
\begin{eqnarray}\label{ztrescubo}
 J_{000}=\span{E_1+E_2+E_3}, && J_{001}=\span{\om E_1+\om^2 E_2+E_3},\cr
 J_{002}=\span{\om^2 E_1+\om E_2+E_3}, &&
 J_{010}=\span{u_3\sup{3}+e_1\sup{2}+v_3\sup{1}},\cr
 J_{011}=\span{\om^2u_3\sup{3}+\om e_1\sup{2}+v_3\sup{1}}, &&
 J_{012}=\span{\om u_3\sup{3}+\om^2e_1\sup{2}+v_3\sup{1}},\cr
 J_{020}=\span{v_3\sup{3}-e_2\sup{2}+u_3\sup{1}}, &&
 J_{021}=\span{\om^2v_3\sup{3}-\om e_2\sup{2}+u_3\sup{1}},\cr
 J_{022}=\span{\om v_3\sup{3}-\om^2 e_2\sup{2}+u_3\sup{1}}, &&
 J_{100}=\span{-v_2\sup{3}-u_2\sup{2}+e_1\sup{1}},\cr
 J_{101}=\span{-\om^2v_2\sup{3}-\om u_2\sup{2}+e_1\sup{1}}, &&
 J_{102}=\span{-\om v_2\sup{3}-\om^2 u_2\sup{2}+e_1\sup{1}},\cr
 J_{110}=\span{e_2\sup{3}-u_1\sup{2}+v_1\sup{2}}, &&
 J_{111}=\span{\om^2e_2\sup{3}-\om u_1\sup{2}+v_1\sup{1}},\cr
 J_{112}=\span{\om e_2\sup{3}-\om^2 u_1\sup{2}+v_1\sup{1}}, &&
 J_{120}=\span{v_1\sup{3}+v_3\sup{2}+v_2\sup{1}},\cr
 J_{121}=\span{\om^2 v_1\sup{3}+\om v_3\sup{2}+v_2\sup{1}}, &&
 J_{122}=\span{\om v_1\sup{3}+\om^2 v_3\sup{2}+v_2\sup{1}}, \cr
 J_{200}=\span{u_2\sup{3}+v_2\sup{2}+e_2\sup{1}}, &&
 J_{201}=\span{\om^2u_2\sup{3}+\om v_2\sup{2}+e_2\sup{1}},\cr
 J_{202}=\span{\om u_2\sup{3}+\om^2 v_2\sup{2}+e_2\sup{1}},&&
 J_{210}=\span{u_1\sup{3}+u_3\sup{2}+u_2\sup{1}},\cr
 J_{211}=\span{\om^2 u_1\sup{3}+\om u_3\sup{2}+u_2\sup{1}},&&
 J_{212}=\span{\om u_1\sup{3}+\om^2 u_3\sup{2}+u_2\sup{1}},\cr
 J_{220}=\span{-e_1\sup{3}-v_1\sup{2}+u_1\sup{1}},&&
 J_{221}=\span{-\om^2e_1\sup{3}-\om v_1\sup{2}+u_1\sup{1}},\cr
 J_{222}=\span{-\om e_1\sup{3}-\om^2 v_1\sup{2}+u_1\sup{1}}. &&
 \end{eqnarray}

\noindent It is produced by the  set of commuting
 diagonalizable automorphisms $\{t_{\om^2,\om^2,\om^2,1},$
 $t_{\om^2,\om,1,\om^2},\varphi\}$, where $\varphi=\theta\circ \Psi_U$
 for $U \in\mathop{\rm O}(C,n)  $  is
 given by the matrix relative to the standard basis
 of $C$
$$ \small{\begin{pmatrix}
 0 & 0 &
    -1 & 0 & 0 & 0 & 0 & 0 \cr 0 & 0 & 0 & 0 & 0 & 1 & 0 & 0 \cr 0 & 0 & 0 & 1 & 0 & 0 & 0 & 0 \cr
   0 & 1 & 0 & 0 & 0 & 0 & 0 & 0 \cr 0 & 0 & 0 & 0 & 0 & 0 & 0 & 1 \cr 0 & 0 & 0 & 0 & 0 & 0 & 1 &
   0 \cr -1 & 0 & 0 & 0 & 0 & 0 & 0 & 0 \cr 0 & 0 & 0 & 0 & 1 & 0 & 0 & 0
   \end{pmatrix}}$$
following the notations in 2.3. This $\Z_3^3$-grading is   fine
and nontoral,  since
 all the homogeneous components are one-dimensional.  Following
 \cite[(7.4) THEOREM, p.\,278]{Griess} for odd prime $p$ there is a unique conjugacy class
of elementary abelian nontoral $p$-subgroup. This is obtained for
$p=3$ and it is isomorphic to $\Z_3^3$ (see also \cite[TABLE II,
p.\,258]{Griess}).
  Therefore
the previous $\Z_3^3$-gradings on the Albert algebra are
isomorphic.

In   Section  6 we shall be able to give another description of
this grading in terms of the Weyl group of $F_4$. Besides, the
uniqueness of the $\Z_3^3$-grading will also be a consequence.
\smallskip

Once we have described the previous gradings   on the Albert
algebra, we can announce the first of our main results:

\begin{te}\label{notenJ}
The unique up to equivalence nontoral gradings on the Albert
algebra are those described in  (\ref{nt1}), (\ref{nt2}),
(\ref{nt3}), (\ref{nt4}), (\ref{nt5}), (\ref{grad1}), (\ref{coar})
and (\ref{ztrescubo}).

The nontoral fine gradings on $J$ are the ones described in
(\ref{nt1}), (\ref{nt5}) and (\ref{ztrescubo}).
\end{te}

In fact, there are four fine gradings, taking into account that
the Cartan decomposition is a toral and fine grading on $\f4$
which induces a toral and fine grading on $J$. The proof of the
above theorem will have to be postponed to a forthcoming section.

\section{Weyl group of $\f4$}

In   next sections we shall   use the Weyl group as an important
tool for our purposes. First of all we must invoke a version of
the  Borel-Serre  theorem (Theorem~\ref{Bose}) asserting that a
  supersolvable subgroup
  of semisimple elements in an algebraic group is contained
   in the normalizer of some maximal torus. In particular, this
   can be applied to finitely generated abelian groups. The point
   of this is that most of our arguments can be carried out within
   the normalizer of a maximal torus, hence the relevance of the
   Weyl group, which in our context is isomorphic to
   the quotient of the normalizer of any maximal torus by the
   torus itself.

In order to describe the abstract Weyl group of $\f4$, we must
begin by fixing a basis $\Delta=\{\a_i\mid i=1,\ldots,4\}$ of a
root system of $\f4$. Its Dynkin diagram is

\vskip .5cm \hskip 5cm \vbox{ \xy <1cm, 0cm>: \POS(0,0)*+{\circ}
\ar @{-} +(1,0)*+{\circ} \POS(1.1,0) \ar @2{-} +(1,0)*+{\circ}
\POS(1.5,0)*+{>} \POS(2.2,0)\ar @{-} +(1,0)*+{\circ}
\POS(3.3,0.3)*+{\a_4} \POS(1,.3)*+{ \a_2} \POS(2.1,.3)*+{\a_3}
\POS(0,.3)*+{\a_1}
\endxy}
\vskip .5cm \noindent and its Cartan matrix is
\begin{equation}\label{matriz}
\small\begin{pmatrix} 2 & -1 & 0 & 0\cr -1 & 2 & -2 & 0\cr 0
& -1 & 2 & -1\cr 0 & 0 & -1 & 2
\end{pmatrix}.
\end{equation}

 Taking the euclidean space
 $E=\sum_{i=1}^4\R\a_i$ with the inner product $\langle\ ,\ \rangle$,
 the Weyl group of $\f4$ is the subgroup $\W$ of $\GL(E)$
generated by the (simple) reflections $s_i$ with $i=1,2,3,4$,
given by $s_i(x):=x-\span{x,\a_i}\a_i$. Identifying $GL(E)$ to
$GL(4,\R)$ by means of the matrices relative to the $\R$-basis
$\Delta$, the reflections $s_i$ are represented by
$$
\small{\begin{array}{ll} s_1=\begin{pmatrix}-1 & 0 & 0 & 0\cr 1 &
1 & 0 & 0\cr 0 & 0 & 1 & 0\cr 0 & 0 & 0 & 1\end{pmatrix},\
&s_2=\begin{pmatrix}1 & 1 & 0 & 0\cr 0 & -1 & 0 & 0\cr 0 & 1 & 1 &
0\cr 0 & 0 & 0 & 1\end{pmatrix},\\\,\\s_3=\begin{pmatrix}1 & 0 & 0
& 0\cr 0 & 1 & 2 & 0\cr 0 & 0 & -1 & 0\cr 0 & 0 & 1 &
1\end{pmatrix},\ &s_4=\begin{pmatrix}1 & 0 & 0 & 0\cr 0 & 1 & 0 &
0\cr 0 & 0 & 1 & 1\cr 0 & 0 & 0 & -1\end{pmatrix},\end{array}}
$$
since the Cartan integers $\span{\a_i,\a_j}$ are the entries of
the Cartan matrix.

 We shall consider $\W\subset GL(4,\R)$ ordered
lexicographically, that is, first for any two different couples
$(i,j), (k,l)$ such that $i,j,k,l\in\{1,2,3,4\}$ we define
$(i,j)<(k,l)$ if and only if either $i<k$ or $i=k$ and $j<l$, and
second, for any two different matrices $\si=(\si_{ij})$,
$\si'=(\si'_{ij})$ in $\W$,  $\si<\si'$ if and only if
$\si_{ij}<\si'_{ij}$ where $(i,j)$ is the least element (with the
previous order in the couples)  such that $\si_{ij}\ne\si'_{ij}$.
One possible way to compute the Weyl group with this particular
enumeration is provided by the following code implemented with
{\sl Mathematica}:
\smallskip

{\parindent=3cm\tt
W=Table[$s_i$,\{i,4\}];

a[L\_,x\_]:=Union[L,

\hskip 2cm Table[L[[i]].x,\{i,Length[L]\}],

\hskip 2cm Table[x.L[[i]],\{i,Length[L]\}]]

Do[W=a[W,$s_i$],\{i,4\}]\hskip 1cm \textrm{(4 times repeated)} }
\smallskip

We get a list of $1152=2^73^2$ elements in the table {\tt W} which
is nothing but the Weyl group $\W$ of $\f4$. We are denoting by
$\si_i$ the $i$-th element of $\W$ lexicographically ordered. The
following result comes from a straightforward computation which
may be done with any matrix multiplication software.
\begin{pr}\label{sigmas}
The $1152$ elements of the Weyl group $\W$ of $\f4$ are distributed
in $25$ orbits (=conjugacy classes) according to the following table
{\small\rm
\begin{center}
\begin{tabular}{|c|c|c|c|}
\hline
order & no. of elements & no. of orbits & representatives \cr
\hline
 1 & 1 & 1 & $\si_{748}=1$\cr
2 & 139 & 7 & $\si_{28}$, $\si_{42}$, $\si_{55}$, $\si_{103},$\cr
  &     &   & $\si_{105}$, $\si_{142}$, $\si_{405}$\cr
3 & 80  & 3 & $\si_7$, $\si_{15}$, $\si_{114}$\cr 4 & 228 & 5 &
$\si_{1}$, $\si_{3}$, $\si_{56}$,\cr
  &     &   &  $\si_{104}$, $\si_{110}$\cr
6 & 464 & 7 &  $\si_4$, $\si_8$, $\si_9$, $\si_{14}$,\cr
  &     &   &  $\si_{30}$, $\si_{78}$, $\si_{106}$\cr
8 & 144 & 1 &  $\si_2$\cr
12 & 96 & 1 & $\si_{10}$\cr
\hline
   & 1152 & 25 & \cr
\hline
\end{tabular}
\end{center}
}
\end{pr}
The column in the left gives the   order of every element  in the
corresponding orbit. We shall denote by $I$ the set of indices of
representatives in the right column: {\small
\begin{equation}\label{ecjn}
1, 2, 3, 4, 7, 8, 9, 10,14, 15, 28, 30, 42, 55, 56, 78, 103,
104, 105,
106, 110, 114,142, 405, 748.
\end{equation}
}
\subsection{The maximal torus of $\aut(\f4)$}\label{ldlb}

If $\frak{h}$ is a Cartan subalgebra of $L=\f4$, consider $L=\frak
h\oplus(\oplus_{\alpha\in\frak h^*}L_\alpha)$ the decomposition in
root spaces relative to $\frak h$, that is, $L_\alpha=\{x\in L\mid
[h,x]=\alpha(h)x\,\forall h\in\frak h\}$ if $\a\in\frak h^*$,
$\Phi=\{\alpha\in\frak h^*\mid L_\alpha\ne0\}$ the root system,
and take a basis $\Delta=\{\alpha_1,\alpha_2,\alpha_3,\alpha_4\}$
of the root system. Identifying the roots to their coordinates
relative to the basis $\Delta$, the $24$ positive roots  of
$\Phi^+$ are:
\begin{equation}\label{mm1}
\begin{matrix}
 (0,0,0,1), & (0,1,1,1), & (1,2,2,1),\cr
 (0,0,1,0), & (0,1,2,0), & (1,1,2,2),\cr
 (0,1,0,0), & (1,1,1,1), & (1,2,3,1),\cr
 (1,0,0,0), & (0,1,2,1), & (1,2,2,2),\cr
 (0,0,1,1), & (1,1,2,0), & (1,2,3,2),\cr
 (0,1,1,0), & (1,1,2,1), & (1,2,4,2),\cr
 (1,1,0,0), & (0,1,2,2), & (1,3,4,2),\cr
 (1,1,1,0), & (1,2,2,0), & (2,3,4,2).
 \end{matrix}
 \end{equation}
 As usual, the nondegeneracy of the Killing form $k$ allows to
 identify $\frak{h}$ to $\frak{h}^*$, calling    $t_\a$ the
unique element in $\frak{h}$ satisfying
$\alpha(h)=k(t_{\alpha},h)$ for all $h\in\frak{h}$, as in
\cite[p.\,37]{Humphreysalg}.

Any automorphism   fixing pointwise $\frak h$ preserves the root
spaces. The set of all such automorphisms is a maximal torus of
$\aut(L)$; more precisely, given $x,y,z,u\in F^\times$ there is an
only automorphism $\Psi$ such that $\Psi\vert_{\frak h}=\id$,
$\Psi\vert_{L_{\alpha_1}}=x\,\id$,
$\Psi\vert_{L_{\alpha_2}}=y\,\id$,
$\Psi\vert_{L_{\alpha_3}}=z\,\id$,
$\Psi\vert_{L_{\alpha_4}}=u\,\id$ (particular case of the
isomorphism theorems in \cite[p.\,75]{Humphreysalg}). Obviously,
if $\alpha=n_1\alpha_1+n_2\alpha_2+n_3\alpha_3+n_4\alpha_4$, then
$\Psi\vert_{L_{\alpha}}=x^{n_1}y^{n_2} z^{n_3}u^{n_4}\id$. Denote
by $\Psi_{xyzu}$ the above automorphism $\Psi$, and by $ T_{\frak
h}$ the maximal torus $\{\Psi_{xyzu}\mid x,y,z,u\in F^\times\}$
($T_{\frak h}$ depends only on $\frak h$, and $\Psi_{xyzu}$
depends on $\frak h$ and $\Delta$).

On the other hand, we have got a concrete  maximal torus of $\F4
$, since we have got an algebraic group isomorphism $\Ad\colon
F_4\to\F4 $ and we have already introduced the maximal torus
$\To_0$ of $F_4$ (see (\ref{maxtor})).   Thus we get a maximal
torus $\To:=\Ad(\To_0)$ in $\F4 $, whose generic element is
$t'_{xyzu}:=\Ad( t_{xyzu})$. Let us take as $\frak h$ the Cartan
subalgebra of the elements fixed by $\To$. Since $\To$ is
contained in $T_{\frak h}$, they necessarily coincide.

 Let us choose a comfortable basis of $\f4$ for
which we know the matrix representation of $t'_{xyzu}$ relative to
it. Let $\omega_i$ ($i\in\{1,\dots,27\}$) be the $i$-th element in
$B$, the standard basis  on $J$ that we chose in Subsection 2.2.
Recall that $t_{xyzu}(\om_i)=\eta_i\om_i$ where
$\eta_i$($=\eta_i(x,y,z,u)$) is the $i$-th entry of the vector
$$ (1,1,1,x,\frac{1}{x},y,z,\frac{u^2}{xyz},\frac{1}{y},\frac{1}{z},
\frac{xyz}{u^2},u,\frac{1}{u},\frac{xy}{u},\frac{xz}{u},\qquad\qquad$$
$$\qquad\frac{u}{yz},
\frac{u}{xy},\frac{u}{xz},\frac{yz}{u},\frac{u}{x},\frac{x}{u},\frac{y}{u},
\frac{z}{u},\frac{u}{xyz},\frac{u}{y},\frac{u}{z},\frac{xyz}{u}).
$$
\noindent For any $v\in J$ define the map $R_v\colon J\to J$ such
that $a\mapsto av$. Since $\f4=\der(J)=[R_J,R_J]$
(\cite[p.\,117]{Schafer}),   we can extract a basis of $\f4$ from
the generators set $\{[R_{\om_i},R_{\om_j}]\}_{i,j= 1}^{27}$.
Taking into account that $$\Ad(t_{xyzu})[R_{\om_i},R_{\om_j}]=
[R_{t_{xyzu}(\omega_i)},R_{t_{xyzu}(\omega_j)}]= \eta_i \eta_j
[R_{\om_i},R_{\om_j}],$$ and  defining $S$ as the set of all
couples $(i,j)\in\{1,\ldots, 27\}^2$ such that
$[R_{\om_i},R_{\om_j}]\ne 0$, we have that the eigenvalues of
$t'_{xyzu}=\Ad(t_{xyzu})$ are those of the set $\{\eta_i\eta_j\mid
(i,j)\in S\}$, which are precisely:

\begin{eqnarray}\label{tpm}
\nonumber\big(1,1,1,1,x,\frac{1}{u},y\,z,\frac{u^2}{x\,y^2\,z},\frac{x}{u},
  \frac{y\,z}{u},\frac{u^2}{x\,y},\frac{u}{x\,y},\frac{x\,y\,z}{u},
  \frac{y\,z}{u^2},\frac{u}{y},\frac{x\,y\,z}{u^2}, \frac{1}{x\,y},
  \frac{1}{y},\frac{x^2\,y\,z}{u^2},\\\frac{z}{x},z,\frac{x}{y},
  \frac{z}{u},
  \nonumber x\,z,\frac{x\,z}{u},\frac{x\,z}{u^2},\frac{x\,y\,z^2}{u^2},
  \frac{z}{y},\frac{1}{x},u,\frac{1}{y\,z},\frac{x\,y^2\,z}{u^2},
  \frac{u}{x},\frac{u}{y\,z},
  \frac{x\,y}{u^2},\frac{x\,y}{u},
  \frac{u}{x\,y\,z},\\\frac{u^2}{y\,z},\frac{y}{u},\frac{u^2}{x\,y\,z},
  x\,y,y,\frac{u^2}{x^2\,y\,z},
  \frac{x}{z},\frac{1}{z},\frac{y}{x},
  \frac{u}{z},\frac{1}{x\,z},\frac{u}{x\,z},\frac{u^2}{x\,z},
  \frac{u^2}{x\,y\,z^2},\frac{y}{z}\big),
  \end{eqnarray}

  \noindent where each eigenvalue is repeated according to its
  multiplicity
(looking only at $S$, we would not know the multiplicities because
the set $\{\eta_i\eta_j\mid (i,j)\in S\}$ has 228 elements, but 1
must appear 4 times, and the remaining values at least once, so by
dimensions those are just the multiplicities). On the other hand,
recalling that $T_{\frak h}=\To$, there must exist rational
fractions  $X,Y,Z,U\in F(x,y,z,u)$ in the list (\ref{tpm}) such
that the whole list agrees with $(1,1,1,1,X^{n_1}Y^{n_2}
Z^{n_3}U^{n_4}\mid (n_1,n_2,n_3,n_4)\in\Phi)$, with
$\Phi=\Phi^+\cup(-\Phi^+)$ and $\Phi^+$ given by (\ref{mm1}). One
solution is, for instance
 \begin{equation}\label{elidlc}
X=\frac{u^2}{xy^2z},\quad Y=yz,\quad Z=\frac{1}{u},\quad U=x.
 \end{equation}

\noindent Next we choose as our reference basis of $\f4$ anyone
extracted of $\{[R_{\om_i},R_{\om_j}]\}_{(i,j)\in S}$
 such that the matrix of
$t'_{xyzu}$ relative to  this basis is diagonal with the list
(\ref{tpm}) as diagonal. One possible choice is {\small
$$\begin{array}{llll}
 b_1=[R_{\om_1}, R_{\om_4}]\quad &b_2=[R_{\om_1}, R_{\om_{13}}]\quad
&b_3=[R_{\om_6}, R_{\om_7}]\quad
  &b_4=[R_{\om_8}, R_{\om_9}]\\
b_5=[R_{\om_2}, R_{\om_{21}}]\quad &b_6=[R_{\om_1},
R_{\om_{19}}]\quad
 &b_7=[R_{\om_7}, R_{\om_8}]\quad
&b_8=[R_{\om_1}, R_{\om_{17}}]\\
  b_9=[R_{\om_2}, R_{\om_{27}}]\quad
 &b_{10}=[R_{\om_5}, R_{\om_{11}}]\quad
&b_{11}=[R_{\om_2}, R_{\om_{25}}]\quad &b_{12}=[R_{\om_1},
R_{\om_{11}}]\\
 b_{13}=[R_{\om_5}, R_{\om_{9}}]\quad
&b_{14}=[R_{\om_1}, R_{\om_{9}}]\quad &b_{15}=[R_{\om_4},
R_{\om_{11}}]\quad &b_{16}=[R_{\om_5}, R_{\om_{7}}]
\end{array}$$
$$\begin{array}{llll}
b_{17}=[R_{\om_1}, R_{\om_{7}}]\quad &b_{18}=[R_{\om_4},
R_{\om_{9}}]\quad &b_{19}=[R_{\om_2}, R_{\om_{23}}]\quad
&b_{20}=[R_{\om_4}, R_{\om_{7}}]\\
b_{21}=[R_{\om_1}, R_{\om_{15}}]\quad &b_{22}=[R_{\om_9},
R_{\om_{11}}]\quad &b_{23}=[R_{\om_7}, R_{\om_{11}}]\quad
&b_{24}=[R_{\om_7}, R_{\om_{9}}]\\
 b_{25}=[R_{\om_1}, R_{\om_5}]\quad &b_{26}=[R_{\om_1}, R_{\om_{12}}]\quad
&b_{27}=[R_{\om_9}, R_{\om_{10}}]\quad
  &b_{28}=[R_{\om_6}, R_{\om_{11}}]\\
 b_{29}=[R_{\om_2}, R_{\om_{20}}]\quad &b_{30}=[R_{\om_1}, R_{\om_{16}}]\quad
&b_{31}=[R_{\om_{10}}, R_{\om_{11}}]\quad
  &b_{32}=[R_{\om_1}, R_{\om_{14}}]\\
 b_{33}=[R_{\om_2}, R_{\om_{24}}]\quad &b_{34}=[R_{\om_4}, R_{\om_{8}}]\quad
&b_{35}=[R_{\om_{2}}, R_{\om_{22}}]\quad
  &b_{36}=[R_{\om_1}, R_{\om_{8}}]\\
 b_{37}=[R_{\om_4}, R_{\om_{6}}]\quad &b_{38}=[R_{\om_1}, R_{\om_{6}}]\quad
&b_{39}=[R_{\om_{5}}, R_{\om_{8}}]\quad
  &b_{40}=[R_{\om_4}, R_{\om_{10}}]\\
 b_{41}=[R_{\om_1}, R_{\om_{10}}]\quad &b_{42}=[R_{\om_5}, R_{\om_{6}}]\quad
&b_{43}=[R_{\om_{2}}, R_{\om_{26}}]\quad
  &b_{44}=[R_{\om_5}, R_{\om_{10}}]\\
 b_{45}=[R_{\om_1}, R_{\om_{18}}]\quad &b_{46}=[R_{\om_6}, R_{\om_{8}}],\quad
&b_{47}=[R_{\om_{8}}, R_{\om_{10}}]\quad
  &b_{48}=[R_{\om_6}, R_{\om_{10}}]
\end{array}$$}

\noindent where these are root vectors relative to $\frak h$, but
the missing elements spanning the Cartan subalgebra must be taken
carefully.
 If we denote
 $\beta_i\in \frak h^*$ for $i=1\dots 48$ such that $[h,b_i]=\beta_i(h)b_i$ for any $h\in\frak
 h$, it is
 easy to check   that $(\beta_4,\beta_3,\beta_2,\beta_1)$ is a basis of $\Phi$
 with Cartan matrix (\ref{matriz}),
 which will be our election for $\Delta$ from now on (and respectively for $\alpha_i$).
Thus, we can denote
\begin{equation}\label{gnrd}
\begin{array}{ll}
   b_1=v_{\alpha _4}   &   b_{25}=v_{-\alpha _4}     \\
   b_2=v_{\alpha _3}   &   b_{26}=v_{-\alpha _3}     \\
   b_3=v_{\alpha _2}   &   b_{27}=v_{-\alpha _2}    \\
   b_4=v_{\alpha _1}   &   b_{28}=v_{-\alpha _1}    \\
   b_5=v_{\alpha _3+\alpha _4}   &   b_{29}=v_{-\alpha _3-\alpha _4  }   \\
   b_6=v_{\alpha _2+\alpha _3}   &   b_{30}=v_{-\alpha _2-\alpha _3 }   \\
   b_7=v_{\alpha _1+\alpha _2}   &   b_{31}=v_{-\alpha _1-\alpha _2 }   \\
   b_8=v_{\alpha _1+\alpha _2+\alpha _3}   &   b_{32}=v_{-\alpha _1-\alpha _2-\alpha _3 }   \\
   b_9=v_{\alpha _2+\alpha _3+\alpha _4}   &   b_{33}=v_{-\alpha _2-\alpha _3-\alpha _4 }   \\
   b_{10}=v_{\alpha _2+2\alpha _3}   &   b_{34}=v_{-\alpha _2-2\alpha _3 }   \\
   b_{11}=v_{\alpha _1+\alpha _2+\alpha _3+\alpha _4}   &
     b_{35}=v_{-\alpha _1-\alpha _2-\alpha _3-\alpha _4 }   \\
   b_{12}=v_{\alpha _2+2\alpha _3+\alpha _4}   &   b_{36}=v_{-\alpha _2-2\alpha _3-\alpha _4 }   \\
   b_{13}=v_{\alpha _1+\alpha _2+2\alpha _3}   &   b_{37}=v_{-\alpha _1-\alpha _2-2\alpha _3 }   \\
   b_{14}=v_{\alpha _1+\alpha _2+2\alpha _3+\alpha _4}   &
     b_{38}=v_{-\alpha _1-\alpha _2-2\alpha _3-\alpha _4 }  \\
   b_{15}=v_{\alpha _2+2\alpha _3+2\alpha _4}   &   b_{39}=v_{-\alpha _2-2\alpha _3-2\alpha _4 }   \\
   b_{16}=v_{\alpha _1+2\alpha _2+2\alpha _3}   &   b_{40}=v_{-\alpha _1-2\alpha _2-2\alpha _3 }   \\
   b_{17}=v_{\alpha _1+2\alpha _2+2\alpha _3+\alpha _4}   &
     b_{41}=v_{-\alpha _1-2\alpha _2-2\alpha _3-\alpha _4 }   \\
   b_{18}=v_{\alpha _1+\alpha _2+2\alpha _3+2\alpha _4}   &
     b_{42}=v_{- \alpha _1-\alpha _2-2\alpha _3-2\alpha _4}   \\
   b_{19}=v_{\alpha _1+2\alpha _2+3\alpha _3+\alpha _4}   &
     b_{43}=v_{-\alpha _1-2\alpha _2-3\alpha _3-\alpha _4 }   \\
   b_{20}=v_{\alpha _1+2\alpha _2+2\alpha _3+2\alpha _4}   &
     b_{44}=v_{-\alpha _1-2\alpha _2-2\alpha _3-2\alpha _4 }   \\
   b_{21}=v_{\alpha _1+2\alpha _2+3\alpha _3+2\alpha _4}   &
     b_{45}=v_{-\alpha _1-2\alpha _2-3\alpha _3-2\alpha _4 }   \\
   b_{22}=v_{\alpha _1+2\alpha _2+4\alpha _3+2\alpha _4}   &
     b_{46}=v_{-\alpha _1-2\alpha _2-4\alpha _3-2\alpha _4 }   \\
   b_{23}=v_{\alpha _1+3\alpha _2+4\alpha _3+2\alpha _4}   &
     b_{47}=v_{-\alpha _1-3\alpha _2-4\alpha _3-2\alpha _4 }   \\
   b_{24}=v_{2\alpha _1+3\alpha _2+4\alpha _3+2\alpha _4}  \quad &
     b_{48}=v_{-2\alpha _1-3\alpha _2-4\alpha _3-2\alpha _4 }
\end{array}
\end{equation}
where each $v_\a$ is a root vector relative to the root $\a$,
  verifying
$$
\Ad(t_{xyzu})v_{m_1\a_1+m_2\a_2+m_3\a_3+m_4\a_4}=X^{m_1}Y^{m_2}Z^{m_3}U^{m_4}
v_{m_1\a_1+m_2\a_2+m_3\a_3+m_4\a_4}.
$$
At last   we choose our \emph{standard basis} of $\f4$ as
$$
B'=(4[b_4,b_{28}],4[b_{27},b_3],8[b_{26},b_2],8[b_{25},b_1],
b_i\mid i=1,\ldots,48),
$$
 formed by root vectors but besides described
with whole precision. This will be needed in the next section to
extend elements from $\W$ to $\F4 $. The first four elements in
$B'$, which of course form a basis of $\frak h$, have not been
chosen arbitrarily, but they are respectively $t_{\alpha_1}$,
$t_{\alpha_2}$, $t_{\alpha_3}$ and $t_{\alpha_4}$ for our election
of $\Delta$.

Notice that   we have algebraic group isomorphisms
$\alpha\colon(F^\times)^4\to \To$ such that
$\alpha(x,y,z,u)=t'_{xyzu}=\Ad(t_{xyzu})$ and
$\beta\colon(F^\times)^4\to \To$ acting as
$\beta(X,Y,Z,U)=\Psi_{XYZU}$ (now we have fixed $\Delta$). The map
$\eta\colon (x,y,z,u)\mapsto (X,Y,Z,U)$ (see (\ref{elidlc})) is an
automorphism of the algebraic group $(F^\times)^4$ and we have a
commutative diagram
 \vskip .5cm \hskip 5cm \vbox{
 \xy <1cm, 0cm>:
 \POS(0,0) $(F^\times)^4$\POS(-0.5,-0.1)\ar @{->} +(0,-0.9)*+
 {(F^\times)^4}
 \POS(0,-1)\ar @{->} +(1,0.5)*+{\To}
 \POS(0,0)\ar @{->} +(0.7,-0,4)
 \POS(0.4,0)*+{\a}
 \POS(0.4,-1.1)*+{\b}
 \POS(-0.8,-0,4)*+{\eta}
\endxy}
\vskip .5cm
 {\noindent in which all the arrows are isomorphisms.
 But $\Ad\colon \To_0\to \To$ is also an isomorphism, which allows to
 write the simultaneous diagonalization of $J$ relative to $\To_0$ as
 a $\Z^4$-grading in such a way that the induced grading on $\f4$ by $\Ad$
 is just the root decomposition indexed in the coordinates of
 $\Phi$ relative to $\Delta$, that is, $[R_{J_{n_1,n_2,n_3,n_4}},R_{J_{n_1',n_2',n_3',n_4'}}]\subset
L_{(n_1+n_1')\a_1+(n_2+n_2')\a_2+(n_3+n_3')\a_3+(n_4+n_4')\a_4}$.
Let us explain this in a practical way: Since the set of
eigenvalues of $t_{xyzu}$ is contained in the set of eigenvalues
of $t'_{xyzu}$, any eigenvalue of $t_{xyzu}$ can be written  in
the form $X^{m_1}Y^{m_2}Z^{m_3}U^{m_4}$. This defines a map  from
the set of eigenvalues of $t_{xyzu}$ to $\Z^4$ such that
$X^{m_1}Y^{m_2}Z^{m_3}U^{m_4}\mapsto (m_1,m_2,m_3,m_4)$, providing
a $\Z^4$-grading on the Albert algebra   $J=\oplus
J_{m_1,m_2,m_3,m_4}$ such that $(m_1,m_2,m_3,m_4)$ is in the image
of the above map.  To determine $J_{m_1,m_2,m_3,m_4}$ we find the
eigenvalue of $t_{xyzu}$ of the form
$X^{m_1}Y^{m_2}Z^{m_3}U^{m_4}$, and then take the element of the
standard basis of the Albert algebra which is an eigenvector for
that eigenvalue. For instance, to find $J_{1,2,2,1}$ we compute
$XY^2Z^2U=z$ and write $J_{1,2,2,1}=\span{u_2\sup{3}}$ since
$u_2\sup{3}$ is the basic element such that
$t_{xyzu}(u_2\sup{3})=z u_2\sup{3}$. The complete description of
the grading is:

\begin{eqnarray}\label{mm2}
&J_{0,0,0,0}=\span{E_1,E_2,E_3},&\cr
J_{0,0,0,1}=\span{e_1\sup{3}}, & J_{0,0,1,0}=\span{e_1\sup{2}}, &
J_{0,0,-1,-1}=\span{e_1\sup{1}}, \cr
 J_{0,0,0,-1}=\span{e_2\sup{3}}, & J_{0,0,-1,0}=\span{e_2\sup{2}},
 & J_{0,0,1,1}=\span{e_2\sup{1}},\cr
 J_{-1,-1,-2,-1}=\span{u_1\sup{3}}, &
 J_{-1,-1,-1,0}=\span{u_1\sup{2}}, &
 J_{-1,-1,-1,-1}=\span{u_1\sup{1}},\cr
 J_{1,2,2,1}=\span{u_2\sup{3}}, & J_{1,2,3,2}=\span{u_2\sup{2}}, &
 J_{1,2,3,1}=\span{u_2\sup{1}},\cr
 J_{0,-1,-2,-1}=\span{u_3\sup{3}}, &
 J_{0,-1,-1,0}=\span{u_3\sup{2}}, &
 J_{0,-1,-1,-1}=\span{u_3\sup{1}},\cr
 J_{1,1,2,1}=\span{v_1\sup{3}}, & J_{1,1,1,0}=\span{v_1\sup{2}}, &
 J_{1,1,1,1}=\span{v_1\sup{1}},\cr
 J_{-1,-2,-2,-1}=\span{v_2\sup{3}}, &
 J_{-1,-2,-3,-2}=\span{v_2\sup{2}},
 & J_{-1,-2,-3,-1}=\span{v_2\sup{1}},\cr
 J_{0,1,2,1}=\span{v_3\sup{3}}, & J_{0,1,1,0}=\span{v_3\sup{2}}, &
 J_{0,1,1,1}=\span{v_3\sup{1}},
\end{eqnarray}

\noindent which is of type $(24,0,1)$. We remark again that this
$\Z^4$-grading on the Albert algebra is toral and fine  and the
induced grading on $\f4$ is precisely the Cartan grading.
Alternatively we could have got this grading in Section 3 directly
from the maximal torus of $\aut (J)$ but in that case its
relationship to the Cartan grading would   have not been so
direct. The group $\Z^4$ is the universal
  group of the grading. Moreover all the toral gradings on
the Albert algebra are coarsenings of this   and so can be
obtained by constructing equivalence classes of epimorphisms
$\Z^4\to G$ module the relation  given as in \cite[4.1]{g2}. This
is another way of understanding Proposition~\ref{toynoto}.
\smallskip

To finish this subsection we must devote a few lines to the action
of the Weyl group $\W$ on the maximal torus of $\F4 $. Since this
is isomorphic to $(F^\times)^4$ we can introduce this action as
the integration of the action of $\W$ on the Cartan subalgebra
  $\W\times\h\to\h$ such that
$\sigma\cdot t_\a=t_{\sigma(\a)}$ for any $\sigma=(a_{ij})\in\W$
and $\a\in\Phi$.  Since $\sigma$ is an endomorphism of the dual
space $\h^*$, then the transposed matrix $\sigma^t$ represents the
dual map $\h\to\h$. So
 identifying the elements
in $\h$ with their coordinates relative to the basis
$(t_{\a_i})_{i=1}^4$,   the action of $\sigma$ on the element
$\sum_i x_i t_{\a_i}\in\h$ is given by {\small
$(x_1,x_2,x_3,x_4)\cdot  [(a_{ij})_{i,j=1}^4 ]^t$} where the
product $\cdot$ is the usual matrix product. The integration of
this, is the desired action $\W\times (F^\times)^4\to
(F^\times)^4$ which consequently acts in the form $\sigma\cdot
(X,Y,Z,U)=(X',Y',Z',U')$ where
\begin{eqnarray}\label{ecdlc}
 X'= & X^{a_{11}}Y^{a_{12}}Z^{a_{13}}U^{a_{14}}\cr
 Y'= & X^{a_{21}}Y^{a_{22}}Z^{a_{23}}U^{a_{24}}\cr
 Z'= & X^{a_{31}}Y^{a_{32}}Z^{a_{33}}U^{a_{34}}\cr
 U'= & X^{a_{41}}Y^{a_{42}}Z^{a_{43}}U^{a_{44}}.
 \end{eqnarray}

\subsection{Extending Weyl group elements to automorphisms of
$\f4$}\label{ewg} In this subsection we shall use the isomorphism
theorem of \cite[p.75]{Humphreysalg} for extending any
$\sigma\in\W$ to an automorphism $\tilde\sigma\in\F4 $. In the
context of the mentioned theorem we can take $L=L'=\f4$, $\h=\h'$
agreeing with the Cartan subalgebra generated by the four first
elements in $B'$ the standard basis of $\f4$, $\Phi=\Phi'$ the
root system relative to $\frak h$, we choose $\Delta=\Delta'$ the
basis of $\Phi$ as in 4.1 ($\alpha_1,\dots\alpha_4$ the roots
corresponding to $b_4,\dots b_1\in B'$ respectively) and finally,
we take as isomorphism $\h\to \h$ the induced by $\sigma$ (as
above, by means of the identification $\frak h\to\frak h^*$
through $t\mapsto k(t,-)$). According to that theorem, for any
choice $x_{\a_i}\in L_{\a_i}\setminus\{0\}$ and
$x'_{\sigma(\a_i)}\in L_{\sigma(\a_i)}\setminus\{0\}$ for
$i=1,2,3,4$, there is only one $\tilde\sigma\in\F4 $ such that
$\tilde\sigma(t_{\a_i})=t_{\sigma(\a_i)}$ and
$\tilde\sigma(x_{\a_i})=x'_{\sigma(\a_i)}$ for every $i=1,2,3,4$.
We choose $x_{\a_i}$ to be the generator $v_{\a_i}\in L_{\a_i}$ as
in (\ref{gnrd}) and also $x'_{\sigma(\a_i)}=v_{\sigma(\a_i)}\in
L_{\sigma(\a_i)}$. The matrix of $\tilde\sigma$ relative to the
standard basis is  block diagonal     $\tiny\begin{pmatrix} A &
0\cr 0 & D\end{pmatrix}  $ where $A$ is just the $4\times 4$
matrix of $\sigma$ relative to the basis $\{\a_i\}_{i=1}^4$ and
$D$ is a $48\times 48$ matrix with only one nonzero element in
each row and in each column.  Thus we have constructed an
injective map $\W\to\F4 $ such that $\sigma\mapsto\tilde\sigma$.
It is important to highlight that this is not a group homomorphism
but only a map. In fact, there does not exist a group monomorphism
$\W\to\F4 $, as it is proved in \cite[p.\,717]{raros}.

  We shall denote by  $\No$ the normalizer
of $\To$ in $\F4 $.  It
is a standard result 
  that $\W\cong\No/\To$. It follows easily,
by construction of $\tilde\sigma$, that $\tilde\si\in\No$ for any
$\si\in\W$. Thus the previous map $\W\to\F4 $ is actually a map
$\W\to\No$, and composing with the universal epimorphism
$\No\to\No/\To$ we get an injective map $\W\to \No/\To$ such that
$\sigma\mapsto \tilde\sigma\To$ (the equivalence class of
$\tilde\sigma$ in the quotient group). Since domain and codomain
of this map share the same finite cardinal, the map is a
bijection. Even more: it can be proved  that
$\widetilde{\si_1\si_2}$ is in the same equivalence class that
$\widetilde{\si_1}\widetilde{\si_2}$ (since
$\widetilde{\si_1\si_2} {\tilde\si_2}^{-1}{\tilde\si_1}^{-1}$ acts
in $\frak h$ as $ {\si_1\si_2}  {\si_2}^{-1} {\si_1}^{-1}=\id$, so
that it belongs to $\To$) which proves that the previous map
$\W\to\No/\To$ is a group isomorphism. In particular
$$\No=\{\tilde\si t\mid \si\in\W, t\in\To\} .$$

We can now revisit the action of the Weyl group $\W$ on the
maximal torus $(F^\times)^4\cong\To$ from another viewpoint.
Identifying $\W$ with $\No/\To$ we can define the action
$\W\times\To\to\To$ given  by $\sigma\cdot t:=\tilde\si
t\tilde\si^{-1}$ for $\si\in\W$ and $t\in\To$. Then the
isomorphism $\beta\colon(F^\times)^4\to\To$ given by
$\beta(X,Y,Z,U)=\Psi_{XYZU}$ is an isomorphism of $\W$-groups in
the sense that $\beta(\si\cdot t)=\si\cdot\beta(t)$. Thus
$\si\cdot \Psi_{XYZU}=\Psi_{X'Y'Z'U'}$ as in (\ref{ecdlc}). And
since $ \Psi_{XYZU}=t'_{xyzu}$ for (\ref{elidlc}), a simple
computation proves that $\si\cdot t'_{xyzu}=t'_{x'y'z'u'}$ where
now {\small\begin{eqnarray}\label{ecdlc1}
\begin{matrix}
 x'= & x^{b_{11}}y^{b_{12}}z^{b_{13}}u^{b_{14}}\cr
 y'= & x^{b_{21}}y^{b_{22}}z^{b_{23}}u^{b_{24}}\cr
 z'= & x^{b_{31}}y^{b_{32}}z^{b_{33}}u^{b_{34}}\cr
 u'= & x^{b_{41}}y^{b_{42}}z^{b_{43}}u^{b_{44}}
 \end{matrix}\,,\quad\hbox{with }
 (b_{ij})=
 m\sigma m^{-1},\quad
 m=\begin{pmatrix}
 0 & 0 & 0 & 1\cr
 -1 & -1 & -2 & -1\cr
 1 & 2 & 2 & 1\cr
 0 & 0 & -1 & 0
 \end{pmatrix}.\cr
 \end{eqnarray}}
 The action of $\W$ on $(F^\times)^4$ given by
 $\si\cdot(x,y,z,u)=(x',y',z',u')$ as above is essential for our
 study, specially   the study of fixed elements in the torus under the action
 of certain elements of $\W$. Denote
 $$
   \tor{j}=\{t\in\To\mid\si_j\cdot t=t\}.
   $$
     It is easily seen that this
is a subgroup of $\To$ such that $\tor{i}\cong\tor{j}$ when
$\si_i$ is conjugated to $\si_j$ in $\W$. The information given by
these subgroups $\tor{i}$ is needed for our study, so we are
calculating them. For this, it suffices to consider the
representatives of conjugacy classes given in the table of
Proposition~\ref{sigmas}. We summarize all this information in the
following table. In it, $\tor{j}$ is the subgroup of all
$t'_{xyzu}\in\To$ such that the element given satisfies the
displayed condition. We also write down the abstract group
isomorphic to $\tor{j}$  in the right column.\smallskip
\begin{center}
\begin{tabular}{|c|c|c|c|}
\hline $j$ & Generic element  of $\tor{j}$& Membership condition &
Isomorphic to\cr \hline
 $1$ & $(u^2y,y,y,u)$ & $y^2=1, u\in F^\times$ & $F^\times\times\Z_2$\cr
 $2$ & $(x,x,x,1)$ & $x^2=1$ & $\Z_2$\cr
 $3$ & $(u^2y,y,u^2,u)$ & $u^4=1=y^2$ & $\Z_4\times\Z_2$\cr
 $4$ & $(x,x^{-1},1,x^{-1})$ & $x\in F^\times$ & $F^\times$\cr
 $7$ & $(u^2y^{-3},y,y,u)$ & $u,y\in F^\times$ & $(F^\times)^2$\cr
 $8$ & $(x,x,x,x^2)$ & $x\in F^\times$ &  $F^\times$\cr
 $9$ & $(1,u^{2/3},u^{2/3},u)$ & $u\in F^\times$ & $F^\times$\cr
 $10$ & $(1,1,1,1)$ &  & $\{1\}$\cr
 $14$ & $(x,x,x,y)$ & $x^2=y^2=1$ & $\Z_2^2$\cr
 $15$ & $(x,y,x^2y,x^2)$ & $x^3=y^3=1$ & $\Z_3^2$\cr
 $28$ & $(x,x^{-1},z,x^{-1})$ & $x,z\in F^\times$ & $(F^\times)^2$\cr
 $30$ & $(1,1,z,1)$ & $z\in F^\times$ & $F^\times$\cr
 $42$ & $(x,x,z,u)$ & $x^2=u^2=1,z\in F^\times$ & $F^\times\times\Z_2^2$\cr
 $55$ & $(x,y,y,u)$ & $x,y,u\in F^\times$ & $(F^\times)^3$\cr
 $56$ & $(x,y,y,xy)$ & $x,y\in F^\times$ & $({F^\times})^2$\cr
 $78$ & $(1,1,1,1)$ &  & $\{1\}$\cr
 $103$ & $(x,y,y,u)$ & $y^2=1, x,u\in F^\times$  & $(F^\times)^2\times\Z_2$\cr
 $104$ & $(x,y,y,xy)$ & $y^2=1, x\in F^\times$  & $F^\times\times\Z_2$\cr
 $105$ & $(x,y,xy,u)$ & $x^2=y^2=1,u\in F^\times$ &
$F^\times\times\Z_2^2$\cr
 $106$ & $(x,y,xy,y)$ & $x^2=y^2=1$ & $\Z_2^2$\cr
 $110$ & $(x,x,x,u)$ & $x^2=u^2=1$  & $\Z_2^2$\cr
 $114$ & $(x,1,z,1)$ & $x,z\in F^\times$  & $(F^\times)^2$\cr
 $142$ & $(x,1,z,u)$ & $x,z,u\in F^\times$  & $(F^\times)^3$\cr
 $405$ & $(x,y,z,u)$ & $x^2=y^2=z^2=u^2=1$ & $\Z_2^4$\cr
 $748$ & $(x,y,z,u)$ & $x,y,z,u\in F^\times$  & $(F^\times)^4$\cr
 \hline
\end{tabular}
\end{center}
\begin{equation}\label{tablon}
\end{equation}

\section{Quasitori in $F_4$}\label{seqt}

Recall that a quasitorus is a commutative algebraic group whose
identity component is a torus \cite[p.\,105]{enci}. An algebraic
linear group is a quasitorus if and only if in some basis its
elements can be expressed simultaneously by diagonal matrices.
Such groups are also called \emph{diagonalizable}. Besides, a
quasitorus $Q$ in an algebraic group $G$ can be written as a
disjoint union $Q={T}\cup {T}a_1\cup\cdots\cup{T}a_k$ where ${T}$
is a torus and $\{1_G,a_1,\ldots,a_k\}$ a finite abelian subgroup
of $G$. We remark that, as  a consequence of the algebraic version
of the  Borel-Serre theorem (Theorem~\ref{Bose} in the Appendix),
any quasitorus in $G$ normalizes some of the maximal tori of $G$.
To see that, define $Z$ as the centralizer in $G$ of ${T}$.
Applying this theorem to $H :=\{1_G,a_1,\ldots,a_k\}$, which is
contained in $Z$,   there is $T'$   some maximal torus of $Z$ that
$H\subset\No_Z(T')$ (the normalizer of $T'$ in $Z$). But ${T}$ is
contained in the center of $Z$ and since all its elements are
semisimple, ${T}$ is contained in the intersection of all maximal
tori of $Z$, hence in $T'$. Thus ${T}\subset T'\subset\No_Z(T')$
and since we had $H\subset\No_Z(T')\subset\No_{G}(T')$ then
$Q\subset \No_{G}(T')$. But   actually $T'$ is a maximal torus in
$G$, because if $T''$ is a maximal torus of $G$ which contains
$T'\supset T$, then $T''\subset Z$ and $T''=T'$.

  As any grading is given by a
 quasitorus,  the above paragraph gives
 the reason why we want to work inside $\No =\W\, \To$
   and we
 have studied in detail the Weyl group and
 its action on the torus $\To$ in the previous section.

 Next we consider a class
of quasitori which is relevant for our study. Define for each
$j\in\{1,\ldots,1152\}$ and each $t\in\To$ the quasitorus $A(j,t)$
as the (closed) subgroup of $\F4 $ generated by $\tor{j}$ and
$\widetilde{\si}_jt$, which of course defines a grading on $\f4$
by the group $\X(A(j,t))$ as in 2.1. But it suffices to consider
the gradings induced by the quasitori $A(j,t)$ with $j\in I$ (the
set $I$ defined in (\ref{ecjn})), taking into account  that if
$\si_i$ and $\si_j$ are conjugated in $\W$, then $A(i,t)\cong
A(j,t')$ for a suitable $t'\in\To$.
 We also have
\begin{pr}\label{ehlc}
If for some $j$ the group $A(j,\id)$ is toral then $A(j,t)$ is
toral for any $t\in\To$.
\end{pr}
Proof. Let $Z=\Co_{\Fs4}(\tor{j})$ and $Z_0$ its unit component.
Since $A(j,\id)$ is toral there is some maximal torus $T$ of $\F4
$ such that $A(j,\id)\subset T$. Then $T\subset Z$ but from
$\tor{j}\subset\To$ we also get $\To\subset Z$. Of course
$\To,T\subset Z_0$ and since $t\in\To$ and $\widetilde{\si}_j\in
A(j,\id)\subset T$ we have $\widetilde{\si}_j,t\in Z_0$ hence
$\widetilde{\si}_j t\in Z_0$. But $\widetilde{\si}_j t$ is a
semisimple element of $Z_0$ and consequently there is some $p\in
Z_0$ such that $p\widetilde{\si}_j t p^{-1}\in\To$. This equality
together with the fact that $p\tor{j}p^{-1}=\tor{j}$ imply that
$pA(j,t)p^{-1}\subset\To$. $\square$
\smallskip

The same proof shows that $A(j,t)$ is toral if and only if
$A(j,t')$ is toral for all $t'\in\To$. Now let us detect the
indices which cause nontorality.

\begin{pr}\label{nonto}
For $j\in I$ the group $A(j,\id)$ is nontoral if and only if $j=3,15,105,106$
or $405$.
\end{pr}
Proof. Let us prove first that the five quasitori are nontoral.
For $A(15,\id)$ we have $\tor{15}\cong\Z_3^2$ (see Table of
Section \ref{ewg}). The grading induced by this quasitorus is
produced by the automorphisms
$\{t'_{\om,\om,1,\om^2},t'_{1,\om,\om,1},\widetilde{\si_{15}}\}$
where $\om$ is a primitive cubic root of $1$. This is a
$\Z_3^3$-grading and computing the subalgebra of fixed elements by
the three previous automorphisms we find that this is null. This
implies that the grading is nontoral since in the toral case, this
should be an algebra of rank four. For $A(405,\id)$ we have
$\tor{405}\cong\Z_2^4$ and the associated grading agrees with the
one produced by
$\{t'_{-1,1,1,1},t'_{1,-1,1,1},t'_{1,1,-1,1},t'_{1,1,1,-1},
\widetilde{\si_{405}}\}$. This is a $\Z_2^5$-grading whose
$0$-homogeneous component is again null, hence the grading is
nontoral. For $A(3,\id)$ we have $\tor{3}\cong\Z_4\times\Z_2$ and
the induced grading is the produced by
$\{t'_{1,-1,-1,i},t'_{-1,-1,1,1},\widetilde{\si_3}\}$ which is a
$\Z_2\times\Z_2\times\Z_8$-grading (see remarks 1 and 2 after this
proof) whose $0$-homogeneous component has dimension $1$. Hence
$A(3,\id)$ is nontoral. The grading induced by $A(106,\id)$ is
also nontoral since $\tor{106}\cong \Z_2^2$ and the grading is the
one induced by $\{t'_{-1,1,-1,1},
t'_{1,-1,-1,-1},\widetilde{\si_{106}}\}$, which is a
$\Z_2^2\times\Z_6$-grading whose $0$-homogeneous component is
one-dimensional. The last grading is the induced by $A(105,\id)$.
We have $\tor{105}\cong F^\times\times\Z_2^2$, so
$A(105,\id)=\langle\tor{105},\widetilde\si_{105}\rangle\cong
F^\times\times\Z_2^3$ which induces a grading over
$\X(F^\times\times\Z_2^3)\cong \Z \times\Z_2^3$. The grading
agrees with the one produced for instance by
$\{\widetilde{\si_{105}},t'_{x,y,xy,1}, t'_{1,1,1,2}\mid
x^2=y^2=1$\}, which can be implemented in a computer.  This is a
$\Z_2^3\times\Z$-grading whose $0$-homogeneous component is
one-dimensional and so $A(105,\id)$ is nontoral. We include also a
table of homogeneous components dimensions for further reference.
These types can be computed with any linear algebra software
allowing simultaneous diagonalization.
\begin{center}
\begin{tabular}{|c|c|}
 \hline
Quasitorus  & Type\cr
 \hline
$A(3,\id)$ & $( 19,6,7)$\cr $A(15,\id)$ & $( 0,26)$\cr
$A(105,\id)$ & $( 31,0,7)$\cr $A(106,\id)$ & $( 3,14,7)$\cr
$A(405,\id)$ & $( 24,0,0,7)$\cr \hline
\end{tabular}
\end{center}
Let us prove now that $A(j,\id)$ is toral in the rest of the
cases. If $\tor{j}=\id$ then $A(j,\id)$ is cyclic and then toral
(this applies to the cases $j=10,78$). In case $\tor{j}$ is cyclic
or $\tor{j}\cong F^\times$, then $A(j,\id)$ has two factors and by
Lemma~\ref{Kasper} (\cite[Lema 1.1.3, p.\,5]{Kasper}) the grading
is toral (this applies to $j=2,4,8,9,30$). Another trivial case is
$j=748$ since $\widetilde{ \si_{748}}=\id$ and $\tor{748}=\To$.
For $j=1,7,14,28,42,55,56,103,104,114$ and $142$, performing a
simultaneous diagonalization of the algebra relative to the set of
automorphisms inducing the grading, one finds that the zero
homogeneous component of the corresponding grading is an abelian
four-dimensional algebra. Thus the grading is toral. Finally, for
$j=110$ we have $A(110,\id)=\span{t'_{-1,-1,-1,1},
t'_{1,1,1,-1},\widetilde{\si_{110}}}$, which produces a
$\Z_2^2\times\Z_4$-grading. In this case the  zero homogeneous
component is a six-dimensional (reductive) algebra
$L_e=\span{y_1,\dots,y_6}$ where
$$\begin{array}{ll}
y_1=b_{3}-b_{10}+b_{27}+b_{34}, \qquad&
y_{4}=-b_{4}+b_{22}+b_{28}+b_{46}, \\
y_{2}=-b_{7}-b_{18}-b_{31}+b_{42},&
y_{5}=b_{13}-b_{23}-b_{37}+b_{47},\\
y_{3}=b_{16}-b_{20}-b_{40}+b_{44}, &
y_{6}=b_{15}-b_{24}-b_{39}+b_{48},
\end{array}$$ which has rank $4$  because $\{y_1-y_6,y_2-y_5,y_3-y_4\}$ is contained
in the center of $L_e$ (there are only two types of
six-dimensional reductive subalgebras, $\frak a_1$ plus a
three-dimensional center and $2\frak a_1$, of ranks 4 and 2
respectively). Hence the grading is toral. $\square$\smallskip

\noindent\textbf{Remark 1.} Notice that the order of $\si_i$ does
not necessarily coincide with the order of $\widetilde\si_i$, but
it is a divisor. That happens, for instance, for $i=3$, since
$\si_3$ has order 4 while $\widetilde\si_3$ has order 8. This is
not because of a bad choice of $\widetilde\si_3$, since all the
possible extensions of $\si_3$ have the same order, as the next
lemma shows.
\medskip

\begin{lm}
Take $j\in\{1,\dots,1152\}$, and $m$  the order of $\si_j\in\W$.
Then the following conditions are equivalent:
\begin{itemize}
\item[i)] $\tor{j}$ is finite,
\item[ii)] $\tor{j}\subset\{t'_{x,y,z,u}\mid x^m=y^m=z^m=u^m=1\}$,
\item[iii)] $(\widetilde\si_jt)^m=\widetilde\si_j^m$ for any
$t\in\To$,
\item[iv)] Every element in $\{f\in\No\mid \pi(f)=\si_j\}$
($\pi\colon\No\to\W$ the canonical projection) has the same order.
\end{itemize}
\end{lm}

Proof.  Take the element in $\To$ given by
$$s_{x,y,z,u}:=t'_{x,y,z,u}(\widetilde\si_jt'_{x,y,z,u}\widetilde
\si_j^{-1})(\widetilde\si_j^2t'_{x,y,z,u}\widetilde\si_j^{-2})
\dots(\widetilde\si_j^{m-1}t'_{x,y,z,u}\widetilde\si_j^{1-m})$$
(product of elements in $\To$). Since $\si_j^m=\id$, we have
$\widetilde\si_j^m\in\To$ and thus
$\widetilde\si_js_{x,y,z,u}\widetilde \si_j^{-1}=s_{x,y,z,u}$,
that is, $s_{x,y,z,u}\in\tor{j}$. Besides it verifies that
$(\widetilde\si_jt'_{x,y,z,u})^m=s_{x,y,z,u}\widetilde\si_j^m$.

The implication $ii)\Rightarrow i)$ is trivial. Now, if we assume
$iii)$, $s_{x,y,z,u}=\id$ for any $x,y,z,u\in F^\times$. But if
$t'_{x,y,z,u}\in\tor{j}$, then $(t'_{x,y,z,u})^m=s_{x,y,z,u}$, and
so $t^m=\id$ for any $t\in\tor{j}$, and we have $ii)$.

Next suppose $i)$. We have
$s_{x,y,z,u}=t'_{f_1(x,y,z,u),f_2(x,y,z,u),f_3(x,y,z,u),f_4(x,y,z,u)}$
for some    rational  fractions $f_i\in F(x,y,z,u)$. But the set
$\{s_{x,y,z,u}\mid x,y,z,u\in F^\times\}\subset\tor{j}$ must be
finite, so the fractions are constant $s_{x,y,z,u}=t'_{a,b,c,d}$,
and, since $s_{1,1,1,1}=\id$, the identity $s_{x,y,z,u}=\id$ holds
for any $x,y,z,u\in F^\times$. This gives $iii)$.

From $iii)$, we conclude that the orders of $\widetilde\si_jt$ and
$\widetilde\si_j$ are the same, since none of them can be less
than $m$ and $(\widetilde\si_jt)^m=\widetilde\si_j^m$, so we have
$iv)$.

Conversely, let $c$ be such that $mc$ is the order of
$\widetilde\si_jt'_{x,y,z,u}$ for any $x,y,z,u$. Thus
$\id=(\widetilde\si_jt'_{x,y,z,u})^{mc}=(s_{x,y,z,u}{\widetilde\si_j}^m)^c=
s_{x,y,z,u}^c$, hence   $f_i(x,y,z,u)^c=1$ and so
$s_{x,y,z,u}=\id$, which is equivalent to $iii)$.
$\square$\smallskip

\noindent\textbf{Remark 2.} Apparently  the grading induced by the
quasitorus  $A(3,\id)=\langle\{\widetilde{\si_3},$
$t'_{1,-1,-1,i}, t'_{-1,-1,1,1} \}\rangle$  is a
$\Z_8\times\Z_4\times\Z_2$-grading, since 8, 4 and 2 are
respectively the orders of the generators. However,   the right
group generated by the support is $\Z_2^2\times\Z_8$ because
$(\widetilde\si_3)^2t'_{1,-1,-1,i}$ has order 2.
\medskip

Returning to our quasitori $A(j,t)$, the toral element $t$ plays a
secondary roll.

\begin{pr}\label{retra1}
If $j\in \{3,15,105,106,405\}$, the quasitorus $A(j,t)$ is
conjugated to $A(j,t')$ for any $t,t'\in\To$.
\end{pr}
Proof. Take
$$
\sor{j}=\{\widetilde{\si}_j^{-1}t\widetilde{\si}_jt^{-1}\mid
t\in\To\}.
$$
Denoting by $\si=\widetilde{\si}_j$, we have  that
$(\sigma^{-1}t\sigma t^{-1})(\sigma^{-1}s\sigma
s^{-1})=\sigma^{-1}ts\sigma t^{-1}\sigma^{-1}\sigma
s^{-1}=\sigma^{-1} ts \sigma(ts)^{-1}$, since $\sigma
t^{-1}\sigma^{-1}\in\To$ and so it commutes with $s$. Thus
$\sor{j}$ is a subgroup. Besides it has the property that
$A(j,\id)$ is conjugated to $A(j,s)$ for any $s\in\sor{j}$.
Indeed,   if $s=\sigma^{-1}t\sigma t^{-1}$, then $\Ad(
t)(\tor{j})=\tor{j}$ and $\Ad( t)(\si)=\si s$.

On the other hand, it is clear that $A(j,s)=A(j,st)$ for any
$t\in\tor{j}$, therefore
$$
A(j,\id)\cong A(j,st)
$$
for any $s\in\sor{j},t\in\tor{j}$. Let us see that
$\tor{j}\sor{j}=\To$. First we observe that the map
$$\begin{array}{lll}
\To/\tor{j}&\to&\sor{j}\\
\protect[t]&\mapsto&\widetilde{\si}_j^{-1}t\widetilde{\si}_jt^{-1}
\end{array}
$$
is a group isomorphism. It is well defined and injective because
  $t\in\tor{j}$ if and only if
$\widetilde{\si}_j^{-1}t\widetilde{\si}_jt^{-1}=\id$. In
particular $\dim\To=\dim\tor{j}+\dim\sor{j}$ (see
\cite[Proposition~2.26, p.\,41]{Milne}). And we have another
isomorphism:
$$
 \sor{j}/\sor{j}\cap\tor{j}\to\sor{j}\tor{j}/\tor{j},
$$
hence
$\dim\tor{j}+\dim\sor{j}=\dim\sor{j}\tor{j}+\dim\sor{j}\cap\tor{j}$.

But $\sor{j}\cap\tor{j}$ is a finite group for  $j$ any of our
indices. Indeed, for $j=3,15,106,405$ the group $\tor{j}$ is
already finite, and for $j=105$, $\sor{105}=\{t'_{x,y,z,u}\mid
u^2=xyz\}$, $\tor{105}=\{t'_{x,y,z,u}\mid x^2=y^2=1,z=xy,u\in
F^\times\}$ and $\sor{105}\cap\tor{105}=\{t'_{x,y,xy,uxy}\mid
x^2=y^2=u^2=1\}\cong\Z_2^3$.

Consequently $\dim\To=\dim\tor{j}+\dim\sor{j}=\dim\sor{j}\tor{j}$
and so $\tor{j} \sor{j}= \To$.
 $\square$
\medskip

Furthermore, $\sor{j}\cap\tor{j}$ is a finite group for all $j\in
I$, so also for any $j\in\{1,\dots1152\}$, and thus $ A(j,t)\cong
A(j,t') $  for all $t,t'\in\To$, although it is unnecessary for
our purposes.

As a consequence of this lemma, if $\tor{j}$ is finite, then
$\sor{j}=\To$, and every $t\in\To$ is in $\sor{j}$, that is, there
is $s\in\To$ such that
$t=\widetilde{\si}_j^{-1}s^{-1}\widetilde{\si}_js$. Thus we have
obtained the following technical result, which will be very useful
in some forthcoming proofs.

\begin{co}\label{coento}
If $\tor{j}$ is finite, for any $t\in\To$ there is $s\in\To$ such
that
$$
s\widetilde{\si}_jts^{-1}=\widetilde{\si}_j
$$
\end{co}

The relevance  of the quasitori $A(j,t)$ is highlighted by the
following result.
\begin{pr}\label{pade1}
Let $F=\{f_1,\ldots,f_n,f_{n+1}\}$ be a nontoral commutative
family of semisimple elements in a connected algebraic group $G$
such that $\{f_1,\ldots,f_n\}$ is toral. Fix $T$   any maximal
torus of $G$. Then,   the subgroup generated by $F$ is conjugated
to some subgroup of the form $\span{t_1,\ldots,t_n,\si}$ where
$t_i\in T$ and $\si$ is a conjugated of $f_{n+1}$ in the
normalizer of $T$ in $G$.
\end{pr}
Proof. Define $Z$ as the centralizer of $\{f_1,\ldots,f_n\}$ in
$G$ and let $T'$ be some maximal torus in $G$ containing
$\{f_1,\ldots,f_n\}$. The subgroup $\overline{\span{F}}\subset Z$
  is
 a quasitorus of $Z$, hence is contained
in the normalizer $\No_Z(T'')$ of some maximal torus $T''$ in $Z$.
Then $T''$ is also a maximal torus in $G$ (since $T'\subset Z$)
and there is some $p\in G$ such that $T''=pTp^{-1}$. On the other
hand $\{f_1,\ldots,f_n\}$ is contained in the center of $Z$, and
since these are semisimple elements, then they are contained in
each maximal torus of $Z$. In particular $f_i\in T''$ for
$i\in\{1,\ldots,n\}$, and $f_{n+1}\in
\No_Z(T'')\subset\No_G(T'')$. Thus, $p^{-1}\span{F}p$ is generated
by    $f'_i=p^{-1}f_ip$ for $i=1,\ldots,n+1$, with $f'_i\in T$ for
$i\le n$ and $f'_{n+1}\in p^{-1}\No_G(T'')p=\No_G(T)$.
$\square$\smallskip

Thus we have proved that any nontoral grading has a coarsening
isomorphic to a grading induced by a subquasitorus of $A(j,t)$ for
some $j\in\{3,15,105,106,405\}$. Furthermore, by
Proposition~\ref{retra1}, the element $t$ can be taken to be the
identity. We can even remove two more possibilities for $j$, as
the following corollary shows.

\begin{co}\label{oroc1}
Each of the subgroups $A(3,\id)$ and $A(106,\id)$ of $\F4 $ is
conjugated to a subgroup of $A(105,\id)$.
\end{co}
Proof. We know that
$A(3,\id)=\span{\{t'_{-1,-1,1,1},t'_{1,-1,-1,i},\widetilde{
\si_3}\}}$ so that making $f_1=t '_{1,-1,-1,i}$,
$f_2=\widetilde{\si_3}$ and $f_3=t'_{-1,-1,1,1}$ we can apply the
previous proposition to $F=\{f_1,f_2,f_3\}$. Of course
$\{f_1,f_2\}$ is toral by  Lemma~\ref{Kasper}, while $A(3,\id)$ is
nontoral as proved in Proposition~\ref{nonto}. Thus $A(3,\id)$ is
conjugated to a  group of the form $\span{t_1,t_2,\si}$ with
$t_1,t_2\in\To$  and $\si\in\No$. Moreover $\si$ has order $2$
since it is conjugated to $t'_{-1,-1,1,1}$. We also know that
$\si=\widetilde{\si_i}t$ for some $t\in\To$ and
$i\in\{1,\ldots,1152\}$. Since $\si$ has order two, the same can
be said about $\si_i$.  That is, $A(3,\id)$ is conjugated to some
subgroup of $A(i,t)$ with $\si_i$ of order two. Furthermore
$A(i,t)$ is nontoral so that, applying Proposition~\ref{ehlc}, the
quasitorus $A(i,\id)$ is nontoral. Then  Proposition~\ref{nonto}
implies that, up to conjugacy, $i=3,15,105,106$ or $405$. We get,
by Proposition~\ref{retra1}, that the quasitorus $A(3,\id)$ is
conjugated to some subgroup of $A(i,\id)$. Since $\si_i^2=1$ the
only possible values of $i$ are $105$ and $405$. If we had a copy
of $A(3,\id)$ within $A(405,\id)$ then   this group (isomorphic to
$\Z_2^5$) should contain an element  of order $8$.

Consider now
$A(106,\id)=\span{\{t'_{-1,1,-1,1},t'_{1,-1,-1,-1},\widetilde{
\si_{106}}\}}$ and apply the previous proposition to
$F=\{f_1,f_2,f_3\}$ with $f_1=\widetilde{\si_{106}}$,
$f_2=t'_{-1,1,-1,1}$ and $f_3=t'_{1,-1,-1,-1}$. As before
$A(106,\id)$ is conjugated to a subgroup of a nontoral $A(i,\id)$
with $\si_i$ of order two. Again, up to conjugacy, $i=105$ or
$405$ but   $A(405,\id)$ contains no order six element.
$\square$\smallskip

Our objective in next subsection is to show that in fact any
nontoral grading is isomorphic to one produced by a quasitorus
contained in $A(j,\id)$ for $j=15,105,405$.

\subsection{Fine gradings}

Next we study the maximality of some of the previous  quasitori.
As a first step, their maximality in $\No$, the normalizer of the
maximal torus $\To$, is clear:

\begin{pr}\label{ladelc}
Let $A=A(j,\id)$ for $j\in\{15,105,405\}$. Then $A$ is its own centralizer in
$\No$, that is $\Co_{\Nos}(A)=A$.
\end{pr}
Proof.  To prove $\Co_{\Nos}(A)\subset A$ take $f\in
\Co_{\Nos}(A)$. Since $f\in \No$ there is some $i$ and some
$t\in\To$ such that $f=\widetilde{\si_i}t$.  Consider first the
possibility $j=405$. Of course $f$ commutes with each element in
$\tor{405}\subset A$ implying that $\widetilde{\si_i}$ does the
same. Consequently $\tor{405}\subset\tor{i}$.  According to table
(\ref{tablon}), this is possible only for $i=405$ or $i=748$ (take
into account that any $\tor{k}$ is isomorphic to some in the
table, that the unique groups in the table which may contain a
copy of $\Z_2^4$ are $\tor{405}$ and $\tor{748}$, and that the
orbits of $\si_{405}$ and $\si_{748}$ have cardinal one). So
 $\widetilde{\si_i}$ (equal to $\widetilde{\si_{405}}$ or $\widetilde{\si_{748}}$)
 commutes with  $\widetilde{\si_{405}}$, and since $f=\widetilde{\si_i}t$ also does, this forces
the commutativity of $t$ and $\widetilde{\si_{405}}$ so that
$t\in\tor{405}$ and $f=\widetilde{\si_i}t\in A(405,\id)$. Let us
consider now the case $j=15$. From $f=\widetilde{\si_i}t$ and the
fact that $f$ commutes with $\tor{15}$ we get that
$\widetilde{\si_i}$ commutes with $\tor{15}$.  Therefore
$\tor{15}\subset\tor{i}$ and as
$\tor{15}=\span{t'_{\om,1,\om^2,\om^2},t'_{1,\om,\om,1}}$, we
conclude $\si_i\cdot t=t$ for both $t=t'_{\om,1,\om^2,\om^2}$ and
$t=t'_{1,\om,\om,1}$. But using $(\ref{ecdlc1})$ we see that the
unique values of $i$ satisfying the above relations are $ 15,748$
and $1075$. Hence either $\si_i=\si_{15}$ or
$\si_i=\si_{1075}=\si_{15}^2$ or $\si_i=\si_{748}=\si_{15}^3$. In
any case $\widetilde{\si_i}$ commutes with $\widetilde{\si_{15}}$
(which is not automatic from the commutativity of $\si_i$ and
$\si_{15}$, but can be proved directly for those indices or
checking the equalities
$\widetilde{\si_{15}^n}=\widetilde{\si_{15}}^n$ for $n=2,3$).
Thus, recalling that $f$ commutes with $\widetilde{\si_{15}}$ we
get that $t$ commutes with $\widetilde{\si_{15}}$, that is,
$t\in\tor{15}\subset A(15,\id)$. Therefore
$f=\widetilde{\si_i}t\in A(15,id)$. Finally, we must investigate
the case $j=105$. Since $f=\widetilde{\si_i}t$ commutes with
$A(105,\id)$ and $\tor{105}$ is generated by $t'_{-1,1,-1,1}$,
$t'_{1,-1,-1,1}$ and $t'_{1,1,1,u}$ with $u\in F^\times$, we must
have that $\widetilde{\si_i}$ commutes with any of the mentioned
generators. So solving in $i$ the system $\si_i\cdot t=t$ for
$t\in\{t_{-1,1,-1,1},t_{1,-1,-1,1},t_{1,1,1,u} \mid u\in F^\times
\}$, we get that $i\in\{105,748\}$. But $\widetilde{\si_{105}}$
and $\widetilde{\si_{748}}$ commute and consequently
$\widetilde{\si_i}$ commutes with $\widetilde{\si_{105}}$. As
before we get that $t\in\tor{105}$, implying $f\in A(105,\id)$.
$\square$\smallskip

We must emphasize that the fact that the centralizer (in $\No$) of
a subgroup agrees with the subgroup itself does not imply that the
centralizer of the  group in $\F4$ agrees with the group. However
this will be the case for some special relevant subgroups in our
work. For the next result we must recall from 2.1 that if we have
a maximal abelian subgroup of semisimple elements in $\aut(L)$ for
a Lie algebra $L$, then the induced grading on $L$ is fine.
\begin{pr}\label{rpfi}
The gradings induced by $A(j,\id)$ for $j=15,105,405$ are fine.
\end{pr}

Proof. Let $A=A(j,\id)$ and suppose that the induced grading is
not fine. Then  there is some semisimple
$f\in\Co_{\Fs4}(A)\setminus A$ such
 the
grading induced by $A\cup\{f\}$ is a proper refinement of the
original grading (proper in the sense that it is different from
the given grading). Define $Z=\Co_{\Fs4}(\tor{j})$. Then the group
$\span{A\cup\{f\}}$ is an abelian subgroup of $Z$ and also its
closure $\overline{\span{A\cup\{f\}}}$ in the Zarisky topology.
But this is again a quasitorus, whence it is contained in the
normalizer of some maximal torus $T$ of $Z$. In particular
$\span{A\cup\{f\}}\subset\No_Z(T)$ and by construction also
$\To\subset Z$ so that there is some $p\in Z$ such that
$pTp^{-1}=\To$. Consequently $p\No_{\Fs4}(T)p^{-1}=
\No_{\Fs4}(\To)=\No$. Thus $p\span{A\cup\{f\}}p^{-1}\subset \No$
and
$$
pfp^{-1}, p\widetilde{\si_{j}}p^{-1}\in
\No\cap\Co_{\Fs4}(\tor{j})
$$
with $ptp^{-1}=t$ for any $t\in\tor{j}$.

For $j=105$ we have
$\No\cap\Co_{\Fs4}(\tor{105})=\To\cup\widetilde{\si_{105}}\To$,
taking into account the previous Proposition. Now we must analyze
different possibilities:
\begin{itemize}
\item  If $p\widetilde{\si_{105}}p^{-1}\in\To$, then $pAp^{-1}\subset\To$
and the grading induced by $A$ would be toral, which is a
contradiction.
\item If $pfp^{-1}\in\To$ and $p\widetilde{\si_{105}}p^{-1}=
\widetilde{\si_{105}}t$ for some $t\in\To$, it is clear that
$pfp^{-1}\in\tor{105}$ (it commutes with $\widetilde{\si_{105}}t$)
and $pAp^{-1}\subset A(105,t)$. Hence $\span{A\cup\{f\}}\subset
p^{-1}A(105,t)p$ and it cannot   produce a proper refinement. This
is another contradiction.
\item Thus we have $pfp^{-1}=\widetilde{\si_{105}}t_1$ and
$p\widetilde{\si_{105}}p^{-1}= \widetilde{\si_{105}}t_2$ for some
$t_1,t_2\in\To$. Since both elements commute we get easily that
$t_1t_2^{-1}\in\tor{105}$. Therefore $pfp^{-1}\in
p\widetilde{\si_{105}}p^{-1}\tor{105}=p\widetilde{\si_{105}}\tor{105}p^{-1}
\subset pAp^{-1}$ implying $f\in A$, a contradiction again.
\end{itemize}
The case $j=405$ is proved analogously since
$\No\cap\Co_{\Fs4}(\tor{405})= \To\cup\widetilde{\si_{405}}\To$.
However for $j=15$ we need some slight modifications since
$\No\cap\Co_{\Fs4}(\tor{15})= \To\cup\widetilde{\si_{15}}\To\cup(
\widetilde{\si_{15}})^2\To$. Let use the notation
$\si:=\widetilde{\si_{15}}$ for simplicity. Since $pfp^{-1},p\si
p^{-1}\in\No \cap\Co_{\Fs4}(\tor{15})$, we have several
possibilities.
\begin{itemize}
\item If $pfp^{-1}$ or $p\si p^{-1}$ is of the form $\si t$ for some $t\in\To$,
applying Corollary~\ref{coento}, there is some $s\in\To$ such that
$s\si t s^{-1}=\si$. This implies that $\tor{15}\cup\{\si\}\subset
sp\span{A\cup\{f\}}(sp)^{-1} \subset \No$ and by
Proposition~\ref{ladelc}
 we have
$\Co_{\No}(A(15,\id))=A(15,\id)$ implying
$sp\span{A\cup\{f\}}(sp)^{-1}= A(15,\id)=A$. Thus the grading
induced by $sp\span{A\cup\{f\}}(sp)^{-1}$, which is equivalent to
a proper refinement of the induced by $A$, is also equivalent to
$A$, a contradiction.
\item If $pfp^{-1}$ or $p\si p^{-1}$ is of the form $\si^2 t$ for some $t\in\To$,
applying again Corollary~\ref{coento}, there is some $s\in\To$
such that $s\si^2 t s^{-1}=\si^2=\si_{1075}$. This implies that
$\tor{15}\cup\{\si^2\}\subset sp\span{A\cup\{f\}}(sp)^{-1} \subset
\No$ and since $\tor{15}\cup\{\si^2\}$ generates
$A(15,\id)=A(1075,\id)$, we conclude as before.
\item And if $pfp^{-1}$ and $p\sigma p^{-1}$ are both in $\To$,
the grading is obviously toral.
 $\square$
\end{itemize}
We finish this section with one of our main results. Before
proceeding with its precise statement and proof we must realize
that any quasitorus $A$ in an algebraic group $G$ has a (not
necessarily unique) {\em maximal toral part}, a toral subgroup $B$
which is not contained in another toral subgroup of $A$. This is
trivial if the quasitorus $A$ has finite cardinal. Otherwise
consider the family $\mathcal F$ of all the toral subgroups of
$A$. The elements in $\mathcal F$ are toral in the sense that each
of them is contained in some torus of $G$.   To see that $\mathcal
F$ has maximal elements we consider a maximal subtorus ${T}$ of
$A$. For any $B\in\mathcal F$ containing ${T}$ (for instance,
${T}$ verifies this)  we have ${T}\subset B_0$ (the unit component
of $B$) and $B_0\subset A_0={T}$. Thus $B_0={T}$ and the quotient
group $B/{T}$ is finite and its cardinal is bounded by that of
$A/{T}$. Hence any $B\in\mathcal F$ with ${T}\subset B$ and such
that $B/{T}$ has maximal cardinal can be proved to be a maximal
element in $\mathcal F$. Such subgroups are what we shall
understand as {\em maximal toral} subgroups of the quasitorus $A$.
\begin{te}\label{main1}
Let $A\subset\F4 $ be a quasitorus, then up to conjugacy, $A$
falls in one of the following cases:
\begin{itemize}
\item $A\subset\To$ (maximal torus).
\item  $A\subset A(105,\id)$.
\item $A\subset A(405,\id)$.
\item $A=A(15,\id)$.
\end{itemize}
\end{te}
Proof. The precise meaning of this theorem is that there is some
$\phi\in\F4$ such that $A':=\phi A\phi^{-1}$ is in some of the
cases above. So we can replace at any moment $A$ with some of its
conjugated subgroups in the group $\F4$. From the beginning we
suppose that $A$ is nontoral. Consider a maximal toral subgroup
$B$ of $A$ which may be supposed to be a subgroup of the maximal
torus $\To$ defined in Subsection \ref{ewg}. Define now the group
$Z:=\Co_{\Fs4}(B)$ which contains to $A$ and to $\To$. Since $A$
is a quasitorus of $Z$ there is a maximal torus $T\subset Z$ such
that $A\subset\No_{Z}(T)$ (see the first paragraph in Section
\ref{seqt}). But $\To,T\subset Z$ are maximal tori in $Z$ so that
there exists $p\in Z$ such that $\To=pTp^{-1}$. Thus
$pAp^{-1}\subset\No_{Z}(\To)\subset\No_{\Fs4}(\To)=\No$ and since
$p$ centralizes $B$, the subgroup $B=pBp^{-1}$ is still a toral
maximal subgroup of $pAp^{-1}$. In this way, replacing $A$ by
$pAp^{-1}$ we can suppose that $A\subset\No$ with $B\subset\To$,
$B$ a maximal toral subgroup of $A$. Furthermore
$A=\span{B\cup\{f_1,\ldots f_k\}}$ with $f_i\in\No$ so that each
$f_i$ is of the form $\widetilde{\si_j}t$ with
$j\in\{1,\ldots,1152\}$ and $t\in\To$. Moreover, for any
$\widetilde{\si_j}t\in A$ we have $B\subset\tor{j}$ since any
element in $B$ centralizes $A$ and $\To$ and so it commutes with
$\widetilde{\si_j}$. On the other hand,
$\span{B\cup\{\widetilde{\si_j}t\}}$ being nontoral and contained
in $A(j,t)$, implies that $\si_j$ is in the orbit (under
conjugation) of $\si_3$, $\si_{15}$, $\si_{105}$, $\si_{106}$ or
$\si_{405}$. Next we analyze the different possibilities arising.
\begin{itemize}
\item
  If some of the $\si_j$'s is in the orbit of
$\si_{15}$ we can suppose that $f_1=\widetilde{\si_j}t$.  Then
$f_1$ is conjugated in $\No$ to $\widetilde{\si_{15}}t'$ for some
$t'\in\To$ and without loss in generality we can take $f_1=
\widetilde{\si_{15}}$, by  Corollary~\ref{coento} (the element
$s\in\To$ does not change neither $\To$ nor $B$). Thus
$B\subset\tor{15}\cong\Z_3^2$. Since $\span{B\cup\{f_1\}}$ is
nontoral, then it has at least three generators (see
Lemma~\ref{Kasper}). Consequently $B$ has at least two generators,
and this implies $B=\tor{15}$. So $A(15,\id)\subset A\subset\No$
and applying Proposition~\ref{ladelc}, we have
$\Co_{\Nos}(A(15,\id))=A(15,\id)$, hence $A=A(15,\id)$.
\item
Some of the $\si_j$'s is conjugated to $\si_{106}$. As before
$f_1$ can be taken to be $f_1=\widetilde{\si_{106}}$ and
$B\subset\tor{106}\cong \Z_2^2$. Since $B$ must have at least two
generators we have $B=\tor{106}$. Now, for any other
$f_i=\widetilde{\si_k}t$ we must have $\tor{106}\subset \tor{k}$.
The commutativity of $f_1$ and $f_i$ implies that of $\si_{106}$
and $\si_{k}$. But a computation of the number of $k$'s satisfying
both conditions
$$\begin{cases}
\tor{106}\subset
\tor{k},\cr
\si_{106}\si_k=\si_k\si_{106},
\end{cases}$$
reveals that there are only six possible values of $k$. Since
$\si_{106}$ has order $6$, obviously the six powers $\si_{106}^n$
with $n\in\{0,\ldots,5\}$ satisfy the conditions so that there is
$n\in\mathbb N$ and $t_1\in\To$ with $f_i=\widetilde{\si_{106}}^n
t_1$. But $t_1\in\tor{106} $ again by the commutativity of $f_i$
with $f_1$, so $f_i\in A(106,\id)$. From here $A=A(106,\id)$ which
up to conjugacy is a subgroup of $A(105,\id)$ by
Corollary~\ref{oroc1}.
\item
 Suppose now that some of the $\si_j$'s is in the
orbit of $\si_3$. So we can suppose $f_1=\si_3$ without loss of
generality, and $B\subset
\tor{3}=\span{\{t'_{-1,-1,1,1},t'_{-1,1,-1,i}\}}$, which is
isomorphic to $\Z_2\times\Z_4$. But again by Lemma~\ref{Kasper},
the subgroup $B$ must have at least two generators since
$\span{B\cup\{f_1\}}$ is nontoral. Whence either $B=\tor{3}$ or
$B=\span{\{t'_{-1,-1,1,1},t'_{1,1,1,-1}\}}$. But this last
possibility can not occur because then $\span{B\cup\{f_1\}}=
\span{\{t'_{-1,-1,1,1},t'_{1,1,1,-1},\widetilde{\si_3}\}}$ would
be toral because $\widetilde{\si_3}^4=t_{1,1,1,-1}$, hence
$\span{B\cup\{f_1\}}$ would have two generators.
 So necessarily
$B=\tor{3}$ and for any other $f_i=\widetilde{\si_k}t$ we must
have $\tor{3}\subset\tor{k}$ and $\si_3 \si_k=\si_k\si_3$. We
obtain only four possible $k$'s satisfying the previous
conditions. Hence $\si_k\in\{\si_3^n\mid n=0,1,2,3\}$ and
$f_i=\widetilde{\si_3^n}t_1$ for some $t_1\in\To$ (in fact $t_1=t$
because $\widetilde{\si_3^n}=\widetilde{\si_3}^n$). Then
$t_1\in\tor{3}$ and as in previous cases $A=A(3,\id)$. But this is
conjugated to a subgroup of $A(105,\id)$ by Corollary~\ref{oroc1}.
\item We may suppose now that each $\si_j$ is in the orbit of $\si_{105}$
or of $\si_{405}$. We prove next that either all the $\si_j$ are
in the orbit of $\si_{105}$, or all of them are in the orbit of
$\si_{405}$. Otherwise, after re-ordering and applying
Corollary~\ref{coento} we can take $f_1=\widetilde{\si_{405}}$ and
$f_2=\widetilde{\si_{105}}t$ for some $t\in\To$. Thus
$B\subset\tor{105}\cap \tor{405}$. Moreover
$B\subsetneq\widetilde{B}:=\span{B\cup\{f_1f_2\}}\subset A$ and
$\widetilde{B}\subset\span{(\tor{105}\cap
\tor{405})\cup\{\widetilde{\si_{405}}\widetilde{\si_{105}}t\}}\subset
A(1048,t_1)$ for some $t_1\in\To$, taking into account that
$\si_{405}\si_{105}=\si_{1048}$. But $A(1048,t_1)$ is toral since
$\si_{1048}$ is in the orbit of $\si_{142}$ (see Propositions
\ref{ehlc} and \ref{nonto}). This contradicts the maximal toral
nature of $B$ within $A$.
\item
If all the $\si_j$'s are in the orbit of $\si_{405}$, since this
orbit has cardinal one, then we can take
$f_1=\widetilde{\si_{405}}$ and any other $f_i$ of the form
$\widetilde{\si_{405}}t$ with $t\in\tor{405}$. So $A\subset
A(405,\id)$.
\item
If all the $\si_j$'s are in the orbit of $\si_{105}$. The unique
elements in this orbit which commute with $\si_{105}$ are
$\si_{105}$, $\si_{403}$, $\si_{429}$ and $\si_{1011}$. So we can
take $f_1=\widetilde{\si_{105}}t_1$ and any other $f_i$   of the
form $\widetilde{\si_{k}}t_2$ for some $t_1,t_2\in\To$ and
$k\in\{105,403,429, 1011\}$. Then
$B\subset\widetilde{B}:=\span{B\cup\{f_1f_i\}}\subset A$ and
depending on the values of $k$ we have:
\begin{itemize}
\item For $k=403$, since $\si_{105}\si_{403}=\si_{1050}$, we have
$f_1f_i=\widetilde{\si_{105}}t_1\widetilde{\si_{403}}t_2=\widetilde{
\si_{1050}}t$ for some $t\in\To$. So $\widetilde{B}\subset
A(1050,t)$ which is toral because $\si_{1050}$ is in the orbit of
$\si_{103}$ (see again Propositions \ref{ehlc} and \ref{nonto}).
This contradicts the maximality of $B$ among the toral subgroups
of $A$.
\item For $k=429$ we have $\si_{105}\si_{429}=\si_{1009}$ so that
$B\subsetneq\widetilde{B}\subset A(1009,t)$ for some $t\in\To$.
This is toral since $\si_{1009}$ is conjugated to $\si_{103}$ in
$\W$. Again a contradiction.
\item For $k=1011$ we have $\si_{105}\si_{1011}=\si_{427}$ and
$B\subsetneq\widetilde{B}\subset A(427,t)$, for $t\in\To$. This is again toral
because $\si_{427}$ is in the orbit of $\si_{103}$.
\end{itemize}
These contradictions lead us to the conclusion that $k=105$
necessarily. So $f_i$ is of the form $\widetilde{\si_{105}}t_2$
and the commutativity of $f_1$ and $f_i$ implies that
$t_1t_2^{-1}\in\tor{105}$. Hence $f_i\in f_1\tor{105}$ (for any
other $i$) implying $A\subset A(105,t_1)$. Thus, after conjugation
$A$ is a subgroup of $A(105,\id)$ by Proposition~\ref{retra1}.
$\square$
\end{itemize}

 \begin{co}\label{fiinas}
 Up to equivalence the unique fine gradings on $\f4$ are: (1) the
Cartan grading, (2) the induced by $A(105,\id)$, (3) the induced by
$A(405,\id)$, and (4) the induced by $A(15,\id)$.
 \end{co}

We summarize the information about these gradings in the following table
\begin{center}
\begin{tabular}{|c|c|c|}
\hline
Quasitorus & Universal & Type\cr
           & Grading group & \cr
\hline\hline $\To$ & $\Z^4$ & $( 48,0,0,1)$\cr $A(15,\id)$ &
$\Z_3^3$  & $( 0,26)$\cr $A(105,\id)$ & $\Z_2^3\times\Z$ & $(
31,0,7)$\cr $A(405,\id)$ & $\Z_2^5$   & $( 24,0,0,7)$\cr \hline
\end{tabular}
\centerline{\small Fine gradings on $\f4$}
\end{center}

Let us observe that the groups given in the table above are the
universal grading groups.
 Because for any of the quasitorus  $Q$ in the table, we have proved in
Proposition~\ref{rpfi} that $Q$ is its own centralizer in $\F4$,
that is, $Q$ is a MAD, and according to
    Subsection 2.1, $\X(Q)$ is the universal   group of the
    related grading.\smallskip

    \noindent \textbf{Remark 3.} The four fine gradings on $\f4$
    are also fine when are considered as Lie gradings (as in
    \cite{Zass}) instead of being considered as group gradings.
    This observation is   pertinent because the question about
    the existence of a grading group is still on the
    table.\smallskip

On the other hand, we observe that every homogeneous element in
the Cartan grading is either semisimple or nilpotent (according to
its membership to the zero component) and that all the homogeneous
elements in the three remaining fine gradings are semisimple, as a
consequence of the nullity of their zero components
 (\cite[Corollary after Th.~3.4]{enci}). This happens not only in $\f4$.

\begin{pr}
In a fine grading on a simple Lie algebra, every homogeneous
element is either semisimple or nilpotent.
\end{pr}
Proof.  If $L=\oplus_{g\in G}L_g$ is a fine grading, the group $Q$
of the automorphisms of $L$ for which  all the homogeneous
components are  eigensubspaces   is just the MAD which produces
the grading.

First of all, if $x\in L_g$  then $x$ is uniquely represented as
$x=x_s+x_n$, where $x_s$ and $x_n$ are respectively semisimple and
nilpotent homogeneous elements in $L_g$ verifying $[x_s,x_n]=0$
\cite[Th.\,3.3]{enci}. Thus, either every element in $L_g$ is
semisimple or there is a nilpotent element $x\in L_g$. Let us see
that in this case every element in $L_g$ is nilpotent. According
to \cite[Th.\,3.4]{enci}, there is a semisimple element $h\in L_e$
and a nilpotent element $y\in L_{-g}$ such that $[h,x]=2x$,
$[h,y]=-2y$ and $[x,y]=h$. Notice then that ${\rm exp}{\rm
ad}(h)\in\aut(L)$ is an automorphism which leaves invariant all
the homogeneous components. Moreover, ${\rm exp}\,{\rm ad}(h) $
commutes with every $f\in Q$, because if $v$ is a homogeneous
element with $f(v)=\alpha v$, then $({\rm exp}\,{\rm
ad}(h))f(v)=\alpha {\rm exp}\,{\rm ad}(h)(v)=f({\rm exp}\,{\rm
ad}(h)(v))$ since ${\rm exp}\,{\rm ad}(h)(v)$ belongs to the same
homogeneous component. By maximality, ${\rm exp}\,{\rm ad}(h)\in
Q$. Since ${\rm exp}\,{\rm ad}(h)(x)=e^2x$, there are some
automorphisms $f_1,\dots, f_s\in Q$ and scalars
$\alpha_1,\dots,\alpha_s\in F^\times$ such that $L_g=\{v\in L\mid
f_1(v)=\alpha_1v,\dots,f_s(v)=\alpha_sv,{\rm exp}\,{\rm ad}
(h)(v)=e^2v\}$. Therefore $g\in G$ is an element of infinite
order, and, by
 \cite[Prop.\,3.5]{enci}, all the elements in $L_g$ are nilpotent.
$\square$\smallskip

In particular, the zero homogeneus component of a fine grading is
an abelian (toral) subalgebra, with dimension equal to the
dimension of the quasitorus producing the grading. More
information on homogeneous semisimple and nilpotent elements of
gradings on semisimple Lie algebras can be found in
\cite{Vinberg}.

\subsection{Nontoral gradings}

Up to the moment we have described the fine gradings on $\f4$. If
we want to find all the nontoral ones, it suffices to describe all
the nontoral coarsenings of the fine ones.  For this, we must find
all the (nontoral) subquasitori $Q$ of the tori producing fine
gradings up to conjugacy.

Consider first the fine grading provided by
$A(105,\id)=\span{\tor{105}\cup\{\widetilde{\si_{105}}\}}$.  Let
us use the following notations $g_1:=t'_{-1,1,-1,1}$,
$g_2:=t'_{1,-1,-1,1}$, $g_3:=\widetilde{\si_{105}}$ and
$g_4:=t'_{1,1,1,-1}$.  Since the number of generators of the group
$H=\X(Q)$ must be at least three, we know that $H$ has the
following possibilities: either $H\cong\Z_2^3\times\Z_m$,     or
$H\cong\Z_2^2\times\Z_{2m}$,  or $H\cong\Z_2^2\times\Z$. Let us
analyze first the case $H\cong\Z_2^3$.  In this case the
quasitorus $Q$ providing the coarsening has three order-two
generators $\varphi_i$ with $i=1,2,3$. Thus $\varphi_i=
g_1^{n_i}g_2^{m_i}g_3^{l_i}g_4^{s_i}$ and
$(\varphi_1,\varphi_2,\varphi_3)=(g_1,g_2,g_3,g_4)\cdot M$ where
$M$ is the $3\times 3$ integer matrix
$$\small\begin{pmatrix}n_1 & n_2 & n_3\cr
m_1
& m_2 & m_3\cr
l_1 & l_2 & l_3\cr
s_1 & s_2 & s_3\cr
\end{pmatrix}.
$$
We are using here the action of $ M_{k\times l}(\Z)$ on any power
group $G^k$ given by $G^k\times M_{k\times l} (\Z)\to G^l$ such
that $(g_1,\ldots,g_k)\cdot(n_{ij}):= (g'_1,\ldots,g'_l)$ where
$g'_i=\prod_j g_j^{n_{ji}}$. Now, there are two actions that we
can use to simplify the matrix $M$ without changing the quasitorus
$Q$ (in the worst case $Q$ changes to some of its conjugated
ones).  First we can act on the columns of the matrix by
elementary operations (exchanging columns or adding one column to
another). All this can be made module two. This reduces the
possible matrices $M$ to a few, but there is a second action on
the matrix coming from conjugacy of elements in $\F4$. Consider
the subgroup $G$ of $\F4$ fixing $A(105,\id)$ by conjugation, that
is, $f\in G$ if and only if $f A(105,\id)f^{-1}=A(105,\id)$. The
appendix contains relevant information on certain elements in $G$.
 By passing from $Q$ to $fQf^{-1}\subset A(105,\id)$ the matrix $M$ transforms
into a new matrix $M'$. We are proving that the joint action by
elementary operations on the columns together with the action of
$G$ by conjugation, either moves $M$ to the matrix
$$\small\begin{pmatrix}
1 &  0 &   0\cr
  0
&   1 &   0\cr
  0 &   0 &   1\cr
  0 &   0 &   0\cr
\end{pmatrix},$$
or the induced grading is toral. First we observe that if this
grading   is nontoral there must be some nonzero entry in the
first two rows of $M$. Second, $\widetilde{\si_{468}}\in G$ acts
on $M$
 permuting  its first two rows (it permutes $g_1$ and
$g_2$ fixing $g_i$ for $i=3,4$). So we can suppose that the first
row in $M$ is $(1\ 0\  0)$ (after operations on columns). Hence
the first two rows are $\tiny\begin{pmatrix}1 & 0 & 0\cr 0 & 1 &
0\end{pmatrix}$, $\tiny\begin{pmatrix}1 & 0 & 0\cr 0 & 0 &
0\end{pmatrix}$ or $\tiny\begin{pmatrix}1 & 0 & 0\cr 1 & 0 &
0\end{pmatrix}$ also after operations on columns. However
$\widetilde{\si_{94}}\in G$ acting on the second matrix produces
the third one. The grading induced by the third matrix is
$\span{g_1g_2g_3^{n_1}g_4^{n_2},g_3^{m_1}g_4^{m_2},g_3^{l_1}g_4^{l_2}}$
where the couples $(m_1,m_2)$ and $(l_1,l_2)$ are linearly
independent in the vector space $\Z_2^2$ (otherwise the grading
has two generators and so is toral). Thus $(n_1,n_2)$ is a linear
combination of them and the grading can be written as
$\span{g_1g_2f^kg^h,f,g}$ where $f,g\in\span{g_3,g_4}$. But
$\span{g_1g_2f^kg^h,f,g}=\span{g_1g_2,f,g}$ hence we have the
(agreeing) possibilities
$\span{g_1g_2,g_3,g_4}=\span{g_1g_2,g_3g_4,g_3}=\span{g_1g_2,g_3g_4,g_4}$.
So in this case we must only worry about the grading
$\span{g_1g_2,g_3,g_4}$, which can be proved to be toral (seen as
a grading on the Albert algebra, this comes from a necessarily
toral $\Z_2\times\Z_2$-grading on the octonions refined with the
$\Z_2$-grading on $H_3(F)$, see Subsection \ref{gfo}, but
alternatively  we next prove that is conjugated to
$\span{g_1,g_2,g_4}$). The first matrix above produces the grading
$\span{g_1g_3^{n_1}g_4^{n_2},g_2g_3^{m_1}g_4^{m_2},g_3^{l_1}g_4^{l_2}}$.
But $\widetilde{\si_{485}}$ permutes $g_1$ and $g_2$ while
$\widetilde{\si_{485}}g_3\widetilde{\si_{485}}^{-1}=g_3g_4$ and
$\widetilde{\si_{485}}g_4\widetilde{\si_{485}}^{-1}=g_4$, so we
can suppose $l_1=1, l_2=0$ or $l_1=0, l_2=1$. This reduces the
possibilities to $\span{g_1g_4^{n_2},g_2g_4^{m_2},g_3}$ or
$\span{g_1g_3^{n_1},g_2g_3^{m_1},g_4}$. The four different
possibilities for the first case are $\span{g_1,g_2,g_3}$,
$\span{g_1g_4,g_2,g_3}$, $\span{g_1,g_2g_4,g_3}$,
$\span{g_1g_4,g_2g_4,g_3}$, but conjugating by
$\widetilde{\si_{103}}\in G$, the first and second gradings turn
out to be isomorphic, so as the third and the forth ones; while
conjugation by $\widetilde{\si_{468}}\in G$ proves the isomorphism
between the second and the third gradings. On the other hand the
four possibilities for $\span{g_1g_3^{n_1},g_2g_3^{m_1},g_4}$ are
$\span{g_1,g_2,g_4}$ (which is obviously toral),
$\span{g_1g_3,g_2,g_4}$, $\span{g_1,g_2g_3,g_4}$ and
$\span{g_1g_3,g_2g_3,g_4}$. But conjugation by
$\widetilde{\si_{468}}$ makes evident the isomorphism between the
second and the third gradings while $\span{g_1g_3,g_2g_3,g_4}=
\span{g_1g_2,g_2g_3,g_4}=\span{g_1g_2,g_1g_3,g_4}$ and conjugation
by $\widetilde{\si_{ 491}}\in G$ proves that this is isomorphic to
$\span{g_1g_3,g_2,g_4}$. Now it is possible to show that there
exists an element $\psi\in G$ such that $\psi g_i\psi^{-1}=g_i$
for $i=1,4$ while $\psi g_2\psi^{-1}=g_3$ and $\psi
g_3\psi^{-1}=g_2$ (the element $\psi$ can be taken of order $4$).
So $\psi\span{g_1g_3,g_2,g_4}\psi^{-1}=\span{g_1g_2,g_3,g_4}$
whose torality has been previously stated. Moreover,
$\span{g_1g_2,g_3,g_4}$ is conjugated to $\span{g_1,g_3,g_4}$ by
$\widetilde{\si_{94}}$ and this is conjugated to
$\span{g_1,g_2,g_4}$ by $\psi$. This last one is obviously toral.
Summarizing the results in this paragraph we have:

\begin{pr}\label{lprdelco}
If $Q\subset A(105,\id)$ is a quasitorus with $\X(Q)\cong\Z_2^3$
then it is conjugated to $\span{g_1,g_2,g_3}$ if $Q$ is nontoral,
and to $\span{g_1,g_2,g_4}$ if it is toral. Moreover the
conjugating element can be taken to fix the subgroup $A(105,\id)$.
The $\Z_2^3$-grading induced by $\span{g_1,g_2,g_3}$ is of type
$(0,0,1,0,0,0,7)$.
\end{pr}
Proof. We have proved in the previous paragraph that $Q$ is
conjugated to either $\span{g_1,g_2,g_3}$ or  $\span{g_1,g_2
,g_4}$. Obviously   $\span{g_1,g_2 ,g_4}$ is toral. On the other
hand, the grading induced by  $\span{g_1,g_2,g_3}$ is nontoral, of
type $( 0,0,7,0,0,1)$ in the Albert algebra and of type $(
0,0,1,0,0,0,7)$ in $\f4$. Its nontoral character is obvious
noticing that the dimension of its $0$-homogeneous component in
$\f4$ is $3$. $\square$

\begin{pr}\label{grandes}
If $Q\subset A(105,\id)$ is a nontoral quasitorus such that
$\X(Q)\cong\Z_2^3\times\Z_m$, with $m>1$, then
 up to conjugacy there is $v\in F^\times$ a primitive $m$-root of the unit such that
 $Q=\span{g_1,g_2,g_3,t'_{1,1,1,v} }$.
\end{pr}
Proof. The quasitorus $Q=\span{\phi_i}_{i=1}^4$ is generated by
$\phi_1,\phi_2,\phi_3$   order-two elements in $\F4$ and $\phi_4$
  an order-$m$ element in $\F4$. We can apply the previous
proposition to $Q':=\span{\phi_1, \phi_2,\phi_3}$. So we can
assume (by conjugation) that either (1) $Q'=\span{g_1,g_2,g_3}$,
or (2) $Q'=\span{g_1,g_2,g_4}$, with
$\phi_4=g_1^ng_2^kg_3^lt'_{1,1,1,u}$ for some $n,k,l\in\{0,1\}$
and some $u\in F^\times$. From $\phi_4^m=\id$, it follows that
$lm$ is even (either $l$ or $m$) and
$t'_{(-1)^n,(-1)^k,(-1)^{n+k},u}$ has order $m$, in particular
$u^m=1$.

In case (1), $Q=\span{g_1,g_2,g_3, g_1^ng_2^kg_3^lt'_{1,1,1,u} }
=\span{g_1,g_2,g_3,t'_{1,1,1,u}}$, but $t'_{1,1,1,u}\notin Q'$, so
if $u$ would have order $m'$ (divisor of $m$), $\X(Q)$ would be
isomorphic to $\Z_2^3\times \Z_{m'}$. Hence   $u$ is a primitive
$m$-root of the unit.

  In case (2), $Q=\span{g_1,g_2,g_4, g_1^ng_2^kg_3^lt'_{1,1,1,u} }
=\span{g_1,g_2,t'_{1,1,1,-1},g_3^lt'_{1,1,1,u}}$. Since $Q$ is
nontoral, $l=1$. Hence $m=2m'$ must be even, with
$u^{m'}\in\{\pm1\}$. If $u^{m'}=-1$, $(g_3t'_{1,1,1,u})^{m'}$
would be $t'_{1,1,1,-1}$ if $m'$ is even, and $g_3t'_{1,1,1,-1}$
if $m'$ is odd. In the first case $Q=\span{g_1,g_2,
g_3t'_{1,1,1,u}}$, a contradiction with the number of generators
of $\X(Q)$. In the other case, $g_3\in Q$, hence
$Q=\span{g_1,g_2,t'_{1,1,1,-1}, g_3,t'_{1,1,1,u} } =\span{g_1,g_2,
g_3,t'_{1,1,1,u} }$ and so $\X(Q)$ would be a subgroup of
$\Z_2^3\times \Z_{m'}$, a contradiction. Thus $u^{m'}=1$ and
$(g_3t'_{1,1,1,u})^{m'}$ would be $t'_{1,1,1,1}=\id$ if $m'$ is
even, and $g_3 $ if $m'$ is odd. In the first case again $\X(G)$
would be a subgroup of $\Z_2^3\times \Z_{m'}$. So that we have the
case $m'$ odd, in which $g_3\in Q$. In this way
$Q=\span{g_1,g_2,t'_{1,1,1,-1},g_3,t'_{1,1,1,u}} =\span{g_1,g_2,
g_3,t'_{1,1,1,-u} }$, since $(t'_{1,1,1,-u})^{m'}=t'_{1,1,1,-1}$.
But now $v=-u$ is the required primitive $m$-root of the unit.

Notice that the obtained quasitori $Q =\span{g_1,g_2,
g_3,t'_{1,1,1,v} }$ are obviously nontoral because they contain $
\span{g_1,g_2, g_3 }$.
  $\square$\smallskip

  The induced gradings by the previous quasitori depend on $m$.
  For $m=2$ we are talking
  about a $\Z_2^4$-grading of type $(1,8,0,0,7)$, for $m=3$ it is
  a $\Z_2^3\times\Z_3$-grading of type $(3,14,7)$, for $m=4$ it is
  a $\Z_2^3\times\Z_4$-grading of type $(17,7,7)$, but for $m\ge5$
  all the gradings are equivalent to the one produced by
  $A(105,\id)$, since they are all of type $(31,0,7)$. In general
     two gradings  having the same type are not necessarily
equivalent,  but of course they are equivalent if the quasitorus
producing one of the gradings is contained in the other one, since
the homogeneous components of the former are pieces of the
homogeneous components of the latter.\smallskip

  \noindent\textbf{Remark 4. }If $Q'\subset Q$ are quasitori whose induced gradings
  are of the same type, then these gradings are
  equivalent.\smallskip

To continue the study of nontoral gradings coming from
subquasitori of $A(105,\id)$ we must analyze those quasitori
$Q\subset A(105,\id)$ with $\X(Q)\cong\Z_2^2\times\Z_m$, where $m$
must be even (otherwise $Q$ would have two generators and the
grading would be toral). Moreover $m$ must be a multiple of 4,
because we have already studied the cases $m=2$ in
Proposition~\ref{lprdelco} and $m=2m'$ with $m'$ odd in
Proposition~\ref{grandes}.

\begin{pr}\label{pequenos}
If $Q\subset A(105,\id)$ is a nontoral quasitorus such that
$\X(Q)\cong\Z_2^2\times\Z_{4m}$ with $m\ge1$,   then
 up to conjugacy there is $v\in F^\times$ a primitive $4m$-root
 of the unit such that
 $Q=\span{g_1,g_2,g_3t'_{1,1,1,v} }$.
\end{pr}
Proof.
 Suppose that $\{\phi_1,\phi_2,\phi_3\}$ is a set of
generators of $Q$ with $\phi_1$ and $\phi_2$ of order two and
$\phi_3$ of   order $2k=4m$. Then $\phi_3^k$ is an order-two
element and we can suppose (by conjugation) that the subgroup
$\span{\phi_1,\phi_2,\phi_3^k}$ is either $\span{g_1,g_2,g_3}$ or
$\span{g_1,g_2,g_4}$, according to Proposition~\ref{lprdelco}.
Besides, $\phi_3=g_1^ng_2^mg_3^lt'_{1,1,1,u}$ for some
$n,m,l\in\{0,1\}$ and $u\in F^\times$. As $k$ is even,
$\phi_3^k=t'_{1,1,1,u^k}$. Since $\phi_3^k$ has exactly order two,
$u^k=-1$  and $\phi_3^k=g_4$.

But $\span{\phi_1,\phi_2,\phi_3^k}=\span{g_1,g_2,g_3}$ is not
possible, since $g_4\notin\span{g_1,g_2,g_3}$, so necessarily
$\span{\phi_1,\phi_2,\phi_3^k}=\span{g_1,g_2,g_4}$. Thus $Q=
\span{\phi_1,\phi_2,\phi_3^k,\phi_3}=\span{g_1,g_2,g_4,g_3^lt'_{1,1,1,u}}$
and $l=1$ since $Q$ is nontoral. But $(g_3t'_{1,1,1,u})^k=g_4$ so
that $Q=\span{g_1,g_2,g_3 t'_{1,1,1,u}}$, where $u$ has order
$2k=4m$. $\square$\smallskip

Notice that in this case all the $\Z_2^3$-coarsenings of the
gradings are toral. In spite of that, these gradings are nontoral,
because    the dimension of the zero homogeneous component is $3$
in case $m=1$ and $1$ in the remaining cases. If $m=1$, we obtain
a $\Z_2^2\times \Z_4$-grading of type $(0,8,2,0,6)$, if $m=2$ we
have a $\Z_2^2\times \Z_8$-grading of type $(19,6,7)$, but for
$m\ge3$ we obtain gradings of   type $(31,0,7)$,   so equivalent
to the fine $\Z_2^3\times \Z$-grading, again by the previous
remark.\smallskip

Finally we analyze the case $\X(Q)=\Z_2^2\times\Z$. The quasitorus
$Q$ must be the direct product  of   a one-dimensional subtorus
$P$ of $A(105,\id)$ times $\span{\phi_1,\phi_2}$, with  $\phi_i$
two elements  of order two. But necessarily $P=\{t'_{1,1,1,u}\mid
u\in F^\times\}$ (it is the unique nontrivial subtorus of
$A(105,\id)$). Changing the generators if necessary (to remove
$g_4$), we can write $(\phi_1,\phi_2)=(g_3,g_2,g_1)M$ where
  $M$ is a $3\times 2$ matrix with entries
in $\Z_2$. The first row is nonzero because otherwise the grading
would be toral, and making column operations we can suppose that
it is $(1\,0)$. Besides there must be   some $1$ in each column.
As $\widetilde{\si_{468}}$ interchanges the second and third row,
and by doing column operations, the first two rows in $M$ can be
taken to be   $\tiny\begin{pmatrix}1 & 0\cr 0 & 1\end{pmatrix}$.
Thus either  $Q=P\span{g_3g_1,g_2}$, or $Q=P\span{g_3,g_1g_2}$, or
$Q=P\span{g_3g_1,g_2g_1}$. By applying the element
$\widetilde{\si_{491}}$ to the third one, we obtain
$P\span{g_1g_2g_3g_4,g_2g_4}=P\span{g_1g_3,g_2g_4}=P\span{g_1g_3,g_2}$,
the first quasitorus. This one is conjugated to the second one by
means of $\psi$ (see Appendix). Besides we can replace $g_1g_2$ by
$g_1$ by using $\widetilde{\si_{94}}$. Summarizing these facts we
have:
\begin{pr}\label{ootra}
If $Q\subset A(105,\id)$ is nontoral and $\X(Q)\cong\Z_2^2\times\Z
$ then $Q$ is conjugated to $\span{g_1, g_3,t'_{111u}\mid u\in
F^\times }$.
\end{pr}

Notice that the induced  grading on $\f4$  is of type $( 31,0,7)$,
so it is again equivalent to the produced by the whole
$A(105,\id)$.
\smallskip

 The following step is to examine the proper subquasitori of
$A(15,\id)$, but all of them are toral since $A(15,\id)$ is
isomorphic to $\Z_3^3$ hence any proper subquasitorus has a system
of generators  of cardinal $\le 2$ (Lemma~\ref{Kasper}).

Thus, to finish our classification of subquasitori of the maximal
quasitori, we must analyze $A(405,\id)$.
 This group is isomorphic to $\Z_2^5$.  Indeed we
have $A(405,\id)=\span{\tor{405}\cup\{\widetilde{\si_{405}}\}}$
where $\tor{405}=\{t'_{xyzu}\mid x^2=y^2=z^2=u^2=1\}\cong\Z_2^4$
and $\sigma:=\widetilde{\si_{405}}$ is an order two element.
Consider a subquasitorus $Q\subset A(405,\id)$ with three
order-two generators. If $Q$ is nontoral, then $\si t\in Q$ for
some order-two $t\in \To$. But applying Corollary~\ref{coento} we
can conjugate $Q$ to some new quasitorus $Q'$ with $\si\in Q'$.
Thus we can suppose since the beginning that $\si\in Q$, and by
elementary operations we can take $Q=\span{\si,t_1,t_2}$ where
$t_i\in\To$ are order-two elements.

It is   known that there are two conjugacy classes of elements of
order 2 in $\f4$, related to the diagrams obtained by removing the
nodes marked with the number 2 in the affine Dynkin diagram
\cite[Ch.\,8]{Kac}. The first automorphism, related to

\vskip .3cm \hskip 3.5cm \vbox{ \xy <1cm, 0cm>:
\POS(-1,0)*+{\xycircle(.1,.1){}} \POS(-.9,0) \ar @{-} +(1,0)
\POS(0.2,0)*+{\xycircle(.1,.1){}} \POS(0.3,0) \ar @{-} +(1,0)
\POS(1.4,0) *+{\xycircle(.1,.1){}} \POS(1.5,0) \ar @{=>} +(.5,0)
\POS(1.9,0)\ar @{=} +(.5,0) \POS(2.5,0) *+{\xycircle(.1,.1){}}
\POS(2.6,0)\ar @{-}+(1,0) \POS(3.7,0)*+{\xycircle(.1,.1){}}
\POS(3.7,0)*+{\times}
\endxy}
\vskip .3cm

\noindent fixes a subalgebra  of type $\frak b_4$, of dimension
36, and the second one, related to

\vskip .3cm \hskip 3.5cm \vbox{ \xy <1cm, 0cm>:
\POS(-1,0)*+{\xycircle(.1,.1){}} \POS(-.9,0) \ar @{-} +(1,0)
\POS(0.2,0)*+{\xycircle(.1,.1){}} \POS(0.3,0) \ar @{-} +(1,0)
\POS(1.4,0) *+{\xycircle(.1,.1){}} \POS(1.5,0) \ar @{=>} +(.5,0)
\POS(1.9,0)\ar @{=} +(.5,0) \POS(2.5,0) *+{\xycircle(.1,.1){}}
\POS(2.6,0)\ar @{-}+(1,0) \POS(3.7,0)*+{\xycircle(.1,.1){}}
\POS(0.2,0)*+{\times}
\endxy}
\vskip .3cm

\noindent  fixes a subalgebra of type $\frak c_3\oplus \frak a_1$,
of dimension 24. But if two elements in a torus are conjugated,
they are conjugated inside the normalizer. This means that the
Weyl group acts on $\tor{405}$ producing two orbits, apart of the
trivial one, characterized by the fact that the dimensions of
their fixed parts (the number of $1$'s in the list (\ref{tpm}))
are 36 and 24, respectively. This could have been checked directly
making act $\W$ by   computer. In the first orbit there are three
elements, $t'_{1,1,1,-1}$, $t'_{-1,-1,-1,-1}$ and
$t'_{-1,-1,-1,1}$, and in the second orbit the remaining 12
elements. Thus we can move the element $t_1$ in $Q$, since not
only $\si_{405}=-\id$ clearly commutes with any element in $\W$,
but $\si\widetilde{\si_j}=\widetilde{\si_j}\si$ for any
$j\in\{1,\dots1152\}$. Therefore the possibilities for $Q$ are (1)
$Q=\span{\si,t'_{-1,1,1,1},-}$ and (2)
$Q=\span{\si,t'_{1,1,1,-1},-}$, as  $t'_{-1,1,1,1}$ and
$t'_{1,1,1,-1}$ are representatives of the two orbits. In the
first case, the third element can be computed by considering the
subgroup of $\W$ fixing $t_{-1,1,1,1}$ and the orbits it produces
on the set of order-two elements different from $t_{-1,1,1,1}$.
Indeed, that subgroup produces three orbits such that any
order-two element different from $t_{-1,1,1,1}$ is conjugated to
either $t'_{1,-1,1,1}$, or $t'_{1,1,1,-1}$ or  $t'_{1,-1,-1,-1}$
by an element in $\W$ fixing $t'_{-1,1,1,1}$. As a consequence the
possibilities for $Q$ are $\span{\si,t_{-1,1,1,1},t_{1,-1,1,1}}$,
 $\span{\si,t_{-1,1,1,1},t_{1,1,1,-1}}$ and
$\span{\si,t_{-1,1,1,1},t_{1,-1,-1,-1}}$. But the second and third
of these are conjugated and toral ($\widetilde{\si_{69}}$ relates
  both of them), hence the unique nontoral $Q$ is, up to conjugacy, $Q=
 \span{\si,t_{-1,1,1,1},t_{1,-1,1,1}}$.

In the  case (2) we can take the third element $t_2$ also in the
orbit of $t'_{1,1,1,-1}$ and so
$Q=\span{\si,t'_{1,1,1,-1},t'_{-1,-1,-1,1}}$, which produces a
toral grading again. It is obviously not conjugated to the
previous toral one, because the number of order-two elements in a
determined orbit is preserved by conjugation (besides the
 zero homogeneous components have different dimensions).

\begin{pr}\label{aprop1}
The unique  proper subquasitorus of $A(405,\id)$ of order $8$
(up to conjugacy) are:
\begin{itemize}
\item
$\span{\si,t'_{-1,1,1,1},t'_{1,-1,1,1}}$, a nontoral
$\Z_2^3$-grading of type $(0,0,1,0,0,0,7)$.
\item $\span{\si,t'_{-1,1,1,1},t'_{1,1,1,-1}}$ which is toral.
\item $\span{\si,t'_{1,1,1,-1},t'_{-1,-1,-1,1}}$ which is also toral.
\end{itemize}
The last two gradings are not isomorphic.
\end{pr}

To finish this subsection we should now describe the nontoral
subquasitori $Q$ of $A(405,\id)$ isomorphic to $\Z_2^4$. For this
we take any subquasitori $Q'$ of $Q$ of cardinal $8$, that is,
isomorphic to $\Z_2^3$ and apply the previous study to it. Thus we
should study the possible refinements of the quasitori given in
Proposition~\ref{aprop1} which give nontoral gradings. The
techniques already used in the previous paragraph give that any
such quasitori is conjugated to
$\span{\si,t'_{-1,1,1,1},t'_{1,-1,1,1},t'_{1,1,1,-1}}$, which
gives a $\Z_2^4$-grading of type $(1,8,0,0,7)$. The agreement with
the type of the $\Z_2^4$-quasitorus contained in $A(105,\id)$
suggests that they could be conjugated and the corresponding
gradings isomorphic. Indeed, up to conjugacy there is only one
abelian nontoral subgroup of $F_4$  isomorphic to $\Z_2^3$ and
only one to $\Z_2^4$ (see for instance \cite[Prop.\,3.2]{Viruf4}).
For completeness and selfcontainedness   we prove it in our
context.

\begin{pr}
Any nontoral proper subquasitorus  of $ A(405,\id)$ is conjugated
to a subquasitorus of $A(105,\id)$.
\end{pr}

Proof. Notice that for $f_1=t'_{-1,1,-1,1}$, $f_2=t'_{1,-1,-1,1}$
and $f_3=t'_{1,1,1,-1}$, both $\widetilde{\si_{105}}$ and
$\widetilde{\si_{405}}$ commute with the subgroup
$\tor{105}\cap\tor{405}=\span{f_1,f_2,f_3}\cong\Z_2^3$. In fact
both gradings induced by
$\span{f_1,f_2,f_3,\widetilde{\si_{105}}}$ and
$\span{f_1,f_2,f_3,\widetilde{\si_{405}}}$ are nontoral. The key
fact is that
$\span{f_1,f_2,f_3,\widetilde{\si_{105}}\widetilde{\si_{405}}}$ is
toral (since the fixed part is a four-dimensional abelian
subalgebra), so we can apply  Proposition~\ref{pade1} to
$F=\{f_i\}_{i=1}^5$, for
$f_4=\widetilde{\si_{105}}\widetilde{\si_{405}}$ and
$f_5=\widetilde{\si_{105}} $. The group $\span{F}$ is obviously
nontoral since it is isomorphic to $\Z_2^5$. As in the proof of
Proposition~\ref{pade1}, there is $p\in \F4$ such that
$pf_ip^{-1}\in \To$ for $i=1\dots 4$ and $pf_5p^{-1}\in\No$, that
is, there are $j\in\{1,\dots1152\}$   and $t\in\To$ such that
$pf_5p^{-1}=\widetilde{\si_{j}}t$. As
$\Z_2^4\cong\span{pf_ip^{-1}\mid i=1\dots4}\subset\tor{j}$, we
have either $j=405$ or $j=748$. But if we had  this last
possibility, the grading $\span{pf_ip^{-1}\mid i=1\dots5}$ would
be toral. Hence $pf_5p^{-1}=\widetilde{\si_{405}}t$. Moreover $p$
can be taken such that $pf_5p^{-1}=\widetilde{\si_{405}}$, by
Corollary~\ref{coento}, so that $p\span{F}p^{-1}\subset
A(405,\id)$.

Thus, the unique (up to conjugation) nontoral subquasitorus of
$A(105,\id)$ isomorphic to $\Z_2^4$, which is
$\span{f_1,f_2,f_3,f_5}$ according to Proposition~\ref{grandes},
is conjugated by means of $p$ to a subquasitorus of $A(405,\id)$,
and consequently, the same can be said about the $\Z_2^3$-nontoral
quasitorus $\span{f_1,f_2, f_5}$. The proof is finished because
there are only two nontoral proper subquasitori of $A(405,\id)$,
by Proposition~\ref{aprop1} and the paragraph above. $\square$

Summarizing the previous propositions, we have proved the
following theorem, which describes all the nontoral quasitori of
$\F4$.

\begin{te}\label{main2}
Any nontoral subquasitorus of $\F4$ is conjugated to some of the
following: \sitem{I)}
$A(15,\id)=\span{\widetilde{\si_{15}},t'_{\om,\om,1,\om^2},t'_{1,\om,\om,1}}
\cong\Z_3^3$, where $\om$ is a primitive cubic root of $1$.
\sitem{II)} $A(105,\id)=\span{ \{\widetilde{\si_{105}},
t'_{-1,1,-1,1},t'_{1,-1,-1,1}\}\cup\{ t'_{111u}\}_{u\in
F^\times}}\cong F^\times\times\Z_2^3$, and its nontoral (proper)
coarsenings which up to conjugacy are: \sitem{II.1)}
$\span{\widetilde{\si_{105}},t'_{-1,1,-1,1},t'_{1,-1,-1,1}}\cong
\Z_2^3$. \sitem{II.2)}
$\span{\widetilde{\si_{105}},t'_{-1,1,-1,1},t'_{1,-1,-1,1},
t'_{1,1,1,-1}}\cong\Z_2^4$. \sitem{II.3)}
$\span{\widetilde{\si_{105}},t'_{-1,1,-1,1},t'_{1,-1,-1,1},
t'_{1,1,1,\rho}}\cong\Z_2^3\times\Z_m$, where $m>2$ and $\rho\in
F^\times$ of order $m$.
 \sitem{II.4)}
$\span{t'_{-1,1,-1,1},t'_{1,-1,-1,1},\widetilde{\si_{105}}
t'_{1,1,1,\rho}}\cong\Z_2^2\times\Z_{4m}$, where $\rho\in
F^\times$ is of order $4m$. \sitem{II.5)} $\span{
\{\widetilde{\si_{105}}, t'_{-1,1,-1,1}\}\cup\{ t'_{111u}\}_{u\in
F^\times}}\cong F^\times\times\Z_2^2$. \sitem{III)}
$A(405,\id)=\span{\widetilde{\si_{405}},t'_{-1,1,1,1},t'_{1,-1,1,1},
t'_{1,1,-1,1},t'_{1,1,1,-1}}\cong\Z_2^5$. Its nontoral coarsenings
are conjugated to the quasitori in II.1) and II.2).
\end{te}

Therefore, the following table gives all the nontoral gradings on
$\f4$ up to equivalence. In it, we have taken into consideration
that no repeated gradings arise  produced by the infinite families
of quasitori in some of the cases of the theorem (computed after
Propositions~\ref{grandes}, \ref{pequenos} and \ref{ootra}). We
give the quasitori, the universal grading groups and the types:

\begin{center}
\begin{tabular}{|c|c|c|c|}
\hline Grading & Group & Type& \cr \hline I & $\Z_3^3$ & $( 0,26)$
& Fine\cr II & $\Z_2^3\times\Z$ & $( 31,0,7)$ &Fine\cr II.1 &
$\Z_2^3$ & $( 0,0,1,0,0,0,7)$ & \cr II.2 & $\Z_2^4$ & $(
1,8,0,0,7)$ & \cr II.3.1 & $\Z_2^3\times\Z_3$ & $( 3,14,7)$ & \cr
II.3.2 & $\Z_2^3\times\Z_4$ & $( 17,7,7)$ & \cr II.4.1 &
$\Z_2^2\times\Z_4$ & $( 0,8,2,0,6)$ &\cr II.4.2 &
$\Z_2^2\times\Z_8$ & $( 19,6,7)$ &\cr III & $\Z_2^5$ & $(
24,0,0,7)$ & Fine\cr \hline
\end{tabular}
\end{center}\smallskip

\noindent\textbf{Remark 5.} We should notice that there are 9
different equivalence classes of nontoral gradings on $\f4$. But
there are only 8 different ones on the Albert algebra $J$. The
device for translating gradings from $J$ to $\f4$ (and conversely)
has this deficiency. However it works well when applied to fine
gradings, because MAD's are preserved by the adjoint map. This
explains the agreement in the number of fine gradings (up to
equivalence).\smallskip

\noindent\textbf{Remark 6.} Theorem~\ref{main1} states that every
quasitorus in $\F4 $ is   contained in some $A(j,\id)$. Thus every
nontoral grading is induced by a set of automorphisms formed by
one element in the Weyl group together with several elements in
$\To$. This result is also true in nontoral gradings on $\g2$, but
it is false when applied to other simple Lie algebras, for
instance $\frak d_4=\text{o}(8,F)$. Then, an alternative approach
is to find \lq \lq the first steps", that is, the minimal nontoral
quasitori instead of the maximal ones. Independently of the simple
Lie algebra under study, they are contained in some $A(j,\id)$ and
any nontoral grading can be obtained by refining one of the
gradings produced by them. In the case of $\f4$, we have found
three nontoral minimal quasitori, namely, I, II.1 and II.4.1 in
Theorem~\ref{main2}, which provide
  $\Z_3^3$, $\Z_2^3$ and $\Z_2^2\times \Z_4$-gradings respectively.

\section{Gradings on the Albert Algebra revisited}
\def\as#1{\widehat{\si_{#1}}}

In this section we prove Theorem~\ref{notenJ}. We start by
considering the quasitorus $Q$ inducing a nontoral grading on the
Albert algebra $J$.  The key tool here is the previously mentioned
fact that the automorphism group $F_4=\aut(J)$ and the
automorphism group $\F4=\aut(\f4)$ are isomorphic via the map
$\Ad\colon F_4\to\F4 $ such that $\Ad(f)d:=fdf^{-1}$ for any $f\in
F_4$ and $d\in\f4$ (\cite{Jac}). Then
 we can apply Theorem~\ref{main2} to
have an immediate view (up to conjugacy) of the quasitory $Q$.
 So we study the
different possibilities I)-III) provided by Theorem~\ref{main2}
and collect them in the following table
\smallskip

\begin{center}
\begin{tabular}{|c|c|c|c|}
\hline Grading & Group & Type& \cr \hline I & $\Z_3^3$ & $(27)$ &
Fine\cr II & $\Z_2^3\times\Z$ & $(25,1)$ &Fine\cr II.1 & $\Z_2^3$
& $( 0,0,7,0,0,1)$ & \cr II.2 & $\Z_2^4$ & $( 7,8,0,1)$ & \cr
II.3.1 & $\Z_2^3\times\Z_3$ & $(21,3)$ & \cr II.3.2 &
$\Z_2^3\times\Z_4$ & $(23,2)$ & \cr II.4 & $\Z_2^2\times\Z_4$ & $(
0,12,1)$ &\cr III & $\Z_2^5$ & $( 24,0,1)$ & Fine\cr \hline
\end{tabular}
\end{center}
\smallskip

It is remarkable the fact that all the quasitori in case II.4 of
Theorem~\ref{main2} induce the same grading on $J$ up to
equivalence (by Remark\,4). This was not the case in $\f4$ where
these quasitori produced two nonequivalent gradings. This is the
reason why we only have eight equivalence classes of nontoral
gradings on $J$ while in $\f4$ we have detected nine. Now, any
nontoral grading on $J$ must be equivalent to any of these, and
these are nonequivalent. But the gradings (\ref{nt1}),
(\ref{nt2}), (\ref{nt3}), (\ref{nt4}), (\ref{nt5}), (\ref{grad1}),
(\ref{coar}) and (\ref{ztrescubo}) described in
Theorem~\ref{notenJ} are nonequivalent since their types
 are  different:
\smallskip

\vbox{
\begin{tabular}{|c|c|}
\hline Grading & Type\cr \hline (\ref{nt1}) & (25,1)\cr
(\ref{nt2}) & (7,8,0,1)\cr (\ref{nt3}) & (21,3)\cr \hline
\end{tabular}
\vskip -1.75cm\hskip 4cm
\begin{tabular}{|c|c|}
\hline Grading & Type\cr \hline (\ref{nt4}) & (23,2)\cr
(\ref{nt5}) & (24,0,1)\cr (\ref{grad1}) & (0,0,7,0,0,1)\cr \hline
\end{tabular}
\vskip -1.75cm\hskip 8.5cm
\begin{tabular}{|c|c|}
\hline Grading & Type\cr \hline (\ref{coar}) & (0,12,1)\cr
(\ref{ztrescubo}) & (27)\cr \hline
\end{tabular}
} \ \
\bigskip

\noindent and so they give a complete system of pairwise
nonequivalent nontoral gradings on the Albert algebra.
\medskip

\noindent{\bf Remark 7.} There is only one equivalence class of
group gradings on $J$ with every nonzero component spanned by an
invertible element, namely the $\Z_3\times\Z_3\times\Z_3$-grading
(19). This condition is equivalent to be a Jordan $\Lambda$-torus,
according to \cite[Remark 9.2.1]{Allison}. Although Jordan tori
have been classified by \cite{chinos}, Theorem~\ref{notenJ}
provides an alternative proof of the uniqueness.

\section{Description of the nontoral gradings on $\f4$}

In this section we would like to give a more detailed description
of the fine gradings on $\f4$. This description is going to be
twofold.  On the one hand by using any software allowing
simultaneous diagonalization we can get the homogeneous spaces of
the grading under consideration in terms of the basis introduced
in \ref{ldlb}. We include this description for its possible use in
applications requiring explicit computations. For instance, the
subject of gradings is closely related to the graded contractions
 \cite{contr}. In the latter, new Lie algebras are obtained by
modifying the commuting relations respecting the grading.

But on the other hand we would like to highlight the fact that the
whole algebraic group stuff used in this work has been needed to
prove that we have captured all the gradings. However this is not
necessary at all to describe these gradings. This can be made in
an independent way and this is why in this section we are giving
natural descriptions of all the fine gradings.   Of course the
term natural is used here in a subjective manner meaning that the
gradings are given with no reference to computer methods neither
Weyl group nor maximal torus. Thus any mathematician could check
the gradings ignoring such tools.

\subsection{$\Z_3^3$-grading on $\f4$}

This grading has some referrings in the literature though mostly
in geometry than in algebra. Geometers have studied elementary
$p$-groups with different purposes other than the study of
gradings. The fine $\Z_3^3$-grading appears for instance in
\cite[THEOREM\,11.13]{Griess} from the viewpoint of compact Lie
groups but \cite[8.1]{Viru}
 gives results showing how to translate the arguments to the algebraic groups
setting. On the other hand this $\Z_3^3$-grading has been studied
from the viewpoint of Jordan groups (see \cite[p.\,127]{enci}).
It was Alekseevskij, in \cite[Table\,1]{Alek}, who classified
Jordan subgroups in the exceptional case.

The description of the gradings in terms of root vectors can be
obtained instantaneously by performing a simultaneous
diagonalization of $\f4$ relative to the set of automorphisms
$\{\widetilde{\si_{15}},t'_{\om,\om,1,\om^2},t'_{1,\om,\om,1}\}$.
Thus we obtain the following fine  $\Z_3^3$-grading of type
$(0,26)$ on $L=\f4$:

\begin{eqnarray}
L_{000}=0,\cr L_{001}=\span{ b_{17} + b_{32} + b_{33}, -b_{22} +
b_{34} + b_{42} },\cr L_{002}=\span{ b_8 + b_9 + b_{41}, b_{10} +
b_{18} + b_{46} },\cr L_{010}=\span{ b_1 + b_2 + b_{29}, b_4 -
b_{23} + b_{48} },\cr L_{011}=\span{ -b_{19} + b_{30} + b_{38},
-b_{20} + b_{31} + b_{39} },\cr L_{012}=\span{ b_3 + b_{13} +
b_{40}, b_{11} - b_{12} + b_{45} },\cr L_{020}=\span{ b_5 + b_{25}
+ b_{26}, -b_{24} - b_{28} + b_{47} },\cr L_{021}=\span{ -b_{21} -
b_{35} + b_{36}, b_{16} - b_{27} + b_{37} },\cr L_{022}=\span{
-b_6 - b_{14} + b_{43}, b_7 - b_{15} + b_{44} },\cr L_{100}=\span{
t_{\alpha _1} + (2 + \om)t_{\alpha _2} + 2t_{\alpha _3}, ((-1 -
\om)t_{\alpha _1}) + (-1 - 2\om)t_{\alpha _2} +
 2t_{\alpha _4}
},\cr L_{101}=\span{ \om b_{17} +\om^2 b_{32} + b_{33}, -\om^2
b_{22} + \om b_{34} + b_{42} },\cr L_{102}=\span{ \om b_8 +\om^2
b_9 + b_{41}, \om^2 b_{10} + \om b_{18} + b_{46} },\cr
L_{110}=\span{ \om b_1 +\om^2 b_2 + b_{29}, \om b_4 -\om^2 b_{23}
+ b_{48} },\cr L_{111}=\span{ -\om^2 b_{19} + \om b_{30} + b_{38},
-\om b_{20} +\om^2 b_{31} + b_{39} },\cr L_{112}=\span{ \om^2 b_3
+ \om b_{13} + b_{40}, \om b_{11} -\om^2 b_{12} + b_{45} },\cr
L_{120}=\span{ \om b_5 +\om^2 b_{25} + b_{26}, -\om b_{24} -\om^2
b_{28} + b_{47} },\cr L_{121}=\span{ -\om b_{21} -\om^2 b_{35} +
b_{36}, \om^2 b_{16} - \om b_{27} + b_{37} },\cr L_{122}=\span{
-\om^2 b_6 - \om b_{14} + b_{43}, \om b_7 -\om^2 b_{15} + b_{44}
},\cr L_{200}=\span{ t_{\alpha _1} + (2 +\om^2) t_{\alpha _2} +
2t_{\alpha _3}, (-1 -\om^2) t_{\alpha _1} + (-1 - 2 \om^2)
t_{\alpha _2} +
 2t_{\alpha _4}
},\cr L_{201}=\span{ \om^2 b_{17} + \om b_{32} + b_{33}, -\om
b_{22} +\om^2 b_{34} + b_{42} },\cr L_{202}=\span{ \om^2 b_8 + \om
b_9 + b_{41}, \om b_{10} +\om^2 b_{18} + b_{46} },\cr
L_{210}=\span{ \om^2 b_1 + \om b_2 + b_{29}, \om^2 b_4 - \om
b_{23} + b_{48} },\cr L_{211}=\span{ -\om b_{19} +\om^2 b_{30} +
b_{38}, -\om^2 b_{20} + \om b_{31} + b_{39} },\cr L_{212}=\span{
\om b_3 +\om^2 b_{13} + b_{40}, \om^2 b_{11} - \om b_{12} + b_{45}
},\cr L_{220}=\span{ \om^2 b_5 + \om b_{25} + b_{26}, -\om^2
b_{24} - \om b_{28} + b_{47} },\cr L_{221}=\span{ -\om^2 b_{21} -
\om b_{35} + b_{36}, \om b_{16} -\om^2 b_{27} + b_{37} },\cr
L_{222}=\span{ -\om b_6 -\om^2 b_{14} + b_{43}, \om^2 b_7 - \om
b_{15} + b_{44} }. \nonumber\end{eqnarray}

The easiest way to visualize this grading intrinsically, that is,
with no reference to a particular basis or computer methods,  may
be   considering Lie algebra models based upon $\Z_3$-gradings.
Perhaps the most natural place where to look at the automorphisms
inducing the grading is the Lie algebra $\e6$. Adams gave a
construction of this algebra from three copies of $\hbox{\got
a}_2$ (\cite[p.\,85]{Adams}). This model has been widely spread
for its nice 3-symmetry. Once the automorphisms have been given in
$\e6$ we could hopefully restrict them to $\f4$.
 Given a three-dimensional $F$-vector
space $X$ in which a nonzero alternate trilinear map $\det\colon
X\times X\times X\to F$ has been fixed, we can identify the
exterior product with the dual space by the map $X\wedge
X\stackrel{\approx}\to X^*$ such that $x\wedge y\mapsto
\det(x,y,-)\in \hom(X,F)$. And in a dual way we can identify
$X^*\wedge X^*$ with $X$ through $\det^*$, the dual map of $\det$.
Consider three tridimensional vector spaces  $X_i$ ($i=1,2,3$),
and define:
$$
\mathcal{L} =\s(X_1)\oplus \s(X_2)\oplus \s(X_3)\oplus\ X_1\otimes
X_2\otimes X_3\ \oplus\ X_1^*\otimes X_2^*\otimes X_3^* ,$$
endowed with a Lie algebra structure with the product
\begin{equation}\label{produ}\begin{array}{ll} \protect [\otimes f_i,\otimes
x_i]&=\sum_{\buildrel{k=1,2,3}\over{i\ne j\ne k}}
 f_i(x_i)f_j(x_j)\big(f_k(-)x_k-\frac13f_k(x_k)\id_{X_k}\big)\\
\protect [\protect \otimes x_i,\otimes y_i]&=\otimes (x_i\wedge y_i)\\
\protect [\protect \otimes f_i,\otimes g_i]&=\otimes (f_i\wedge
g_i)
\end{array}\end{equation}
for any   $x_i,y_i\in X_i$, $f_i,g_i\in X^*_i$, with the wedge
products as above, and where the actions of the Lie subalgebra
$\sum \s(X_i)$ on $ X_1\otimes X_2\otimes X_3$ and $X_1^*\otimes
X_2^*\otimes X_3^*$  are the natural ones (the $i$-th simple ideal
acts on the $i$-th slot). The Lie algebra $\mathcal{L}$ is
isomorphic to $\e6$.

Following Hesselink (\cite{Hesse}), we say that a grading is {\em
special} if and only if its $0$-homogeneous component is zero.
 We can observe that $\e6$ admits a {\em special\,} $\Z_3^3$-grading with
the nonzero homogeneous components of the same dimension
($78/26=3$). Since $\mathcal{L} =\L_{\bar 0}\oplus \L_{\bar
1}\oplus \L_{\bar 2}$ is a $\Z_3$-grading for $\L_{\bar
0}=\s(X_1)\oplus \s(X_2)\oplus \s(X_3)$, $\L_{\bar 1}=  X_1\otimes
X_2\otimes X_3$ and $\L_{\bar 2}=X_1^*\otimes X_2^*\otimes X_3^*$,
take $\phi_1$ the automorphism which induces the grading, that is,
$\phi_1\vert_{\L_{\bar i}}=\omega^i \text{id}_{\L_{\bar i}}$ where
$\omega=e^{\frac{2\pi I}3}$ is a primitive cubic root of the unit.

 In order to provide the other automorphisms,
 take into account the following observation.  If $\rho_i\colon X_i\to X_i$,
$i=1,2,3$, are linear maps preserving $\det\colon X_i^3\to F$
(that is,
$\det(x_i,y_i,z_i)=\det(\rho_i(x_i),\rho_i(y_i),\rho_i(z_i))$, or
equivalently, $\det \rho_i=1$), the linear map  $\rho_1\otimes
\rho_2\otimes \rho_3\colon \L_{\bar 1}\to \L_{\bar 1}$ can be
uniquely extended to an automorphism of $\mathcal{L}$ such that
its restriction to $\s(V_i)\subset \L_{\bar 0}$ is the conjugation
map $g\mapsto \rho_i g \rho_i^{-1}$.

Next fix basis $\{u_0,u_1,u_2\}$ of $X_1$, $\{v_0,v_1,v_2\}$ of
$X_2$, and $\{w_0,w_1,w_2\}$ of $X_3$ with
 $\det(u_0,u_1,u_2)=\det(v_0,v_1,v_2)=\det(w_0,w_1,w_2)=1$.
Consider now $\phi_2$ the unique automorphism of $\e6$ extending
the map
$$
u_i\otimes v_j\otimes w_k\mapsto u_{i+1}\otimes v_{j+1}\otimes
w_{k+1}.
$$
(indices module 3). This is of course an order three semisimple
automorphism commuting with $\phi_1$. Finally let $\phi_3$ be the
unique automorphism of $\e6$ extending the map
$$
u_i\otimes v_j\otimes w_k\mapsto \omega^iu_{i}\otimes
\omega^jv_{j}\otimes \omega^kw_{k}=\omega^{i+j+k}u_i\otimes
v_j\otimes w_k.
$$
This is also semisimple and commutes with $\phi_1$ and $\phi_2$.
Thus the set $\{\phi_i\}_{i=1}^3$ induces a $\Z_3^3$-grading on
$\e6$. A computation of its $0$-homogeneous component will suffice
to prove that this grading is nontoral. A first calculation
reveals that the subalgebra of elements fixed by $\phi_1$ and
$\phi_2$ is the linear span of
$$\begin{array}{rll}
 &\{u_1\otimes u_2^*+u_2\otimes
u_3^*+u_3\otimes u_1^*,&u_1\otimes u_3^*+u_2\otimes
u_1^*+u_3\otimes u_2^*,\\
&\,\,v_1\otimes v_2^*+v_2\otimes v_3^*+v_3\otimes
v_1^*,&v_1\otimes
v_3^*+v_2\otimes v_1^*+v_3\otimes v_2^*,\\
&\,\, w_1\otimes w_2^*+w_2\otimes w_3^*+w_3\otimes
w_1^*,&w_1\otimes w_3^*+w_2\otimes w_1^*+w_3\otimes
w_2^*\},\end{array}
$$
where we have taken $u_i^*=u_{i+1}\wedge u_{i+2}$,
$v_i^*=v_{i+1}\wedge v_{i+2}$ and $w_i^*=w_{i+1}\wedge w_{i+2}$
the dual basis of $X_1^*$, $X_2^*$ and $X_3^*$ respectively. While
the corresponding fixed subalgebra for $\phi_1$ and $\phi_3$ is
the linear span of
$$\begin{array}{rll}
 &\{u_1\otimes u_1^*-u_2\otimes u_2^*,\
&u_1\otimes u_1^*-u_3\otimes u_3^*,\\
&\,\,v_1\otimes v_1^*-v_2\otimes v_2^*,\
&v_1\otimes v_1^*-v_3\otimes v_3^*,\\
&\,\,w_1\otimes w_1^*-w_2\otimes w_2^*,\ &w_1\otimes
w_1^*-w_3\otimes w_3^*\},
\end{array}
$$
again a six-dimensional abelian subalgebra. The intersection of
both Cartan subalgebras is obviously zero and therefore our
grading is special and nontoral.  Similar computations prove that
the rest of the homogeneous components are three-dimensional.
  The grading produced by any of
the automorphisms $\phi_i$ is of type $(24,27)$, since they are in
the same conjugacy class. The grading produced by any of the three
couples $\{\phi_i,\phi_j\}$ has one six-dimensional component and
eight nine-dimensional ones. The grading induced by the three
automophisms together is of type $(0,0,26)$. Now we must go down
to see this grading in $\f4$.

To a certain extend, the nice $3$-symmetry described in $\e6$ is
inherited by $\f4$. Indeed graphically speaking, $\f4$ arises by
folding $\e6$. More precisely, taking $X_2=X_3$ we can consider on
$\e6$ the unique automorphism $\tau\colon\e6\to\e6$ extension of
$u\otimes v\otimes w\mapsto u\otimes w\otimes v$. This is an order
two automorphism commuting with the previous $\phi_i$ for
$i=1,2,3$. The subalgebra of elements fixed by $\tau$ is
\begin{equation}
\s(X_1)\oplus \s(X_2)\oplus\ X_1\otimes \Sym^2(X_2)\ \oplus \
X_1^*\otimes \Sym^2(X_2^*), \nonumber\end{equation} where
$\Sym^nX_i$ denotes the symmetric powers (as in
\cite[p.\,473]{Fulton}). This is a simple Lie algebra of dimension
$52$, hence $\f4$. Furthermore, denoting also by
$\phi_i\colon\f4\to\f4$ the restriction of the corresponding
automorphisms of $\e6$ to $\fix\tau$, the set $\{\phi_i\}_{i=1}^3$
is also a set of commuting semisimple order three automorphisms of
$\f4$ with no fixed points other than $0$. So it induces a special
nontoral $\Z_3^3$-grading on $\f4$ with all its homogeneous
components of the same dimension ($52/26=2$), of type $(0,26)$.

\subsection{$\Z_2^5$-grading on $\f4$}

The group of automorphisms inducing this grading is also an
elementary $p$-group so that it is   described in Griess work
\cite[Th.\,7.3,\,p.\,277]{Griess}. Moreover, the grading is
\emph{pure} (that is, there is some homogeneous component which
contains a Cartan subalgebra), therefore it also appears in
Hesselink paper \cite[Table\,1,\,p.\,146]{Hesse}. As before an
instantaneous computer calculation provides its type,
 $(24,0,0,7)$, as well as the description of its homogeneous components
in  terms of the fixed basis:

\begin{eqnarray}
L_{0,0,0,0,0}=0 ,\cr L_{0,0,0,0,1}=\span{ b_2+b_{26} },\cr
L_{0,0,0,1,0}=\span{ b_{41}-b_{17} },\cr L_{0,0,0,1,1}=\span{
b_{43}-b_{19} },\cr L_{0,0,1,0,0}=\span{ b_{38}-b_{14} },\cr
L_{0,0,1,0,1}=\span{ b_{35}-b_{11} },\cr L_{0,0,1,1,0}=\span{
b_3+b_{27} , b_{34}-b_{10} , b_{39}-b_{15} , b_{48}-b_{24} },\cr
L_{0,0,1,1,1}=\span{ b_{30}-b_6 },\cr L_{0,1,0,0,0}=\span{
b_1+b_{25} },\cr L_{0,1,0,0,1}=\span{ b_5+b_{29} },\cr
L_{0,1,0,1,0}=\span{ b_{28}-b_4 , b_{40}-b_{16} , b_{44}-b_{20} ,
b_{22}+b_{46} },\cr L_{0,1,0,1,1}=\span{ b_{45}-b_{21} },\cr
L_{0,1,1,0,0}=\span{ b_7+b_{31} , b_{37}-b_{13} , b_{42}-b_{18} ,
b_{47}-b_{23} },\cr L_{0,1,1,0,1}=\span{ b_{32}-b_8 },\cr
L_{0,1,1,1,0}=\span{ b_{36}-b_{12} },\cr L_{0,1,1,1,1}=\span{
b_{33}-b_9 },\cr L_{1,0,0,0,0}=\span{ t_{\alpha _1} , t_{\alpha
_2} , t_{\alpha _3} , t_{\alpha _4} },\cr L_{1,0,0,0,1}=\span{
b_{26}-b_2 },\cr L_{1,0,0,1,0}=\span{ b_{17}+b_{41} },\cr
L_{1,0,0,1,1}=\span{ b_{19}+b_{43} },\cr L_{1,0,1,0,0}=\span{
b_{14}+b_{38} },\cr L_{1,0,1,0,1}=\span{ b_{11}+b_{35} },\cr
L_{1,0,1,1,0}=\span{ b_{27}-b_3 , b_{10}+b_{34} , b_{15}+b_{39} ,
b_{24}+b_{48} },\cr L_{1,0,1,1,1}=\span{ b_6+b_{30} },\cr
L_{1,1,0,0,0}=\span{ b_{25}-b_1 },\cr L_{1,1,0,0,1}=\span{
b_{29}-b_5 },\cr L_{1,1,0,1,0}=\span{ b_4+b_{28} , b_{16}+b_{40} ,
b_{20}+b_{44} , b_{46}-b_{22} },\cr L_{1,1,0,1,1}=\span{
b_{21}+b_{45} },\cr L_{1,1,1,0,0}=\span{ b_{31}-b_7 ,
b_{13}+b_{37} , b_{18}+b_{42} , b_{23}+b_{47} },\cr
L_{1,1,1,0,1}=\span{ b_8+b_{32} },\cr L_{1,1,1,1,0}=\span{
b_{12}+b_{36} },\cr L_{1,1,1,1,1}=\span{ b_9+b_{33} }.\cr
 \nonumber\end{eqnarray}
But   any of the standard models of $\f4$ would also do for
describing the grading in a basis-free manner.
 One of the best known models is that of
Schafer (\cite[p.112]{Schafer}) as derivations of the Albert
algebra. Starting from the nontoral $\Z_2^3$-grading on the Cayley
algebra $C=\oplus_{g\in\Z_2^3}C_g$ we obtained the
$\Z_2^5$-grading on the Albert algebra $J$ given by
$$\begin{array}{lll}
J_{e,0,0}=<E_1,E_2,E_3>,\quad&J_{g,0,0}=0\ (g\ne e), &\ \\
J_{g,0,1}=C_g^{(1)},&J_{g,1,1}=C_g^{(2)},&J_{g,1,0}=C_g^{(3)},
\end{array}
$$
with $g,e=(0,0,0)\in\Z_2^3$. Obviously the grading induced in
$L=\der(J)=[R_J,R_J]$ has as homogeneous components
$L_a=\{d\in\der( J)\mid d(J_{b})\subset J_{a+b}\ \forall
b\in\Z_2^5\}$, therefore $L_{e,0,0}=0$ and
$$\begin{array}{l}
L_{g,0,1}=\{[R_{x^{(1)}},R_{E_2-E_3}]\mid x\in C_g\}\\
L_{g,1,1}=\{[R_{x^{(2)}},R_{E_3-E_1}]\mid x\in C_g\}\\
L_{g,1,0}=\{[R_{x^{(3)}},R_{E_1-E_2}]\mid x\in C_g\}\quad\forall g\in\Z_2^3\\
L_{g,0,0}=\{D_U\mid U\in N_g\}\oplus\{D_{r_x}\mid x\in
(C_0)_g\}\oplus\{D_{l_x}\mid x\in (C_0)_g\} \quad\forall e\ne
g\in\Z_2^3
\end{array}
$$
where
\begin{itemize}
\item
$N_g=\{d\in\der(C)\mid d(C_{b})\subset C_{ g+b}\ \forall
b\in\Z_2^3\}$ are the components of the grading induced on $\g2$,
all of them two-dimensional and Cartan subalgebras, except for
$N_e=0$,
\item $r_x$ and $l_x$ are the right
and left multiplication operators on $C$,
\item if $U\in
\textrm{o}(C,n)$, $D_U\in\der(J)$ is the derivation given by
$$E_i\mapsto0,\,x^{(1)}\mapsto\overline{U(\bar x)}^{(1)},\,
x^{(2)}\mapsto {U'(x)}^{(2)},\,x^{(3)}\mapsto {U''(x)}^{(3)},$$
where $U'$ and $U''$ are the elements in $\textrm{o}(C,n)$ given
by the local triality principle \cite[p.\,88]{Schafer}, that is,
$U(xy)=U'(x)y+xU''(y)$ for all $x,y\in C$. \end{itemize} Then
clearly
   $h_4=7$ ($\dim L_{g,0,0}=2+1+1=4$), also
    $h_1=8\cdot 3=24$ and the grading is of type $(24,0,0,7)$.

Anyway, we think that there is a more intuitive way of looking at
this grading, as well as at the gradings obtained by crossing
gradings on the Cayley algebra $C$ with gradings on $H_3(F)$.
Recall from 3.3 that if   $H\equiv H_3(F)=\{x\in M_3(F)\mid
x=x^t\}$ and $K\equiv K_3(F)=\{x\in M_3(F)\mid x=-x^t\}$, we could
write
$$J=H\oplus K\otimes C_0.$$
 Since $\f4=\der(J)$ there must exist some model of $\f4$ in these terms. In fact we
can see $\f4$ as
$$ L=\der( C)\oplus K\oplus H_0\otimes C_0$$
identifying the Lie algebra $K$ (subalgebra of $M_3(F)^-$) to
$\der( H_3(F))$  in the known Tits unified construction for the
Lie exceptional algebras (for instance, see
\cite[p.\,122]{Schafer}).

 Consider a $G_1$-grading on the Jordan algebra
$H=\oplus_{g\in G_1}H_g$. This grading will come from a grading on
$M_3(F)$ so that the Lie algebra $K$ will also have an induced
grading.
 Take now the
 $\Z_2^3$-grading on the Cayley algebra
 $C=\oplus_{g\in
G_2=\Z_2^3}C_g$ and the induced grading $\der( C)=\oplus_{g\in
G_2}N_g$. All this material induces a
 $G_1\times G_2$-grading on
$J$ and also on $L$ given by
$$\begin{array}{lll}
J_{g_1,e}=H_{g_1},\quad &J_{g_1,g_2}=K_{g_1}\otimes (C_0)_{g_2},&\\
L_{g_1,e}=K_{g_1},&L_{e,g_2}=N_{g_2}\oplus
(H_0)_e\otimes(C_0)_{g_2},&L_{g_1,g_2}=(H_0)_{g_1}\otimes(C_0)_{g_2}.
\end{array}
$$
In the case of the $\Z_2^5$-grading we have $G_1=\Z_2^2$,  with
the gradings on $H$ and $K$   given by
$$\begin{array}{llll}
H_{0,0}=\span{E_1,E_2,E_3}&H_{0,1}=\span{e_{12}+e_{21}}&H_{1,1}=\span{e_{23}+e_{32}}&H_{1,0}=\span{e_{13}+e_{31}}\\
K_{0,0}=0&K_{0,1}=\span{e_{12}-e_{21}}&K_{1,1}=\span{e_{23}-e_{32}}&K_{1,0}=\span{e_{13}-e_{31}}
\end{array}
$$
and $\dim(C_0)_g=1$,  $\dim N_g=2$ for all
$g\in\Z_2^3\setminus\{(0,0,0)\}$. Therefore
$$
\begin{array}{ll}
\dim J_{e,e}=\dim H_e=3\quad &\dim L_{e,e}=0\\
\dim J_{e,g_2}=0          &\dim L_{e,g_2}=\dim
(N_{g_2}+(H_0)_e\otimes (C_0)_{g_2})=4\\
 \dim J_{g_1,e}=\dim
H_{g_1}=1      & \dim L_{g_1,e}=\dim
K_{g_1}=1\\
\dim J_{g_1,g_2}=\dim K_{g_1}\otimes (C_0)_{g_2}=1   &\dim
L_{g_1,g_2}=\dim (H_0)_{g_1}\otimes (C_0)_{g_2}=1
\end{array}
$$
and so the grading is of type $(24,0,1)$ on $J$,  and $(24,0,0,7)$
on $L=\f4$, as we knew from previous sections.

\subsection{$\Z_2^3\times\Z$-grading of $\f4$}

By contrast with previous gradings, as long as we know this one
does not appear in the mathematical literature. Again a simple
computer aided calculation reveals that the fine
 $\Z_2^3\times\Z$-grading of type $(31,0,7)$ is

$$\begin{array}{rr} L_{0,0,0,-2}=0, & L_{0,0,0,-1}=\span{ b_2+b_9 },\\
L_{0,0,0,0}=\span{ t_{\alpha _2}+2 t_{\alpha _3}+t_{\alpha _4} },
& L_{0,0,0,1}=\span{ b_{33}-b_{26} },\\ L_{0,0,0,2}=0 ,&
L_{0,0,1,-2}=\span{ b_{22}+b_{28} },\\ L_{0,0,1,-1}=\span{
b_{35}-b_{21} },& L_{0,0,1,0}=\span{ b_{38}-b_{14} , b_{40}-b_{16}
, b_{44}-b_{20} },\\ L_{0,0,1,1}=\span{ b_{45}-b_{11} },&
L_{0,0,1,2}=\span{ b_{46}-b_4 },\\ L_{0,1,0,-2}=\span{
b_{15}-b_{10} },& L_{0,1,0,-1}=\span{ b_6-b_5 },\\
L_{0,1,0,0}=\span{ b_1+b_{25} , b_3+b_{27} , b_{48}-b_{24} },&
L_{0,1,0,1}=\span{ b_{29}+b_{30} },\\ L_{0,1,0,2}=\span{
b_{39}-b_{34} },& L_{0,1,1,-2}=\span{ b_{23}+b_{31} },\\
L_{0,1,1,-1}=\span{ b_{19}+b_{32} },& L_{0,1,1,0}=\span{
b_{37}-b_{13} , b_{41}-b_{17} , b_{42}-b_{18} },\\
L_{0,1,1,1}=\span{ b_8+b_{43} },& L_{0,1,1,2}=\span{ b_{47}-b_7
},\\ L_{1,0,0,-2}=\span{ b_{12} },& L_{1,0,0,-1}=\span{ b_9-b_2
},\\ L_{1,0,0,0}=\span{ t_{\alpha _2} , \frac{t_{\alpha
_1}}{2}+t_{\alpha _3} , t_{\alpha _4} },& L_{1,0,0,1}=\span{
b_{26}+b_{33} },\\ L_{1,0,0,2}=\span{ b_{36} },&
L_{1,0,1,-2}=\span{ b_{28}-b_{22} },\\ L_{1,0,1,-1}=\span{
b_{21}+b_{35} },& L_{1,0,1,0}=\span{ b_{14}+b_{38} , b_{16}+b_{40}
, b_{20}+b_{44} },\\ L_{1,0,1,1}=\span{ b_{11}+b_{45} },&
L_{1,0,1,2}=\span{ b_4+b_{46} },\\ L_{1,1,0,-2}=\span{
b_{10}+b_{15} },& L_{1,1,0,-1}=\span{ b_5+b_6 },\\
L_{1,1,0,0}=\span{ b_{25}-b_1 , b_{27}-b_3 , b_{24}+b_{48} },&
L_{1,1,0,1}=\span{ b_{30}-b_{29} },\\ L_{1,1,0,2}=\span{
b_{34}+b_{39} },& L_{1,1,1,-2}=\span{ b_{31}-b_{23} },\\
L_{1,1,1,-1}=\span{ b_{32}-b_{19} },& L_{1,1,1,0}=\span{
b_{13}+b_{37} , b_{17}+b_{41} , b_{18}+b_{42} },\\
L_{1,1,1,1}=\span{ b_{43}-b_8 },& L_{1,1,1,2}=\span{ b_7+b_{47} }.
\end{array}$$

But we can detect this grading without reference to explicit
computations using again the model $\f4=L=\der( C)\oplus K\oplus
H_0\otimes C_0$. We can write the $\Z$-grading on $H$ given in
(\ref{gr1}), jointly with the induced one in $K$ by its extension
to $M_3(F)$, in the following equivalent form:
$$
\begin{array}{ll}
H_2=\span{ e_{23}},\quad &K_2=0,\\
H_1=\span{e_{13}+e_{21}},  &K_1=\span{e_{13}-e_{21}},\\
H_e=\span{E_1,E_2+E_3},   &K_e=\span{E_2-E_3},\\
H_{-1}=\span{e_{12}+e_{31}} ,    &K_{-1}=\span{e_{12}-e_{31}},\\
H_{-2}=\span{e_{32}},    &K_{-2}=0.
\end{array}
$$
Thus, by crossing it with the $\Z_2^3$-grading on $C$ we  obtain
$$
\begin{array}{ll}
\dim J_{2,e}=\dim H_2=1\quad &\dim L_{2,e}=0\\
\dim J_{1,e}=\dim H_1=1       &\dim L_{1,e}=\dim K_1=1\\
\dim J_{e,e}=\dim H_e=2 &\dim L_{e,e}=\dim K_e=1\\
\dim J_{2,g}=0&\dim L_{2,g}=\dim H_2\otimes (C_0)_g=1\\
\dim J_{1,g}=\dim K_1\otimes (C_0)_g=1&\dim L_{1,g}=\dim
H_1\otimes (C_0)_g=1\\
\dim J_{e,g}=\dim K_e\otimes (C_0)_g=1&\dim L_{e,g}=\dim N_g\oplus
(H_0)_e\otimes (C_0)_g=3
\end{array}
$$
and the same dimensions of the homogeneous components with
opposite indices. So the grading is of type $(25,1)$ on $J$ and
$(31,0,7)$ on $\f4$.

\subsection{The remaining nontoral gradings}

Any of the detected gradings on $\f4$    can be obtained by
projections given by epimorphisms from the universal grading
groups corresponding to the fine gradings. However, it is worth to
point out that all of them can also be described by using
different models of the algebra. More precisely, consider a
$G$-grading  on a simple Lie algebra $L$ without outer
automorphisms (equivalently, without automorphisms of the Dynkin
diagram). Then $G$ is a product of cyclic groups $G_i$, each of
which produces certain cyclic $G_i$-grading. The zero component of
most of these $G_i$-gradings is a direct sum of   Lie subalgebras
of type either $\s(V)$ or $\mathop{\rm so}(V)$.  As modules for
the zero component, the remaining components are isomorphic to
either tensor of natural modules or spin ones, respectively (see
\cite{mio} for more details). In the first case, it is possible to
describe the $G$-grading in a similar way that just illustrated
with the $\Z_3^3$-grading on $\f4$.
 We now give some sketches of how this
basis-free method works for the other gradings found in $\f4$, for
instance the $\Z_2^2\times \Z_4$ and $\Z_2^2\times \Z_8$-gradings.

Let $V$ and $W$ be $F$-vector spaces of dimensions $2$ and $4$
respectively. According to \cite{mio}, $\f4$ can be seen in the
way
$$
\mathcal{L} =\s(V)\oplus \s(W) \oplus\   V\otimes W\ \oplus\
\Sym^2V\otimes\bigwedge^2W\oplus V\otimes\bigwedge^3W
$$
and its product given in a similar way that (\ref{produ}).

One of the advantages of this model is that given linear maps
$\varphi\colon V\to V$ and $\tilde\varphi\colon W\to W$ so that
the first one preserves $\det\colon\bigwedge^2V\to F$ and the
second preserves $\det\colon\bigwedge^4W\to F$, the map
$\varphi\otimes \tilde\varphi\colon V\otimes W\to V\otimes W$ can
be uniquely extended to an automorphism of the algebra
$\mathcal{L}\cong\f4$, in such a way that the restriction of this
automorphism to $\s(V)$ is the conjugation $g\mapsto \varphi g
\varphi^{-1}$, and to $\s(W)$ is the conjugation $g\mapsto
\tilde\varphi g \tilde\varphi^{-1}$. Now, denote by $\mathcal
L_{\bar 0}=\s(V)\oplus \s(W)$, $\mathcal L_{\bar 1}=V\otimes W$,
$\mathcal L_{\bar 2}=\Sym^2V\otimes\bigwedge^2W$ and  $\mathcal
L_{\bar 3}= V\otimes \bigwedge^3W$. As $\mathcal{L} =\mathcal
L_{\bar 0}\oplus \mathcal L_{\bar 1}\oplus \mathcal L_{\bar
2}\oplus \mathcal L_{\bar 3}$ is a $\Z_4$-grading,
  take $\phi_1$ the automorphism which gives the grading,
that is, $\phi_1\vert_{\mathcal L_{\bar i}}=I^i
\text{id}_{\mathcal L_{\bar i}}$ where $I$ is a primitive fourth
root of the unit.
Let us fix $\{u_0,u_1\}$ a basis of $V$ with $\det(u_0,u_1)=1$,
and $\{w_0,w_1,w_2,w_3\}$ a basis of $W$ with
$\det(w_0,w_1,w_2,w_3)=1$. Take $\varphi_2\colon V\to V,
u_0\mapsto Iu_0,u_1\mapsto -Iu_1$ and $\tilde\varphi_2\colon W\to
W,w_0\mapsto w_0,w_1\mapsto w_1,w_2\mapsto -w_2,w_3\mapsto -w_3$.
Define $\phi_2\in\F4 $ as the extension of
$\varphi_2\otimes\tilde\varphi_2$. Take now $\varphi_3\colon V\to
V, u_0\mapsto u_1,u_1\mapsto -u_0$ and $\tilde\varphi_3\colon W\to
W,w_0\mapsto w_2,w_1\mapsto w_3,w_2\mapsto  w_0,w_3\mapsto w_1$.
Define then $\phi_3\in\F4 $ as the extension of
$\varphi_3\otimes\tilde\varphi_3$. Take finally $\xi^8=1$ a
primitive eighth root, and define $\tilde\varphi_4\colon W\to W$
by $\tilde\varphi_4(w_0)=\xi^5 w_2$, $\tilde\varphi_4(w_1)=\xi^3
w_3$, $\tilde\varphi_4(w_2)=\xi^5 w_0$,
$\tilde\varphi_4(w_3)=\xi^3 w_1$. Consider $\phi_4\in\F4 $ the
extension of $\varphi_3\otimes\tilde\varphi_4$. The set
$\{\phi_i\}_{i=1}^4$ is a commutative set of semisimple
automorphisms and one can see with some easy though boring hand
computations that the grading induced by
$\{\phi_1,\phi_2,\phi_3\}$ is a $\Z_2^2\times\Z_4$-grading of type
$( 0,8,2,0,6)$,  the grading produced by
$\{\phi_1,\phi_2,\phi_4\}$ is a $\Z_2^2\times\Z_8$-grading of type
$( 19,6,7)$ and the grading produced by
$\{\phi_1,\phi_2,\phi_3,\phi_4\}$ is a $\Z_2^3\times\Z_8$-grading
of type $( 31,0,7)$, equivalent to the fine $\Z_2^3\times
\Z$-grading.
 We display their homogeneous
components.

The $\Z_2^2\times\Z_4$-grading induced by
$\{\phi_1,\phi_2,\phi_3\}$ is:
$$
\begin{array}{l}
\mathcal{L}_{\bar0}=\left\{
\begin{array}{llr}L_{1,1,1}&
=\hbox{diag}(a,a)_W\qquad &\dim 3\\
L_{1,1,-1}&=\hbox{diag}(1,-1)_V \oplus\hbox{diag}(b,-b)_W
\quad&\dim 5\\
L_{1,-1,1}&=\hbox{antidiag}(1,-1)_V\oplus\hbox{antidiag}(b,b)_W
&\dim 5\\
L_{1,-1,-1}&=\hbox{antidiag}(1,1)_V\oplus\hbox{antidiag}(b,-b)_W
&\dim 5\end{array}
\right. \\
\end{array}
$$
with $ a\in\s(2)$ and $ b\in M_2(F)$,
$$
\begin{array}{l}
\mathcal{L}_{\bar1}=\left\{
\begin{array}{llr}L_{I,I,I}&=\span{u_0\otimes w_0-Iu_1\otimes
w_2,u_0\otimes
w_1-Iu_1\otimes w_3}\quad&\dim 2\\
L_{I,I,-I}&=\span{u_0\otimes w_0+Iu_1\otimes w_2,u_0\otimes
w_1+Iu_1\otimes w_3}&\dim 2\\
L_{I,-I,I}&=\span{u_0\otimes w_2-Iu_1\otimes w_0,u_0\otimes
w_3-Iu_1\otimes w_1}
&\dim 2\\
L_{I,-I,-I}&=\span{u_0\otimes w_2+Iu_1\otimes w_0,u_0\otimes
w_3+Iu_1\otimes
w_1}&\dim 2\end{array} \right.\\
\end{array}$$
$$\begin{array}{l}\small
\mathcal{L}_{\bar2}=\left\{ \begin{array}{lll}
L_{-1,1,1}&=\langle(u_0\cdot u_0+u_1\cdot u_1)\otimes (w_0\wedge
w_3-w_1\wedge w_2),&\\&\ (u_0\cdot u_0-u_1\cdot
u_1)\otimes\span{w_0\wedge w_3+w_1\wedge w_2,w_0\wedge
w_2,w_1\wedge w_3},&\\&\  u_0\cdot u_1\otimes (w_0\wedge
w_1-w_2\wedge w_3)\rangle&\dim 5\\
L_{-1,1,-1}&=\langle(u_0\cdot u_0+u_1\cdot u_1)\otimes (w_0\wedge
w_3+w_1\wedge w_2,w_0\wedge w_2,w_1\wedge w_3),&\\&\ (u_0\cdot
u_0-u_1\cdot u_1)\otimes(w_0\wedge w_3-w_1\wedge w_2),&\\&\
u_0\cdot u_1\otimes (w_0\wedge w_1+w_2\wedge w_3)
\rangle&\dim5\\
L_{-1,-1,1}&=\langle(u_0\cdot u_0+u_1\cdot u_1)\otimes (w_0\wedge
w_1+w_2\wedge w_3),&\\&\ (u_0\cdot u_0-u_1\cdot
u_1)\otimes(w_0\wedge w_1-w_2\wedge w_3),&\\&\ u_0\cdot u_1\otimes
\langle w_0\wedge w_3+w_1\wedge w_2,,w_0\wedge w_2,w_1\wedge
w_3\rangle\rangle&\dim5\\
L_{-1,-1,-1}&=\langle(u_0\cdot u_0+u_1\cdot u_1)\otimes (w_0\wedge
w_1-w_2\wedge w_3),&\\&\ (u_0\cdot u_0-u_1\cdot
u_1)\otimes(w_0\wedge w_1+w_2\wedge w_3),&\\&\ u_0\cdot u_1\otimes
(w_0\wedge w_3-w_1\wedge w_2) \rangle&\dim3
\end{array} \right.
\end{array}$$
 and $\mathcal{L}_{\bar3}$   dual to $\mathcal{L}_{\bar1}$,
hence
 $h_5=6$, $h_3=2$ and $h_6=8$, as we wondered. The notation
$\hbox{antidiag}(x_1,\ldots,x_n)$ stands for the $n\times n$
matrix $(a_{ij})$ with all entries zero except for
$a_{i,n-i+1}=x_i$. The subindices in $  L_{ijk}$ indicate that
$\phi_1$, $\phi_2$, $\phi_3$ act with eigenvalues $i,j,k$
respectively. Notice that, although $\phi_i$ has order $4$ for
$i=1,2,3$, $\phi_1\phi_2$ and $\phi_1\phi_3$ have order 2.

The $\Z_2^2\times\Z_8$-grading induced by
$\{\phi_1,\phi_4,\phi_2\}$ is:
$$
\begin{array}{l}
\mathcal{L}_{\bar0}=\left\{
\begin{array}{llr}L_{1,1,1}&=\span{\hbox{diag}(1,-1,1,-1)}_W\quad
 &\dim 1\\
L_{1,1,-1}&=\span{\begin{pmatrix}0&1\\-1&0\end{pmatrix}}_V\oplus
\span{\hbox{antidiag}\left(
\begin{pmatrix}a&0\\0&b\end{pmatrix},\begin{pmatrix}a&0\\0&b\end{pmatrix}\right) }_W\quad
 &\dim 3\\
L_{1,I,1}&= \span{\hbox{diag}\left(
\begin{pmatrix}0&a\\b&0\end{pmatrix},\begin{pmatrix}0&a\\-b&0\end{pmatrix}\right)}_W&\dim 2\\
L_{1,I,-1}&= \span{\hbox{antidiag}\left(
\begin{pmatrix}0&a\\b&0\end{pmatrix},\begin{pmatrix}0&a\\-b&0\end{pmatrix}\right)}_W&\dim
2\\
L_{1,-1,1}&=\span{ \begin{pmatrix}1&0\\0&-1\end{pmatrix} }_V\oplus
\span{\hbox{diag}\left(
\begin{pmatrix}a&0\\0&b\end{pmatrix},\begin{pmatrix}-a&0\\0&-b\end{pmatrix}\right)}_W\quad
 &\dim 3\\
L_{1,-1,-1}&=\span{ \begin{pmatrix}0&1\\1&0\end{pmatrix} }_V\oplus
\span{\hbox{antidiag}\left(
\begin{pmatrix}a&0\\0&b\end{pmatrix},\begin{pmatrix}-a&0\\0&-b\end{pmatrix}\right)}_W\quad
 &\dim 3\\
L_{1,-I,1}&= \span{\hbox{diag}\left(
\begin{pmatrix}0&a\\b&0\end{pmatrix},\begin{pmatrix}0&-a\\b&0\end{pmatrix}\right)}_W&\dim 2\\
L_{1,-I,-1}&= \span{\hbox{antidiag}\left(
\begin{pmatrix}0&a\\b&0\end{pmatrix},\begin{pmatrix}0&-a\\b&0\end{pmatrix}\right)}_W&\dim 2\end{array}
\right. \\
\end{array}
$$
with $a,b\in F$,
$$
\begin{array}{l}
\mathcal{L}_{\bar1}=\left\{
\begin{array}{llr}L_{I,\xi^7,I}&=\span{u_0\otimes w_0-Iu_1\otimes
w_2}\quad&\dim 1\\
L_{I,\xi^7,-I}&=\span{u_0\otimes w_2-Iu_1\otimes w_0}&\dim 1\\
L_{I,\xi^3,I}&=\span{u_0\otimes w_0+Iu_1\otimes w_2}
&\dim 1\\
L_{I,\xi^3,-I}&=\span{u_0\otimes w_2+Iu_1\otimes w_0}&\dim 1\\
L_{I,\xi^5,I}&=\span{u_0\otimes w_1-Iu_1\otimes
w_3}\quad&\dim 1\\
L_{I,\xi^5,-I}&=\span{u_0\otimes w_3-Iu_1\otimes w_1}&\dim 1\\
L_{I,\xi,I}&=\span{u_0\otimes w_1+Iu_1\otimes w_3}
&\dim 1\\
L_{I,\xi,-I}&=\span{u_0\otimes w_3+Iu_1\otimes w_1}&\dim 1\end{array} \right.\\
\end{array}$$
$$\begin{array}{l}
\mathcal{L}_{\bar2}=\left\{ \begin{array}{llr}
L_{-1,1,1}=&\langle(u_0\cdot u_0+u_1\cdot u_1)\otimes (w_0\wedge
w_3-w_1\wedge w_2),\\
& u_0\cdot u_1\otimes (w_0\wedge w_1-w_2\wedge w_3),\\&\ (u_0\cdot
u_0-u_1\cdot
u_1)\otimes(w_0\wedge w_3+w_1\wedge w_2)\rangle&\dim 3\\
L_{-1,1,-1}=&\langle(u_0\cdot u_0+u_1\cdot u_1)\otimes
(w_0\wedge w_1+w_2\wedge w_3),\\
& u_0\cdot u_1\otimes (w_0\wedge w_3+w_1\wedge w_2),\\&\ (u_0\cdot
u_0-u_1\cdot u_1)\otimes(w_0\wedge w_1-w_2\wedge w_3)
\rangle&\dim3\\
L_{-1,I,1}=&\langle(u_0\cdot u_0+u_1\cdot u_1)\otimes (w_1\wedge
w_3),\\
&(u_0\cdot u_0-u_1\cdot u_1)\otimes(w_0\wedge w_2)\rangle&\dim2\\
L_{-1,I,-1}=&\langle  u_0\cdot u_1\otimes w_0\wedge w_2
\rangle&\dim1\\
L_{-1,-1,1}=&\langle(u_0\cdot u_0+u_1\cdot u_1)\otimes
(w_0\wedge w_3+w_1\wedge w_2),\\
& u_0\cdot u_1\otimes (w_0\wedge w_1+w_2\wedge w_3),\\&\ (u_0\cdot
u_0-u_1\cdot
u_1)\otimes(w_0\wedge w_3-w_1\wedge w_2) \rangle&\dim 3\\
L_{-1,-1,-1}=&\langle(u_0\cdot u_0+u_1\cdot u_1)\otimes
(w_0\wedge w_1-w_2\wedge w_3),\\
& u_0\cdot u_1\otimes (w_0\wedge w_3-w_1\wedge w_2),\\&\ (u_0\cdot
u_0-u_1\cdot u_1)\otimes(w_0\wedge w_1+w_2\wedge w_3)
\rangle&\dim3\\
L_{-1,-I,1}=&\langle(u_0\cdot u_0+u_1\cdot u_1)\otimes (w_0\wedge
w_2),\\
& (u_0\cdot u_0-u_1\cdot u_1)\otimes(w_1\wedge w_3)\rangle&\dim2\\
L_{-1,-I,-1}=&\langle u_0\cdot u_1\otimes w_1\wedge w_3
\rangle&\dim1
\end{array} \right.
\end{array}
$$
and $\mathcal{L}_{\bar3}$ dual to $\mathcal{L}_{\bar1}$, hence
$h_3=7$, $h_2=6$ and $h_1=19$, and so this is the grading we are
looking for.

Along this lines all the gradings can be located by using models
of $\f4$.

\section*{Appendix.}

Some  results which play a fundamental roll in our work are
included here. Their references may not be so easily accesible and
so we state the results (without proof) for the seek of
selfcontainedness.  So for instance the well known Borel-Serre
theorem for Lie groups has a version for algebraic groups which is
owed to V. P. Platonov in the following terms:

\begin{te} \label{Bose}(\cite[Theorem~3.15, p.\,92]{Platonov}) A
supersoluble subgroup of semisimple elements of an algebraic group
$G$ is contained in the normalizer of a maximal torus.
\end{te}

Here we must recall that a  group is called supersolvable (or
supersoluble) if it has an invariant normal series whose factors
are all cyclic. Any finitely generated abelian group is
supersolvable.
\newline

Another result that we are applying since the beginning is related
to the number of generators of the quasitorus inducing a grading
on $\f4$.  By abuse of notation we speak of the number of
generators of a quasitorus $Q$ instead of the number of generators
of the related finitely generated abelian group $\X(Q)$. When this
number is $\le 2$ we can say for sure that the grading is toral. A
first approach to this is the fact that every cyclic grading on
$\f4$ is toral. Indeed, a grading is \emph{cyclic} if it is
induced by a diagonalizable automorphism $f$ of $\F4$. This is a
semisimple element and since $\F4$ is a connected algebraic group,
$f$ is in a maximal torus of $\F4$ \cite[Theorem\,11.10,
p.\,151]{Borel}. It is possible to strengthen the previous result
to the case in which the grading group has two generators (stated,
for instance, in \cite[Lemma~1.1.3, p.\,5]{Kasper}).

\begin{lm}\label{Kasper}
Every subquasitorus $Q$ of $\F4$ such that $\X(Q)$ has two
generators is toral.
\end{lm}
Proof. It is known (\cite[Th\,3.5.6,\,p.\,93]{Carter}) that if $G$
is a connected reductive group whose derivated subgroup is simply
connected, then the centralizer of every semisimple element in $G$
is connected. Let $Q=\overline{\span{f_1,f_2}}$ be  a (closed)
abelian subgroup of semisimple elements in $\F4$.  As any
semisimple element belongs to a torus, we replace $f_1$ by
$t_1\in\To$ by conjugation. Now take $Z=\Co_{\F4}(t_1)$, which is
a connected group. Applying \cite[Theorem~1, p.\,94]{g2} for $n=1$
we finish the proof. But without using this fact, we can go on
considering that the element $f_2\in Z$ must be in some maximal
torus of $Z$, say $T$.  But $t_1$ is in the center of $Z$ and
hence in all the maximal torus of $Z$. We have finished since
$\span{t_1,f_2}\subset T$, implying $Q\subset T$.  $\square$
\newline

We finish this appendix describing the action by conjugation of
some elements in $\F4$ on the automorphisms $g_1:=t'_{-1,1,-1,1}$,
$g_2:=t'_{1,-1,-1,1}$, $g_3:=\widetilde{\si_{105}}$ and $g_4
:=t'_{1,1,1,-1}$.
\newline
\begin{center}
\begin{tabular}{|c|c|c|c|}
\hline Element $f$ & $fg_1f^{-1}$ & $fg_2f^{-1}$ & $fg_3f^{-1}$
\cr \hline $\widetilde{\si_{94}}$ & $g_1g_2$ & $g_2g_4$ & $g_3$
\cr $\widetilde{\si_{103}}$ & $g_1g_4$ & $g_2$ & $g_3$ \cr
$\widetilde{\si_{468}}$ & $g_2$ & $g_1$ & $g_3$ \cr
$\widetilde{\si_{485}}$ & $g_2$ & $g_1$ & $g_3g_4$ \cr
$\widetilde{\si_{491}}$ & $g_1g_2$ & $g_1g_4$ & $g_3g_4$ \cr
$\psi$ & $g_1$ & $g_3$ & $g_2$ \cr \hline
\end{tabular}
\end{center}
\smallskip

\noindent Moreover $ft'_{1,1,1,u}f^{-1}=t'_{1,1,1,u}$ for any
$u\in F^\times$ and $f$ in the first column on the left of the
table. We have employed this information to find the nontoral
subquasitori of $A(105,\id)$.

\section*{Thanks}
We would like to thank to Professors Alberto Elduque Palomo for
his suggestions on models of exceptional Lie algebras, and to
Antonio Viruel Arb\'aizar for his help with the geometric
references and also with some results on Lie groups whose
transcription to our setting has been relevant for the work.

\end{document}